\documentclass[reqno,12pt]{amsart}
\usepackage{amsmath, amssymb, amsthm,amsfonts} 
\usepackage[english]{babel}
\usepackage{bbm}
\usepackage{graphicx}
\usepackage{url}
\usepackage{epstopdf}
\usepackage[a4paper,bindingoffset=0.5cm,left=2cm,right=2cm,top=2.5cm,bottom=2cm,footskip=.8cm]{geometry}

\usepackage{tikz}
\usetikzlibrary{arrows, automata,positioning,calc,decorations.pathreplacing,decorations.markings,shapes.misc,petri,topaths}
 \usepackage{pgfplots}
 \pgfplotsset{compat=newest}
 \usetikzlibrary{plotmarks}
 \usepackage{grffile}
\newlength\figureheight
 \newlength\figurewidth
\pgfplotsset{%
  tick label style={font=\scriptsize},
  label style={font=\footnotesize},
  legend style={font=\footnotesize},
     every axis plot/.append style={very thick}
}

\usepackage{rotating}
\usepackage{amsbsy,enumerate}
\usepackage{graphicx}
\usepackage{ccaption}
\usepackage{comment}
\usepackage{mathrsfs} 
\usepackage{diagbox}

\newcommand{\vb}{\vspace{3.2mm}}
\renewcommand{\hat}{\widehat}

\DeclareMathOperator*{\argmin}{arg\,min}

\newcommand{\diag}{\mathrm{diag}}
\newcommand{\tr}{\mathrm{tr}}

\newcommand{\ind}{\BFone}

\newcommand{\BFA}{{\boldsymbol A}}
\newcommand{\BFB}{{\boldsymbol B}}
\newcommand{\BFD}{{\boldsymbol D}}
\newcommand{\BFd}{{\boldsymbol d}}
\newcommand{\BFL}{{\boldsymbol L}}
\newcommand{\BFT}{{\boldsymbol T}}
\newcommand{\BFS}{{\boldsymbol S}}
\newcommand{\BFX}{{\boldsymbol X}}
\newcommand{\BFY}{{\boldsymbol Y}}
\newcommand{\BFy}{{\boldsymbol y}}
\newcommand{\BFn}{{\boldsymbol n}}
\newcommand{\BFmu}{{\boldsymbol \mu}}
\newcommand{\BFSigma}{{\boldsymbol \Sigma}}
\newcommand{\BFr}{{\boldsymbol r}}
\newcommand{\BFzero}{{\boldsymbol 0}}
\newcommand{\BFone}{{\boldsymbol 1}}
\newcommand{\BFH}{{\boldsymbol H}}
\newcommand{\BFI}{{\boldsymbol I}}
\newcommand{\BFJ}{{\boldsymbol J}}
\newcommand{\BFtheta}{{\boldsymbol \theta}}
\newcommand{\BFe}{{\boldsymbol e}}
\newcommand{\BFE}{{\boldsymbol E}}
\newcommand{\BFM}{{\boldsymbol M}}
\newcommand{\BFP}{{\boldsymbol P}}
\newcommand{\BFOmega}{{\boldsymbol \Omega}}
\newcommand{\BFomega}{{\boldsymbol \omega}}
\newcommand{\BFDelta}{{\boldsymbol \Delta}}
\newcommand{\BFdelta}{{\boldsymbol \delta}}
\newcommand{\BFh}{{\boldsymbol h}}
\newcommand{\BFell}{{\boldsymbol \ell}}
\newcommand{\vertiii}[1]{{\left\vert\kern-0.25ex\left\vert\kern-0.25ex\left\vert #1 
		\right\vert\kern-0.25ex\right\vert\kern-0.25ex\right\vert}}

\allowdisplaybreaks

\newtheorem{lemma}{Lemma}
\newtheorem{corollary}{Corollary}

\newtheorem{theorem}{Theorem}

\newtheorem{proposition}{Proposition}

\makeatletter
\renewcommand{\fnum@figure}[1]{\textbf{\figurename~\thefigure}. }
\renewcommand{\fnum@table}[1]{\textbf{\tablename~\thetable}. }
\makeatother

\begin{document}

\title[Road traffic estimation and distribution-based route selection]{Road traffic estimation and \\distribution-based route selection}

\author{Rens Kamphuis, Michel Mandjes, and Paulo Serra}

\begin{abstract}
In route selection problems, the driver's personal preferences will determine whether she prefers a route with a travel time that has a relatively low mean and high variance over one that has relatively high mean and low variance. In practice, however, such risk aversion issues are often ignored, in that a route is selected based on a single-criterion Dijkstra-type algorithm. In addition, the routing decision typically does not take into account the uncertainty in the estimates of the travel time's mean and variance. This paper aims at resolving both issues by setting up a framework for travel time estimation. 

\noindent
In our framework, the underlying road network is represented as a graph.
Each edge is subdivided into multiple smaller pieces, so as to naturally model the statistical similarity between road pieces that are spatially nearby. Relying on a Bayesian approach, we construct an estimator for the joint per-edge travel time distribution, thus also providing us with an uncertainty quantification of our estimates.  Our machinery relies on establishing limit theorems, making the resulting estimation procedure robust in the sense that it effectively does not assume any distributional properties. We present an extensive set of numerical experiments that demonstrate the validity of the estimation procedure and the use of the distributional estimates in the context of data-driven route selection. 

\vb


\vb

\noindent
{\sc Affiliations.} \emph{Rens Kamphuis} and \emph{Michel Mandjes} are with the Korteweg-de Vries Institute for Mathematics, University of Amsterdam, Science Park 904, 1098 XH Amsterdam, the Netherlands. MM is also with E{\sc urandom}, Eindhoven University of Technology, Eindhoven, the Netherlands, and Amsterdam Business School, Faculty of Economics and Business, University of Amsterdam, Amsterdam, the Netherlands. Their research is partly funded by the NWO Gravitation project N{\sc etworks}, grant number 024.002.003. 

\noindent
\emph{Paulo Serra} is with the
Department of Mathematics,
Vrije Universiteit,
De Boelelaan 1111,
1081 HV Amsterdam,
the Netherlands. 

\vb

\noindent
\emph{Date}: {\today}.

\end{abstract}

\maketitle

\newpage

\section{Introduction}



A central problem drivers in a road network are faced with concerns the choice between multiple possible routes in order to travel from their current location to some desired destination. The analysis of such shortest-path problems on a network has a long tradition in operations research. A typical procedure is to consider the per-edge mean travel times, and to apply a Dijkstra-type \cite{DIJK} algorithm to find the fastest route from origin to destination, i.e., the route that minimizes the expected travel time. A conceptual drawback of this approach, however, is that travel times are inherently stochastic. This means that the route that has the shortest expected travel time could also have a substantial standard deviation -- in fact, there may be a route with a higher expected travel time with virtually no variability. 
In such a situation it is up to the driver to make a choice: depending on her personal preferences (in terms of risk aversion) and the importance of the planned trip, she will choose the best alternative. A convenient framework facilitating making such decision uses the concept of \emph{utility functions} \cite{SEN}; see also {e.g.,} \cite{LOUI,SIVA}. Such a utility function could encompass both mean and standard deviation of the travel time, but in principle any distribution-based quantity. A risk averse driver could for instance pick the route that minimizes the 95\%-quantile of the travel time. 

A second conceptual difficulty concerns the way statistical uncertainty is dealt with. If one would aim at identifying the route that optimizes the utility, expressed in terms of a given distribution-based feature of the travel time, it is implicitly assumed that one knows the underlying distribution \emph{with certainty}. In reality, however, the travel times pertaining to the various routes have to be estimated from historic data, necessarily leaving us with some amount of uncertainty. Ignoring this uncertainty, the objective would be to find the route the optimizes the chosen utility function. In a framework accounting for parameter uncertainty, however, the ambition would be to add an uncertainty quantification to this claim. In this context a meaningful statement could be of the type `The probability is $x$\% that the travel time distribution of route A corresponds to a higher utility than the one of route B.'

The main contribution of this paper lies in the development of a broadly applicable framework for travel time estimation in any road traffic network, which is rich enough to also assess the inherent estimation uncertainty. In addition, the performance of our estimation procedure is quantified through a series of numerical experiments, some of them featuring (data-driven) route selection.
Evidently, to determine the optimal route, one has to have a good description of the current state of the road network in terms of the congestion level. As is commonly done, we will treat the network as an undirected graph, where the vertices denote intersections and where the edges connecting these vertices imply the existence of a road between these intersections. In this case, describing the state of the network amounts to estimating the joint per-edge travel time distribution. It is clear that one should not assume that on an edge the level of congestion is evenly spread. Instead, within edges one expects a strong similarity between the congestion levels of spatially nearby road pieces. Moreover, as drivers typically slow down when approaching an intersection, it is to be expected that the velocities near an intersection will be similar for all roads that cross at this intersection but may otherwise vary along the roads represented by those edges. It is part of our approach to incorporate these basic features into our estimation procedure, so as to obtain more accurate estimates of the travel time distribution. 

There is a vast body of literature on estimation techniques for the travel time distribution. Clearly, mean travel times (of a path in the graph, that is) can be derived directly from the mean travel times of the constituent edges, but a major complication is that this property does not carry over to the full travel time distribution or to higher moments. As a consequence, one cannot straightforwardly use techniques for per-edge travel time distribution estimation, such as those discussed in e.g.\ \cite{HELL,HOFL,ZHEN}, to develop an estimation procedure for the path-level travel time distribution. This issue has been resolved in e.g.\ \cite{JENE, RAME, WEST}, but typically at the expense of imposing relatively firm assumptions on the functional form of the per-edge travel time distributions as well as the underlying correlation structure. In \cite{MA} a generalized Markov chain approach has been proposed that estimates the path-level travel time distribution incorporating correlations in time and space. In \cite{RAHM} a non-parametric method is developed that is particularly suited to scenarios in which the travel time distributions vary over time.
Ideally, one would like to have a method of (i)~relatively low computational complexity, that is (ii)~robust in the sense that it does not rely on heavy distributional assumptions, that is (iii) applicable to graphs of any form, also exploiting evident intrinsic properties (such as the ones discussed in the preceding paragraph), and that (iv)~provides us with an uncertainty quantification of the resulting estimates.

Shortest-path problems have a long history in the operations research and combinatorics literature, with Dijkstra's seminal contribution \cite{DIJK} as an important landmark. Various extensions followed. Without aiming at providing an exhaustive overview, we mention a few important contributions; for an in-depth account see e.g.\ \cite{AARD}.
A notable generalization, due to Bellman and Ford \cite{BELL,FORD}, concerns graphs with negative edge weights (assuming, for obvious reasons, no negative cycle can be reached from the source vertex). The so-called $A^\star$ algorithm aims at reducing the subgraph that must be explored \cite{ASTAR}. In e.g.\ \cite{HALP,ORDA} the focus is on networks in which the edges have a time-dependent length. Variants in which the edge lengths attain random values can be found in for instance \cite{BERT,HASS}. In \cite{ADEL} the focus is on adapting Dijkstra's algorithm to the setting of so-called {\sc and-or} graphs.

We proceed with a more detailed account of our contributions. In our modeling framework we represent the road network as a graph, but in our estimation approach we use a version of this graph that is endowed with a higher resolution, i.e., a graph in which each edge is broken up into multiple smaller pieces, each representing a segment of a road. The idea behind working with this high-resolution graph is that it allows us to naturally model the statistical similarity between road pieces that are spatially nearby. Following a Bayesian approach, we construct an estimator for the joint per-edge travel time distribution, thus also providing us with an uncertainty quantification of our estimates. The framework used relies on establishing various limit theorems, making the estimation procedure robust (in the sense that it only very mildly relies on distributional assumptions). The underlying numerics involve basic computational algorithms, predominantly standard routines stemming from linear algebra. Our proposed estimation procedure thus fulfils the desirable properties (i)--(iv). The paper also includes an extensive set of numerical experiments by which we 
thoroughly validate our approach. In addition, we demonstrate the use of the distributional estimates in the context of data-driven route selection: in a series of examples we determine the optimal route from a set of given potential routes, and illustrate how this route is affected by the choice of the utility function, and hence by the driver's preferences.

The remainder of this paper has been organized as follows. Section \ref{sec:model} introduces our notation and model, and defines what our dataset is. Then, in Section \ref{sec:inference} we detail our inference procedure, subsequently considering the mean, covariance structure, and a smoothing parameter~$\lambda$. As pointed out in Section \ref{sec:parameter_models} assumptions on the mean, covariance, and graph Laplacian need to be imposed to make sure that the procedure of Section \ref{sec:inference} is consistent. Section \ref{sec:validation} discusses an extensive set of numerical experiments that have been set up so as to validate the estimation procedure. Then in Section~\ref{sec:route} it is pointed out how our approach can be applied in the context of route selection. Finally, Section~\ref{sec:disc} includes a discussion and concluding remarks. Technical proofs are collected in an appendix.

\section{Notation, model, observations}
\label{sec:model}

In this section we introduce the road traffic network considered, including the notation that we use throughout this paper. In addition, we provide a model for the data collected from this network.

\subsection{Some notation}
\label{sec:model:notation} In this subsection we introduce the graph representation of our road network, including its high-resolution version. 

\subsubsection{Notation for the traffic network}
\label{sec:model:notation:graphs}

We represent the road network by an undirected graph, consisting of vertices that are connected by edges. This graph is, as usual, denoted by $G=(V, E)$ with $V=\{v_1, \dots, v_p\}$ being the set of $p=|V|\in\mathbb{N}$ vertices and $E=\{e_1, \dots, e_q\}$ the set of $q=|E|\in\mathbb{N}$ edges, where $|\cdot|$ denotes the cardinality of the underlying set.
For obvious reasons, we throughout assume that the graph $G$ is connected.
We order the vertices and edges so that we can also identify each vertex $v_i$ and edge $e_j$ with their indices $i$ and $j$, respectively, with $i\in\{1,\ldots,p\}$ and $j\in\{1,\ldots,q\}.$
If this is convenient we sometimes write 
$V=\{1,\dots,p\}$ so that $E \subseteq \{\{i,j\}\in V^2\}$, but we also use the notation
$E=\{1,\dots, q\}$.
For an edge $e=\{i,j\}\in E$ we write $v_a(e)=\min\{i,j\}$ and $v_o(e)=\max\{i,j\}$ for the corresponding vertices of the edge.

We denote by $\BFA\in\{0,1\}^{p\times p}$ the \emph{adjacency matrix} of the graph $G$, i.e., a $p\times p$ matrix whose entries 
indicate whether the corresponding pair of vertices is adjacent or not. More concretely, for $i,j=1,\ldots,p$,
\begin{equation}
\label{eq:adjacency}
A_{i,j} = \left\{\begin{array}{ll}
1&\mbox{ if $\{i,j\}\in E$,}\\
 0&\mbox{ otherwise.}\end{array}\right.
\end{equation}
We write $\BFd = \{d_1, \dots, d_p\} = \{d_v:v\in V\} \in\mathbb{N}^p$ to represent the degrees of the vertices in $V$, so that
\begin{equation}
\label{eq:degrees}
d_i = \sum_{j=1}^p A_{i,j},
\end{equation}
and define $\BFD = \diag\{\BFd\}$.
The \emph{Laplacian matrix} of the graph $G$ is defined by 
$\BFL = \BFD - \BFA\in\mathbb{Z}^{p\times p}$. The diagonal of this matrix consists of the vertices' degrees, 
and the $(i,j)$-th non-diagonal entry is $-1$ if vertices $i$ and $j$ are adjacent and $0$ otherwise. 
At several occasions we want to emphasize in the notation the dependence of certain objects on the underlying graph in the notation. We then write
$V = V(G)$,
$E = E(G)$,
$\BFA = \BFA(G)$,
$\BFd = \BFd(G)$, 
$\BFD = \BFD(G)$ or $\BFD = \diag\{\BFd(G)\}$, and
$\BFL = \BFL(G) = \BFD(G) - \BFA(G)$. 

For a graph $G$ we define the \emph{line graph} of $G$ as another graph $\bar G=(\bar V, \bar E)$, where 
each vertex in $\bar V$ now corresponds to an edge in $G$, and where 
two vertices in $\bar G$ are connected by an edge if, and only if, the corresponding edges in $G$ are incident.
Quantities relating to the line graph of $G$ are denoted with a bar above the quantities; for instance, in our approach we intensively make use of the Laplacian matrix of the line graph of $G$ which we denote as $\bar{\BFL}$.

The graph $G$ represents a traffic network across which particles, to be thought of as cars, are flowing.
In this network, each particle enters the system at some vertex $v\in V$, follows a path to some $v'\in V$, and then leaves the system.
Each particle takes a certain amount of time to cross each edge $e\in E$ on its path, reflecting the current congestion level of that edge.
As pointed out in the introduction it is our objective to infer the travel time distribution pertaining to a given route. Importantly, we wish to do so without a priori assuming that all edges have the same level of congestion and that particles traverse edges at a constant velocity.
In addition we wish to work in a framework by which we can naturally model the statistical similarity between road pieces that are spatially nearby.
To facilitate these requirements it is convenient to subdivide the edges in the traffic network into smaller pieces.
For this we consider a \emph{higher resolution} version of the traffic network.

%

\subsubsection{The traffic network in higher resolution}
\label{sec:model:notation:resolution}

\begin{figure}[!ht]
\begin{center}
\includegraphics[height=2in]{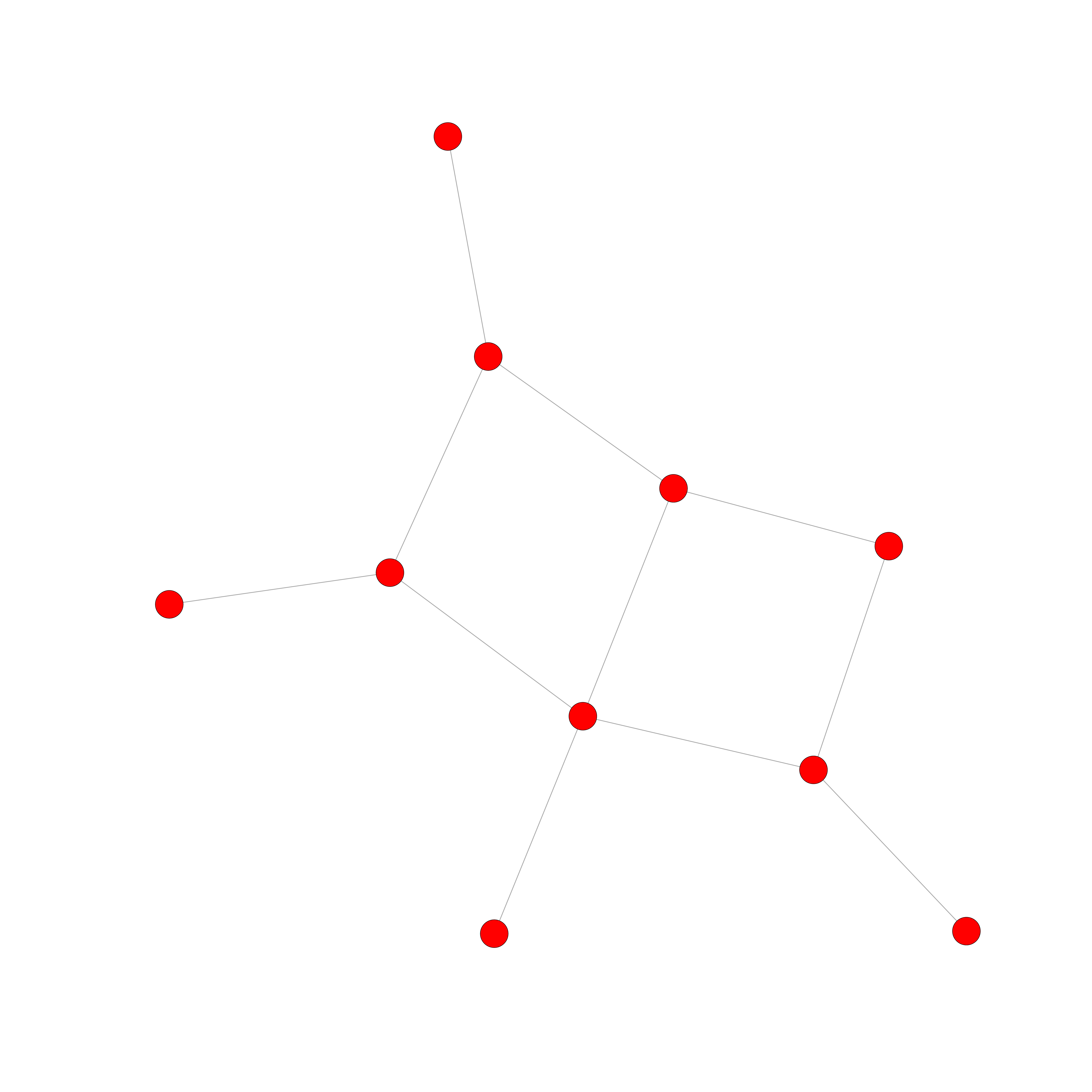}
\includegraphics[height=2in]{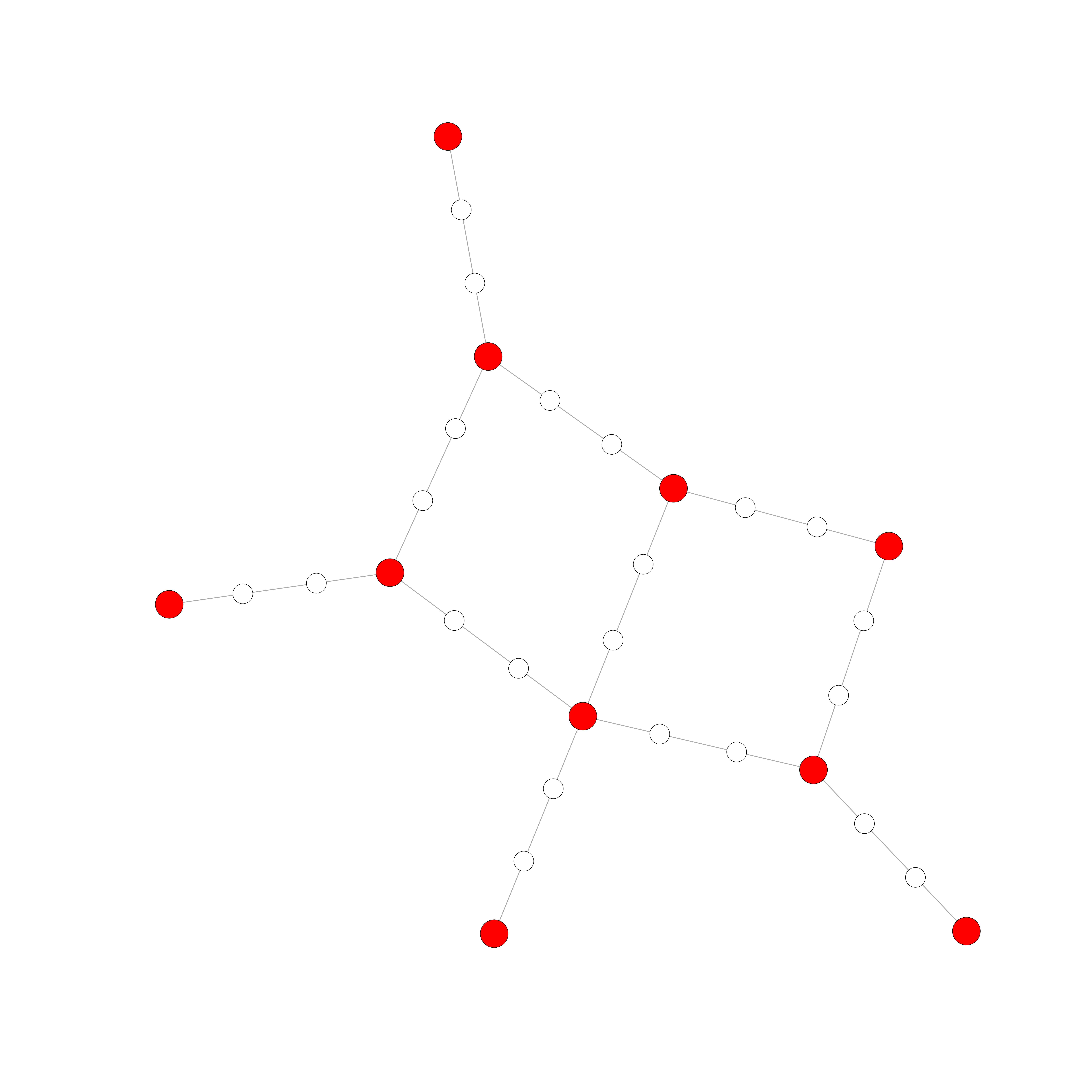}
\caption{Graph $G$ (left) and a higher resolution version $G_\BFr$ of $G$ (right). In this example $r_e=2$, for all $e\in E$.} 
\label{fig:high_res_example}
\end{center}
\end{figure}

We proceed by pointing out how we increase the resolution of the edges. 
To this end, for the graph $G=(V, E)$ we consider a collection $\BFr=\{r_e: e\in E\} =\{r_1, \dots, r_q\}\in\mathbb{N}_0^q$ of user specified \emph{resolution parameters}.
For each such collection $\BFr$, consider the graph $G_\BFr=(V_{\BFr}, E_{\BFr})$ where
\begin{equation}
\label{eq:high_res_sets}
V_{\BFr} :=
V \cup \bigcup_{i : e_i\in E}\big\{v_{i,1}, \dots, v_{i, r_i}\big\}
\qquad\text{and}\qquad
E_{\BFr} := 
\bigcup_{i : e_i\in E} \big\{e_{i,1}, \dots, e_{i,r_i+1}\big\},
\end{equation}
where the edges $e_{i,j}$ in $E_{\BFr}$ are given by, for $j=2, \ldots, r_i$ and $i =1,\ldots, q$,
\begin{equation}
\label{eq:high_res_edges}
e_{i,1}	= \big\{v_a(e_i), v_{i,1}\big\},\;
e_{i,j} 	= \big\{v_{i,j-1}, v_{i,j}\big\},\;
e_{i,r_i+1}	= \big\{v_{i,r_i}, v_o(e_i)\big\}.
\end{equation}
We also denote the vertices of $G_\BFr$ by
$V_{\BFr} = \{v_1, \dots, v_p\} \cup\{v_{i,j}: j=1,\dots, r_i,\, i=1,\dots,p\}$ and its edges by
$E_{\BFr} = \{e_{i,j}: j=1,\dots, r_{i}+1,\, i=1,\dots q\}$.
In the graph $G_\BFr$ we define the cardinalities
\begin{equation}
\label{eq:sizes}
\begin{aligned}
p_\BFr &:=
|V_{\BFr}| = 
|V| + \sum_{e\in E} r_e =
p + \sum_{i=1}^q r_i, \\
q_\BFr &:=
|E_{\BFr}| =
\sum_{e\in E} (r_e+1) =
|E| + \sum_{e\in E} r_e =
q + \sum_{i=1}^q r_i.
\end{aligned}
\end{equation}
We think of the graph $G_\BFr$ as a \emph{higher resolution} version of $G$: 
$G_\BFr$ is constructed from $G$ by replacing each edge $e\in E$ from $G$ by a path graph with $r_e+1$ new edges connecting the original vertices $v_a(e)$ and $v_o(e)$ from $G$. 
We also assume that $\BFr$ is such that the lengths of the road segments corresponding to any edge in $G_\BFr$ is (approximately) the same;
it will become clear in Section~\ref{sec:parameter_models:mean} why we impose this requirement.
Figure~\ref{fig:high_res_example} provides a conceptual illustration of the graph $G$ and its high-resolution version $G_\BFr$. The higher resolution traffic network $G_\BFr$ thus allows us to model the time that a particle takes to traverse each edge $e$ in $G$ in more detail by breaking it down into $r_e+1$ smaller travel times.
The number $r_e+1$ encodes the number of measurements we can collect while the particle moves along $e$, and so we can think of it as a \emph{resolution parameter}.

Finally, in order to translate results for the higher resolution graph into results for the original graph, define a matrix $\BFS_\BFr\in\{0,1\}^{q\times q_\BFr}$ where, for $i = 1, \dots, q$ and $j = 1,\dots,r_{i+1}$,
\begin{equation}
\label{eq:projection}
S_{\BFr,i,j} := \left\{\begin{array}{ll}
1& \mbox{ if $v_{i,j} \in E_{\BFr}$,}\\
0& \mbox{ otherwise.}\end{array}\right.
\end{equation}

\medskip

\subsection{The data format}
\label{sec:model:data}

Consider, for some resolution instance $\BFr\in\mathbb{N}^q$, the corresponding graph $G_\BFr$.
We assume to have access to the average time to traverse each edge in $G_\BFr$.
More explicitly, we assume that we know
\begin{equation}
\label{eq:measurements_def}
X_e^{(n_e)} = 
\frac1{n_e} \sum_{i=1}^{n_e} X_{e,i},
\qquad n_e\in\mathbb{N},\; e\in E_{\BFr},
\end{equation}
where, for each $e\in E_{\BFr}$, $X_{e,i}$ represents the amount of time it took some arbitrary particle to traverse edge $e$, and $n_e$ is the total number of measurements collected at edge $e\in E_{\BFr}$.
Thus, each $X_e^{(n_e)}$ represents the average time it takes particles  to cross the edge $e\in E_{\BFr}$.
In the sequel we abbreviate
\begin{equation}
\label{eq:measurements_abbrev}
\BFn = \{n_e: e\in E_{\BFr}\}
\qquad\text{and}\qquad
\BFX^{(\BFn)} = 
\{X_e^{(n_e)}: e\in E_{\BFr}\}.
\end{equation}

For each $\BFr$ and $\BFn$, our modeling assumption on the corresponding data vector $\BFX^{(\BFn)}$ is that, with ${\mathscr N}({\boldsymbol a},{\boldsymbol b})$ denoting a normally distributed random variable with mean vector ${\boldsymbol a}$ and variance-covariance matrix ${\boldsymbol b}$,
\begin{equation}
\label{eq:model}
\BFX^{(\BFn)} \sim 
{\mathscr N}\Big(\BFmu_\BFr,\; \BFSigma_\BFr^{(\BFn)}\Big).
\end{equation}
Thus, the unknown model parameters are the non-negative vector $\BFmu_\BFr$ and the positive-definite matrix $\BFSigma_\BFr^{(\BFn)}$,
\begin{equation}
\label{eq:model_parameters}
\BFmu_\BFr =
\{\mu_{\BFr, e}: e\in E_{\BFr}\}
\in\mathbb{R}^{q_\BFr},
\qquad 
\BFSigma_\BFr^{(\BFn)} 
\in\mathbb{R}^{q_\BFr\times q_\BFr}.
\end{equation}
We are in the setting that we have ample observations, in that the entries of $\BFn$ are fairly large. In case the dependence between the observations is not excessively strong,
the use of the proposed Gaussian model is justified due to a central-limit type argumentation; we provide more discussion on this issue in Section \ref{sec:disc}. 
Each entry $\mu_{\BFr, e}$ of $\BFmu_\BFr$ represents the expected time for an arbitrary particle to traverse a small segment in the traffic network, corresponding to the edge $e\in E_{\BFr}$. These $\mu_{\BFr, e}$ constitute our main objects of interest, in that we develop a technique to estimate them.
In our setup, the variance-covariance matrix $\BFSigma_\BFr^{(\BFn)}$ plays an important role as well, and is also estimated from the data.

Note that what we would like to infer is actually the travel times in the original graph $G$ (rather than those in the high-resolution graph $G_\BFr$). This means that we want to find
$\BFmu = \BFS_\BFr\,\BFmu_\BFr = \{\mu_e: e\in E\}$, where $\mu_e$ is the expected time for an arbitrary particle to traverse edge $e\in E$.
In the next sections we explain in great detail why it is convenient to collect data at a higher resolution.
We also outline a procedure to infer the model parameters in~\eqref{eq:model_parameters} using a Bayesian approach.

\section{Inference on the model parameters}
\label{sec:inference}

In this section we develop a method to infer the model parameters in~\eqref{eq:model_parameters} from the observations in~\eqref{eq:measurements_def} using a Bayesian approach.
We do so by putting an appropriate prior on $\BFmu_\BFr$ conditional on $\BFSigma_\BFr^{(\BFn)}$, and then estimate $\BFSigma_\BFr^{(\BFn)}$ from the corresponding marginal likelihood for $\BFSigma_\BFr^{(\BFn)}$ using the \emph{empirical Bayes} approach.
In Sections~\ref{sec:inference:estimator}--\ref{sec:inference:smoothness_pars} we give a detailed, step-by-step outline of the estimation procedure, which is then summarized in Section~\ref{sec:inference:final_mean}. The assessment of the performance of this estimation procedure requires some delicate analysis; this is done in Theorem~\ref{theo:moments}.
Numerical illustrations of features of the estimation procedure can be found in Section~\ref{sec:validation}.

\subsection{Estimation of the expected travel times}
\label{sec:inference:estimator}

The objective of this subsection is to propose an estimator for the mean travel times (the vector $\BFmu_\BFr$, that is), and provide an appealing interpretation of it. 

\subsubsection{A Bayesian estimator}
\label{sec:inference:estimator:definition}
To estimate $\BFmu_\BFr$ we follow a so-called \emph{frequentist Bayes} approach, in that 
we assume that the data $\BFX^{(\BFn)}\,|\,(\BFmu_\BFr, \BFSigma_\BFr^{(\BFn)})$ comes from a Bayesian model. 
This means that we endow $\BFmu_\BFr$ with a prior distribution, and 
we use the respective posterior (which is the conditional distribution of $\BFmu_\BFr$, given the data) to produce estimates for $\BFmu_\BFr$.
However, we still see~\eqref{eq:model} as the actual data generating mechanism for fixed $\BFmu_\BFr$ and $\BFSigma_\BFr^{(\BFn)}$ when we study the behavior of the resulting estimates.
Concretely,
we endow $\BFmu_\BFr\,|\,\big(\lambda, \BFSigma_\BFr^{(\BFn)}\big)$ with the following improper\footnote{An improper prior is one whose density integrates to infinity. The corresponding posterior is, however, a proper distribution. See~\cite{SEPK} for more details on such priors.} prior:
\begin{equation}
\label{eq:prior}
\BFmu_\BFr\,|\,\big(\lambda, \BFSigma_\BFr^{(\BFn)}\big)\sim
{\mathscr N}\Big(\BFzero,\; \frac1\lambda \bar{\BFL}_\BFr^-\Big), 
\qquad \lambda> 0,
\end{equation}
where $\bar{\BFL}_\BFr = \BFL(\bar G_\BFr)$ is the Laplacian matrix of the line graph of $G_\BFr$, and where ${\boldsymbol M}^-$ denotes a pseudo-inverse of a matrix ${\boldsymbol M}$. 
As for the prior parameter $\lambda> 0$, for now it suffices to mention that it can be considered as a quantity that controls the concentration of the prior, but in the following section we make the role that it plays in the inference procedure more explicit.

In Proposition~\ref{prop:posterior} (see Appendix A) we show that the posterior distribution for $\BFmu_\BFr\,|\,\big(\lambda, \BFSigma_\BFr^{(\BFn)}\big)$ corresponding to the prior in~\eqref{eq:prior} is normal:
\begin{equation}
\label{eq:posterior}
\BFmu_\BFr\,|\,\big(\lambda, \BFSigma_\BFr^{(\BFn)}, \BFX^{(\BFn)}\big)\sim
{\mathscr N}\Big(\hat{\BFmu}_\BFr,\; \big(\{\BFSigma_\BFr^{(\BFn)}\}^{-1} + \lambda \bar{\BFL}_\BFr \big)^{-1}\Big), 
\qquad \lambda> 0,
\end{equation}
where the posterior mean of $\BFmu_\BFr\,|\,\big(\lambda, \BFSigma_\BFr^{(\BFn)}, \BFX^{(\BFn)}\big)$ is
\begin{equation}
\label{eq:mean_estimator}
\hat{\BFmu}_\BFr(\lambda, \BFSigma_\BFr^{(\BFn)}) := 
\BFH(\lambda, \BFSigma_\BFr^{(\BFn)})\, \BFX^{(\BFn)},
\end{equation}
for a so-called \emph{smoother matrix} $\BFH(\lambda,\BFSigma_\BFr^{(\BFn)})$ defined by
\begin{equation}
\label{eq:smoother}
\BFH(\lambda, \BFSigma_\BFr^{(\BFn)}) =
\Big(\{\BFSigma_\BFr^{(\BFn)}\}^{-1} + \lambda \bar{\BFL}_\BFr \Big)^{-1}\{\BFSigma_\BFr^{(\BFn)}\}^{-1} =
\Big(\BFI_{q_\BFr} + \lambda\, \BFSigma_\BFr^{(\BFn)}\bar{\BFL}_\BFr \Big)^{-1}.
\end{equation}
For the moment, $\lambda$ can be thought of as being fixed, and later we estimate it via generalized cross-validation, as will be pointed out in Section~\ref{sec:inference:smoothness_pars}.

Since the distribution of the prior and the distribution of the posterior belong to the same family of distributions (namely multivariate normal distributions), we say that the prior and the posterior are \emph{conjugate}. It is also common to phrase this as saying that the prior is conjugate for the likelihood of $\BFX^{(\BFn)}\,|\,(\BFmu_\BFr, \BFSigma_\BFr^{(\BFn)})$.
Conjugacy is a desirable property, since it leads to closed-form expressions for the estimator of $\BFmu_\BFr$.

\subsubsection{Interpretation of the estimator}
\label{sec:inference:estimator:interpretation}

The estimator $\hat{\BFmu}_\BFr(\lambda, \BFSigma_\BFr^{(\BFn)})$, being the expectation of a Gaussian distribution, maximizes the posterior density corresponding to~\eqref{eq:posterior}.
This posterior distribution, which is the product of the likelihood and the prior, is proportional (as a function of $\BFmu_\BFr$) to
\begin{equation}
  \label{eq:exp}
\exp\left\{
-\frac12 \big(\BFX^{(\BFn)} - \BFmu_\BFr\big)^\top\big\{\BFSigma_\BFr^{(\BFn)}\big\}^{-1}\big(\BFX^{(\BFn)} - \BFmu_\BFr\big)
-\frac\lambda2 \BFmu_\BFr^\top \bar{\BFL}_\BFr \BFmu_\BFr
\right\}.
\end{equation}
As such, we conclude that the $\BFmu_\BFr$ that maximizes \eqref{eq:exp} solves
\[
\min_{\BFmu_\BFr\in\mathbb{R}^{q_\BFr}}\big(\BFX^{(\BFn)} - \BFmu_\BFr\big)^\top\big\{\BFSigma_\BFr^{(\BFn)}\big\}^{-1}\big(\BFX^{(\BFn)} - \BFmu_\BFr\big)
+ \lambda\, \BFmu_\BFr^\top \bar{\BFL}_\BFr \BFmu_\BFr.
\]
The estimator $\hat{\BFmu}_\BFr(\lambda, \BFSigma_\BFr^{(\BFn)})$ thus solves a so-called \emph{penalized weighted least squares criterion}.
The criterion above, being quadratic in $\BFmu_\BFr$, can be easily solved to yield the estimator $\hat{\BFmu}_\BFr(\lambda, \BFSigma_\BFr^{(\BFn)})$ that is defined in~\eqref{eq:mean_estimator}.
The so-called penalty term, which we call $P(\BFmu_\BFr)$ from now on, can be rewritten as
\[
P(\BFmu_\BFr) =
\BFmu_\BFr^\top \bar{\BFL}_\BFr \BFmu_\BFr =
\sum_{i\leftrightsquigarrow j}\big(\mu_{\BFr,i} - \mu_{\BFr,j}\big)^2, 
\qquad \BFmu_\BFr\in\mathbb{R}^{q_\BFr},
\]
where $i\leftrightsquigarrow j$ denotes that vertices $i$ and $j$ in the line graph $\bar G$ are neighbors, which is equivalent to saying that the corresponding edges in the high-resolution graph $G_\BFr$ are incident.
Informally, a vector $\BFmu_\BFr$ for which $P(\BFmu_\BFr)$ is small is a vector that is \emph{smooth} in the sense that entries of $\BFmu_\BFr$ corresponding to incident edges in the network graph $G_\BFr$ are of a similar magnitude.
In Section~\ref{sec:parameter_models:mean} we explicitly define the class of signals $\BFmu_\BFr$ that we consider.

For a fixed value of $\lambda> 0$, it is clear that among two solutions that fit the data equally well, we always prefer the smoothest solution and, reciprocally, for the solutions that are equally smooth, we always prefer the solution that fits the data best.
The parameter $\lambda$ of the estimator is therefore meant to control the tradeoff between fitting the data well and producing a solution that is smooth in the sense specified above.

The penalized optimization above can be seen as the dual to the (primal) constrained minimization problem
\[
\min_{\BFmu_\BFr\in\mathbb{R}^{q_\BFr}: P(\BFmu_\BFr)\leqslant \ell}\big(\BFX^{(\BFn)} - \BFmu_\BFr\big)^\top\big\{\BFSigma_\BFr^{(\BFn)}\big\}^{-1}\big(\BFX^{(\BFn)} - \BFmu_\BFr\big).
\]
It is known that the two problems are equivalent for an appropriate correspondence between $\lambda$ and $\ell$.
This formulation gives another interpretation for $\hat{\BFmu}_\BFr(\lambda, \BFSigma_\BFr^{(\BFn)})$: it is an estimate that optimally fits the data subject to a maximal `{smoothness budget}' $\ell>0$ (depending on $\lambda$.)

\medskip

The prior distribution in~\eqref{eq:prior} can be interpreted in light of the above discussion. The prior density, being proportional to
\[
\exp\left\{
-\frac\lambda2 \sum_{i\leftrightsquigarrow j}\big(\mu_{\BFr,i} - \mu_{\BFr,j}\big)^2
\right\}, 
\qquad \BFmu_\BFr\in\mathbb{R}^{q_\BFr},
\]
assigns more mass to vectors $\BFmu_\BFr$ that are smooth in terms of the topology of the network.
The larger the \emph{smoothing parameter} $\lambda$ is, the stronger this effect is: 
as $\lambda$ grows, the prior density becomes more tightly concentrated around vectors $\BFmu_\BFr$ in $\mathbb{R}^{q_\BFr}$ that are smooth.

\medskip

Note that the posterior in~\eqref{eq:posterior} provides us with more information than just the estimate $\hat{\BFmu}_\BFr$, as it can also be used to quantify the uncertainty in the produced estimate. More concretely, the extent up to which the posterior measure is concentrated around $\hat{\BFmu}_\BFr$ reflects lack of statistical uncertainty.
This useful information is further exploited in Section~\ref{sec:route}, where we consider path selection problems relying on the information contained in the posterior~\eqref{eq:posterior}.

\subsection{Inferring the covariance structure}
\label{sec:inference:covariance}

We proceed by specifying how the variance-covariance matrix $\BFSigma_\BFr^{(\BFn)}$ can be estimated from the data in a convenient way.

\subsubsection{The empirical Bayes approach for the covariance structure}
\label{sec:inference:covariance:empirical_Bayes}

It turns out to be convenient to parametrize the variance-covariance matrix $\BFSigma_\BFr^{(\BFn)}$ in terms of its eigenvalues.
The matrix $\BFSigma_\BFr^{(\BFn)}$, which is real and symmetric, is rewritten as
\begin{equation}
\label{eq:eigen-decomposition}
\BFSigma_\BFr^{(\BFn)} = 
\BFSigma_\BFr^{(\BFn)}(\BFtheta) =
\sum_{i=1}^{q_\BFr} \theta_i\, \BFe_i\BFe_i^\top =
\sum_{i=1}^{q_\BFr} \theta_i\, \BFE_i,
\end{equation}
where $\BFe_i$ is the $i$-th (column) eigenvector of $\BFSigma_\BFr^{(\BFn)}$, and $\BFtheta=(\theta_1, \dots, \theta_{q_\BFr})$ are the eigenvalues of $\BFSigma_\BFr^{(\BFn)}$.
This entails that we can (trivially) re-parametrize the prior in terms of the eigenvalues $\BFtheta$ as
\begin{equation}
\label{eq:prior_parametrised}
\BFmu_\BFr\,|\,\big(\lambda, \BFtheta\big)\sim
{\mathscr N}\Big(\BFzero,\; \frac1\lambda \bar{\BFL}_\BFr^-\Big), 
\quad \lambda> 0,\; \BFtheta\in\mathbb{R}^{q_\BFr}.
\end{equation}
Combined with the model in~\eqref{eq:model}, this leads to the posterior
\begin{equation}
\label{eq:posterior_parametrised}
\BFmu_\BFr\,|\,\big(\lambda, \BFtheta, \BFX^{(\BFn)}\big)\sim
{\mathscr N}\Big(
\hat{\BFmu}_\BFr(\lambda,\BFtheta),\; 
\big(\BFSigma_\BFr^{(\BFn)}(\BFtheta)^{-1} + \lambda \bar{\BFL}_\BFr \big)^{-1}\Big), 
\quad \lambda> 0,\; \BFtheta\in\mathbb{R}^{q_\BFr}.
\end{equation}

In the \emph{full Bayes} approach we endow $\BFtheta\,|\,\lambda$ with a prior, which together with the (conditional) prior $\BFmu_\BFr\,|\,\big(\lambda, \BFtheta\big)$ results in a joint prior $(\BFmu_\BFr, \BFtheta)\,|\,\lambda$.
Unfortunately, no conjugate priors are known for this parametrization of the model.
In this case, one can still make inference from the respective posterior on $(\BFmu_\BFr, \BFtheta)\,|\,\lambda$ via sampling methods such as Markov Chain Monte Carlo (MCMC); see e.g.\ \cite{GAME}.
However, this might not be computationally attractive since it may require that we have to perform this MCMC procedure for a substantial number of values of $\lambda$.
We therefore opt for an alternative approach.

The \emph{empirical Bayes} method works in the following way.
Suppose that one has an estimator of the eigenvalues, say $\hat\BFtheta(\lambda)$, that depends on $\BFX^{(\BFn)}$ and (eventually) on $\lambda$.
We can obtain a so-called \emph{empirical marginal posterior distribution} by plugging the estimator $\hat{\BFtheta}(\lambda)$ into $\BFtheta$ in the posterior distribution for $\BFmu_\BFr\,|\,\big(\lambda, \BFtheta\big)$, so as to obtain
\begin{equation}
\label{eq:empirical_marginal_posterior_parametrised}
\BFmu_\BFr\,|\,\big(\lambda, \BFX^{(\BFn)}\big)\sim
{\mathscr N}\Big(
\hat{\BFmu}_\BFr(\lambda,\hat{\BFtheta}(\lambda)),\; 
\big(\BFSigma_\BFr^{(\BFn)}(\hat{\BFtheta}(\lambda))^{-1} + \lambda \bar{\BFL}_\BFr \big)^{-1}\Big), 
\quad \lambda> 0;
\end{equation}
cf.\ the posterior \eqref{eq:posterior_parametrised}.
Note that \eqref{eq:empirical_marginal_posterior_parametrised} is not the marginal posterior distribution for $\BFmu_\BFr\,|\,\lambda$, but rather a proxy for it.

From this empirical marginal posterior distribution, which is based on $\hat\BFtheta(\lambda)$, we can obtain the estimators 
\[
\hat{\BFSigma}_\BFr^{(\BFn)}(\lambda) = 
\BFSigma_\BFr^{(\BFn)}(\hat{\BFtheta}(\lambda)),
\] for 
$\BFSigma_\BFr^{(\BFn)}$ 
and 
\[
\hat{\BFmu}_\BFr(\lambda) = 
\hat{\BFmu}_\BFr(\lambda,\hat{\BFtheta}(\lambda)),
\] 
for 
$\BFmu_\BFr$;
one could think of $\hat{\BFmu}_\BFr(\lambda)$ as an \emph{empirical marginal posterior mean}.
For an appropriate choice of the estimator $\hat\BFtheta(\lambda)$, one should still expect the empirical posterior in~\eqref{eq:empirical_marginal_posterior_parametrised} to provide good uncertainty quantification for $\BFmu_\BFr$; 
we refer to e.g.~\cite{ROUS} for a general account of uncertainty quantification for Bayesian estimators based on empirical posteriors.

In the empirical Bayes approach we work with a particular estimator of the eigenvalues $\BFtheta$.
We estimate $\BFtheta$ as the maximizer of the \emph{marginal likelihood} for $\BFtheta$.
With 
$p\big(\BFX^{(\BFn)}\,|\,\BFmu_\BFr, \BFtheta\big)$ representing the likelihood of the data and
$p\big(\BFmu_\BFr\,|\, \lambda, \BFtheta \big)$ being the density of the prior on 
$\BFmu_\BFr\,|\,\big(\lambda, \BFtheta\big)$, the marginal likelihood for $\BFtheta$ is 
\begin{align}\nonumber
p\big(\BFX^{(\BFn)}\,|\, \lambda, \BFtheta \big) &=
\mathbb{E}_{\BFmu_\BFr\,|\,(\lambda, \BFtheta)}\; p\big(\BFX^{(\BFn)}\,|\,\BFmu_\BFr, \BFtheta\big) \\\label{eq:marginal_likelihood_definition}&=
\int\cdots\int p\big(\BFX^{(\BFn)}\,|\,\BFmu_\BFr, \BFtheta\big)\;p\big(\BFmu_\BFr\,|\, \lambda, \BFtheta \big)\, {\rm d}\BFmu_\BFr.
\end{align}
(Note that~\eqref{eq:marginal_likelihood_definition} is simply the normalizing constant for the marginal posterior for $\BFmu_\BFr$.)
It is straightforward to check that in our case the marginal likelihood $p\big(\BFX^{(\BFn)}\,|\, \lambda, \BFtheta \big)$ for $\BFtheta$ can be written as
\begin{equation}
\label{eq:marginal_likelihood_specific}
(2\pi)^{-{q_\BFr}/2}\, \Big|\BFSigma_\BFr^{(\BFn)}(\BFtheta) + \frac1\lambda \bar{\BFL}_\BFr^-\Big|^{-1/2} 
\exp\left\{
-\frac12 \big\{\BFX^{(\BFn)}\big\}^\top\Big(\BFSigma_\BFr^{(\BFn)}(\BFtheta) + \frac1\lambda \bar{\BFL}_\BFr^-\Big)^{-1}\BFX^{(\BFn)}
\right\},
\end{equation}
so that the empirical Bayes estimate of $\BFtheta$ minimizes $-2\ln p\big(\BFX^{(\BFn)}\,|\, \lambda, \BFtheta \big)$. As a consequence, it satisfies
\begin{equation}
\label{eq:empirical_bayes_estimator}
\hat{\BFtheta}(\lambda) =
\arg\min_{\BFtheta}\, 
\big\{\BFX^{(\BFn)}\big\}^\top\Big(\BFSigma_\BFr^{(\BFn)}(\BFtheta) + \frac1\lambda \bar{\BFL}_\BFr^-\Big)^{-1}\BFX^{(\BFn)}
 + \ln \Big|\BFSigma_\BFr^{(\BFn)}(\BFtheta) + \frac1\lambda \bar{\BFL}_\BFr^-\Big|.
\end{equation}
The empirical marginal posterior distribution corresponding to plugging in this particular estimator of $\BFtheta$ is called the \emph{empirical Bayes marginal posterior distribution} for $\BFmu_\BFr$; 
likewise, all of the resulting empirical quantities that we mentioned before just take on the extra qualifier `empirical Bayes' instead of `empirical'.

\subsubsection{Computing the empirical Bayes estimate}
\label{sec:inference:covariance:parametrisation}

In Proposition~\ref{prop:covariance_estimate} (see Appendix A) 
we show that the solution $\hat{\BFtheta}(\lambda)$ in~\eqref{eq:empirical_bayes_estimator} satisfies a certain identity.
Note that the proof of that proposition relies only on the fact that $\BFSigma_\BFr^{(\BFn)}$ is linear in $\BFtheta$, and not on the fact that the entries of $\BFtheta$ are the eigenvalues of $\BFSigma_\BFr^{(\BFn)}$.
It means that we can as well conclude that the empirical Bayes estimate $\hat\BFSigma_\BFr^{(\BFn)}$ of $\BFSigma_\BFr^{(\BFn)}$ satisfies the fixed point equation
\begin{equation}
\label{eq:covariance_estimate_alternative}
\hat\BFSigma_\BFr^{(\BFn)} =
\frac{
\big(\BFX^{(\BFn)} - \BFH(\lambda,\hat\BFSigma_\BFr^{(\BFn)})\BFX^{(\BFn)}\big)\big(\BFX^{(\BFn)} - \BFH(\lambda,\hat\BFSigma_\BFr^{(\BFn)})\BFX^{(\BFn)}\big)^\top}{\tr\big(\BFI_{q_\BFr} - \BFH(\lambda,\hat\BFSigma_\BFr^{(\BFn)})^\top\big)}.
\end{equation}
As a pragmatic way to solve this matrix-valued equation, we approximate the solution to~\eqref{eq:covariance_estimate_alternative} by starting with an arbitrary guess $\hat\BFSigma_{\BFr,0}^{(\BFn)}$ (for instance $\hat\BFSigma_{\BFr,0}^{(\BFn)}=\BFI_{q_\BFr}$), and then iterating 
\[
\hat\BFSigma_{\BFr,i}^{(\BFn)} = 
\frac{
\big(\BFX^{(\BFn)} - \BFH(\lambda,\hat\BFSigma_{\BFr,i-1}^{(\BFn)})\BFX^{(\BFn)}\big)\big(\BFX^{(\BFn)} - \BFH(\lambda,\hat\BFSigma_{\BFr,i-1}^{(\BFn)})\BFX^{(\BFn)}\big)^\top}{\tr\big(\BFI_{q_\BFr} - \BFH(\lambda,\hat\BFSigma_{\BFr,i-1}^{(\BFn)})^\top\big)}, 
\qquad i\in\mathbb{N}.
\]
We continue until a given stopping criterium is met. 
If we stop at iteration $N$, then we approximate the solution to~\eqref{eq:covariance_estimate_alternative} by $\hat\BFSigma_{\BFr,N}^{(\BFn)} $.

\medskip

The above (approximate) solution to~\eqref{eq:covariance_estimate_alternative} is typically problematic: 
a completely unspecified variance-covariance matrix is simply a positive-semi-definite matrix and, as such, it requires a substantial dataset to reliably estimate $\BFSigma_\BFr^{(\BFn)}$'s many entries simultaneously. 
To deal with this complication, a natural remedy is to impose some structure on $\BFSigma_\BFr^{(\BFn)}$. A convenient way to do this is by using a (relatively) low-dimensional parametric model for the variance-covariance matrix.
We thus assume that $\BFSigma_\BFr^{(\BFn)}$ belongs to a low-dimensional family of matrices, say $\{\BFSigma_\BFr^{(\BFn)}(\BFtheta): \BFtheta\in\mathbb{R}^{\tau},\; \tau\in\mathbb{N}\}$.
We can then solve~\eqref{eq:covariance_estimate_alternative} subject to $\BFSigma_\BFr^{(\BFn)}$ belonging to this family;
see Section~\ref{sec:parameter_models:covariance} for a discussion of this approach for our model on $\BFSigma_\BFr^{(\BFn)}$.

Once a particular parametric form $\BFSigma_\BFr^{(\BFn)}(\BFtheta)$ and some initial guess $\hat\BFtheta_0$ has been picked, we can evaluate the following iteration scheme until a given stopping criterium is met:
\[
\BFSigma_\BFr^{(\BFn)}(\hat\BFtheta_i) = 
\frac{
\big(\BFX^{(\BFn)} - \BFH(\lambda,\hat\BFtheta_{i-1})\BFX^{(\BFn)}\big)\big(\BFX^{(\BFn)} - \BFH(\lambda,\hat\BFtheta_{i-1})\BFX^{(\BFn)}\big)^\top}{\tr\big(\BFI_{q_\BFr} - \BFH(\lambda,\hat\BFtheta_{i-1})^\top\big)},
\qquad i\in\mathbb{N},
\]
where \[\BFH(\lambda,\BFtheta) := \Big(\BFSigma_\BFr^{(\BFn)}(\BFtheta)^{-1} + \lambda \bar{\BFL}_\BFr \Big)^{-1}\BFSigma_\BFr^{(\BFn)}(\BFtheta)^{-1}.\]
Importantly, for a specific relevant choice of $\BFSigma_\BFr(\BFtheta)$, the underlying fixed point can be solved explicitly, as pointed out in Corollary~\ref{cor:covariance_estimate} (see Appendix A).

\subsection{Estimation of the parameter $\lambda$}
\label{sec:inference:smoothness_pars}

Now that we have specified how the parameters $\BFmu_\BFr$ and $\BFSigma_\BFr^{(\BFn)}$ can be estimated, it remains to estimate the smoothing parameter $\lambda> 0$; observe that $\BFmu_\BFr$ and $\BFSigma_\BFr^{(\BFn)}$ were estimated {\em for a given value of $\lambda$}.

The procedure from the previous section provides us with estimates $\hat\BFtheta(\lambda)$ for each choice of $\lambda$.
Since our estimator of $\BFmu_\BFr$ is linear in the data, the generalized cross-validation criterion \cite{GOL} for $\lambda> 0$ can be written as
\begin{equation}
\label{eq:GCV}
{\rm GCV}(\lambda) = 
 \frac{{q_\BFr}^{-1}\cdot\big(\BFX^{(\BFn)} - \BFH\big(\lambda,\hat\BFtheta(\lambda)\big)\BFX^{(\BFn)}\big)^\top\big(\BFX^{(\BFn)} - \BFH\big(\lambda,\hat\BFtheta(\lambda)\big)\BFX^{(\BFn)}\big)}
{\big({q_\BFr}^{-1}\cdot\tr\big(\BFI_{q_\BFr}-\BFH(\lambda,\hat\BFtheta(\lambda))\big)\big)^2}.
\end{equation}
This criterium is known to provide an unbiased estimator of the risks $\lambda\mapsto\mathbb{E}\|\hat{\BFmu}_\BFr(\lambda)-\BFmu_\BFr\|^2$ of the estimator $\hat{\BFmu}_\BFr(\lambda)=\hat{\BFmu}_\BFr(\lambda,\hat\BFtheta(\lambda))$ and is fully data-driven.
As such, we can use the $\lambda$ that minimizes the criterion~\eqref{eq:GCV} as an estimate for the $\lambda$ that minimizes the risk.


\subsection{Summary of the procedure}
\label{sec:inference:final_mean}

We conclude this section by summarizing the estimation procedure outlined above, combining the elements from the previous three subsections.
\begin{itemize}
\item[$\circ$]
Our model is
\[
\BFX^{(\BFn)} \sim 
{\mathscr N}\Big(\BFmu_\BFr,\; \BFSigma_\BFr^{(\BFn)}(\BFtheta)\Big), 
\qquad \BFmu_\BFr\in\mathbb{R}^{q_\BFr},\; \BFtheta\in\mathbb{R}^{\tau},\; \tau\in\mathbb{N},
\]
where $\{\BFSigma_\BFr^{(\BFn)}(\BFtheta): \BFtheta\in\mathbb{R}^{\tau},\; \tau\in\mathbb{N}\}$ is some user-specified model for the covariance structure of the data.
\item[$\circ$]
We endow $\BFmu_\BFr$ with the prior~\eqref{eq:prior} with density $p(\BFmu_\BFr\,|\, \lambda, \BFtheta)$ leading to the marginal posterior
\[
\BFmu_\BFr\,|\,\big(\lambda, \BFtheta, \BFX^{(\BFn)}\big)\sim
{\mathscr N}\Big(
\hat{\BFmu}_\BFr(\lambda,\BFtheta),\; 
\big(\BFSigma_\BFr^{(\BFn)}(\BFtheta)^{-1} + \lambda \bar{\BFL}_\BFr \big)^{-1}\Big), 
\quad \lambda> 0,\; \BFtheta\in\mathbb{R}^{q_\BFr}.
\]
\item[$\circ$]
The parameter $\lambda$ can be estimated using generalized cross-validation as a minimizer of \eqref{eq:GCV}.
For each $\lambda$, to obtain $\hat\BFtheta(\lambda)$ we start from some $\hat\BFtheta_0(\lambda)$ and iterate
\[
\BFSigma_\BFr^{(\BFn)}(\hat\BFtheta_i(\lambda)) = 
\frac{
\big(\BFX^{(\BFn)} - \BFH(\lambda,\hat\BFtheta_{i-1}(\lambda))\BFX^{(\BFn)}\big)\big(\BFX^{(\BFn)} - \BFH(\lambda,\hat\BFtheta_{i-1}(\lambda))\BFX^{(\BFn)}\big)^\top}{\tr\big(\BFI_{q_\BFr} - \BFH(\lambda,\hat\BFtheta_{i-1}(\lambda))^\top\big)},
\]
until convergence, where $\BFH(\lambda,\BFtheta) := \Big(\BFSigma_\BFr^{(\BFn)}(\BFtheta)^{-1} + \lambda \bar{\BFL}_\BFr \Big)^{-1}\BFSigma_\BFr^{(\BFn)}(\BFtheta)^{-1}$. For specific choices of the model for $\BFSigma_\BFr^{(\BFn)}(\BFtheta)$, the underlying fixed point may admit an explicit solution (see Corollary~\ref{cor:covariance_estimate}).
\item[$\circ$]
At this point the estimates $\hat\lambda$, and $\hat\BFtheta(\hat\lambda)$ can be plugged into the marginal posterior distribution for $\BFmu_\BFr$, so as to obtain the empirical Bayes marginal posterior distribution for $\BFmu_\BFr$:
\[
\BFmu_\BFr\,|\,\big(\hat\lambda, \hat\BFtheta(\hat\lambda), \BFX^{(\BFn)}\big)\sim
{\mathscr N}\big(
\BFH(\hat\lambda, \hat\BFtheta(\hat\lambda))\BFX^{(\BFn)},\; 
\Big(\BFSigma_\BFr^{(\BFn)}(\hat\BFtheta(\hat\lambda))^{-1} + \hat\lambda \bar{\BFL}_\BFr \big)^{-1}\Big).
\]
\item[$\circ$]
Finally, we should translate results for the higher resolution graph into results for the original graph, so as to translate the estimates for $\BFmu_\BFr$ into estimates for $\BFmu$.
Trivially,
the empirical Bayes marginal posterior distribution for $\BFmu$ is given by
\begin{equation}
\label{eq:posterior_mu}
\BFmu\,|\,\big(\hat\lambda, \hat\BFtheta(\hat\lambda), \BFX^{(\BFn)}\big)\sim
{\mathscr N}\Big(
\BFS_\BFr\,\BFH(\hat\lambda, \hat\BFtheta(\hat\lambda))\BFX^{(\BFn)},\; 
\BFS_\BFr\big(\BFSigma_\BFr(\hat\BFtheta(\hat\lambda))^{-1} + \hat\lambda \bar{\BFL}_\BFr \big)^{-1}\BFS_\BFr^\top\Big).
\end{equation}
\item[$\circ$]
Note that based on this posterior distribution we can compute the posterior distribution of the expected travel time for any path in $G$.
\end{itemize}

\section{Assumptions on the mean, covariance, and graph Laplacian}
\label{sec:parameter_models}

Since the number of parameters that we have to estimate grows with the entries of $\BFr$, the model presented in Section~\ref{sec:model} is in a sense too general.
Indeed, the procedures presented in Section~\ref{sec:inference} can only possibly be consistent if we impose some more constraints on our parameters $\BFmu_\BFr$ and $\BFSigma_\BFr^{(\BFn)}$; 
we do this in Sections~\ref{sec:parameter_models:mean}--\ref{sec:parameter_models:covariance}.
In Section \ref{sec:parameter_models:spectra} we specify the spectrum of $\bar{\BFL}_\BFr$ and its relation with that of $\BFSigma_\BFr^{(\BFn)}$.

\subsection{A model for $\BFmu_\BFr$}
\label{sec:parameter_models:mean}

The vectors $\BFmu$ and $\BFmu_\BFr$ represent the expected travel times for each edge in the graphs $G$ and $G_\BFr$, respectively.
Intuitively, it is clear that the expected times necessary to traverse road segments corresponding to two incident edges in $G_\BFr$ are likely to be close since we assume that the lengths of these road segments are (approximately) the same.
We model this by considering only signals $\BFmu_\BFr$ in the family, for some $C>0$, 
\[
\mathscr{M}_\BFr(C) = 
\left\{\BFmu_\BFr: \max_{i=1,\dots,q} r_i\, P(\BFmu_{\BFr,i}) \leqslant C^2\right\},
\]
where $\BFmu_{\BFr,i}$ corresponds to the entries of the vector $\BFmu_\BFr$ associated with the $i$-th edge of $G$.

The interpretation is the following.
We think of splitting the edge $e_i$ into sub-edges $e_{i,j}$, $j=1,\dots,r_i+1$, as corresponding to splitting the interval $[0,1]$ into $r_i+1$ sub-intervals of equal length.
Suppose that associated with each edge $e_i\in E$ we have a function $m_i:[0,1]\mapsto\mathbb{R}_+$, $i=1,\dots,q$.
If we then see each expected travel time $\mu_{i,j}$ as being obtained from the respective function $m_i$, 
then, recognizing a Riemann sum,
\begin{align*}
r_i \,P(\BFmu_{\BFr,i}) &= 
r_i \sum_{j=1}^{r_i}\big(\mu_{i,j+1} - \mu_{i,j}\big)^2 =
r_i \sum_{j=1}^{r_i}\Big\{m_i\Big(\frac{j}{r_i}\Big) - m_i\Big(\frac{j-1}{r_i}\Big)\Big\}^2 \\ & \approx
\frac1{r_i}\sum_{j=1}^{r_i}m_i'\Big(\frac{j}{r_i}\Big)^2 \approx
\int_0^1 \big\{m_i'(t)\big\}^2\,{\rm d}t.
\end{align*}
For any $\BFmu_\BFr\in\mathscr{M}_\BFr(C)$ it is then the case that
\begin{equation}
\label{eq:bound_P}
P(\BFmu_\BFr) \leqslant 2 \sum_{i=1}^q P(\BFmu_{\BFr,i}) \leqslant \frac{2q\, C^2}{\min_{i=1,\dots,q}r_i}.
\end{equation}
Upon combining the above, we conclude that assuming that our signal $\BFmu_\BFr\in\mathscr{M}_\BFr(C)$ amounts to requiring that the underlying travel times for the different edges in $G_\BFr$ are well represented by a smooth function, meaning a differentiable function whose derivative is square integrable.
A similar class of functions has been relied upon in~\cite{KIRI}. 

\subsection{A simple model for $\BFSigma_\BFr^{(\BFn)}$}
\label{sec:parameter_models:covariance}
In the sequel we partition the variance-covariance matrix into blocks according to the entries in $\BFr$.
More concretely, this means that we write
\[
\BFSigma_\BFr^{(\BFn)} = 
\left(\begin{array}{ccc}\BFSigma_{\BFr, 1,1}^{(\BFn)}&\cdots
&\BFSigma_{\BFr, 1,q}^{(\BFn)}\\
\vdots&\ddots&\vdots\\
\BFSigma_{\BFr, q,1}^{(\BFn)}&\cdots
&\BFSigma_{\BFr, q,q}^{(\BFn)}
\end{array}\right),\]
with the $(k,\ell)$-th entry of the $(i,j)$-th block being
\[
\big(\BFSigma_{\BFr, i,j}^{(\BFn)}\big)_{k,\ell} = 
\mathbb{V}\Big(X_{e_{i,k}}^{(n_i)}, X_{e_{j,\ell}}^{(n_j)}\Big),
\quad
k=1,\dots,r_i+1,\;
\ell=1,\dots, r_j+1,
\]
where $\mathbb{V}(X,Y)$ denotes the covariance between $X$ and $Y$.
As mentioned before, the resulting model for $\BFSigma_\BFr^{(\BFn)}$ is too high-dimensional (and so not amenable to inference) without further assumptions on these covariances.

The structure of $\BFSigma_\BFr^{(\BFn)}$ can be simplified by making more assumptions, thus reducing the dimension of the model on the variance-covariance matrix.
Such assumptions can be motivated in different ways.
For example:
a) a priori knowledge, may allow us to make certain (auto-)covariance assumptions for the particles traversing the network;
b) practical considerations, may lead us to work with a model for which the fixed point equation~\eqref{eq:covariance_estimate_alternative} admits an explicit solution;
c) the observation scheme that is used to collect the data may also validate certain independence assumptions.
Below we mention a few concrete examples.

If we assume that measurements collected at (sub-edges of) edge $i$ are uncorrelated with measurements collected at (sub-edges of) edge $j$ (with $i\not=j$), then off-diagonal blocks in $\BFSigma_\BFr^{(\BFn)}$ are zero.
If we assume that the travel time of each particle is being measured at every sub-edge of edge $i$, then we can model the situation of constant
autocorrelations in these travel times by taking each block $\BFSigma_{\BFr,i,i}^{(\BFn)}$ equal to a Toeplitz matrix.
If we assume that measurements collected at different sub-edges of edge~$i$ are uncorrelated (because, for instance, at edge $e_{i,k}$ we measure $n_i$ particles at random among all those that pass that edge), then
\[
\big(\BFSigma_{\BFr, i,i}^{(\BFn)} \big)_{k,\ell}=
\frac1{n_i^2}\sum_{s=1}^{n_i}\sum_{t=1}^{n_i} \mathbb{V}\Big(X_{e_{i,k},s},\, X_{e_{i,\ell},t}\Big),
\quad
k,\ell=1,\dots,r_i+1,
\]
where $X_{e_{i,k},s}$ represents the $s$-th measurement collected at edge $e_{i,k}$.

\medskip

In the sequel, to derive an asymptotic result, we make a few choices ultimately leading to the form~\eqref{eq:our_covariance} for $\BFSigma_\BFr^{(\BFn)}$.
Asymptotic results for other choices can be worked out in a similar manner (but at the expense of more involved computations.)
Nonetheless, the form in~\eqref{eq:our_covariance} should be appropriate in many circumstances.
With $\delta_{i,j}:=\ind\{i=j\}$, we assume
\[
\mathbb{V}\Big(X_{e_{i,k},s}, X_{e_{j,l},t}\Big) =
\frac{\delta_{i,j}\delta_{k,l}\delta_{s,t}}{r_i+1} \sigma_i^2.
\]
This means that all measurements $X_{e,i}$, for $e\in E_\BFr$ and $i=1,\dots,n_e$, are uncorrelated.
The scaling $1/(r_i+1)$ ensures that, with $X_{e_i}$ representing the total time for a particle to traverse all of the sub-edges $X_{e_{i,j}}$ of some edge $e_i \in E$, by assumption,
\[
\mathbb{V}\Big(X_{e_i}\Big) =
\mathbb{V}\Big(\sum_{i=1}^{r_e+1}X_{e,i}\Big) =
\sum_{i=1}^{r_e+1}\mathbb{V}\Big(X_{e,i}\Big) =
\sigma_e^2.
\]
We also assume that $\BFn = (n,\dots,n)^\top $, so that, upon combining the above,
\begin{equation}
\label{eq:our_covariance}
\BFSigma_\BFr^{(\BFn)} =
\frac1n \BFSigma_\BFr,
\qquad\text{for}\qquad
\BFSigma_\BFr := \diag\{\sigma_1^2\BFone_{r_1+1}, \dots, \sigma_q^2\BFone_{r_q+1}\},
\end{equation}
where $\BFone_r:=(1,\dots,1)^\top \in\mathbb{R}^r$ represents an all-ones vector of length $\BFr$.

This model for $\BFSigma_\BFr^{(\BFn)}$ arises in a situation where we assume that the different particles in the network behave independently, where at each edge the travel times being measured correspond to a random subset of $n$ particles traversing that edge, while allowing the travel times for each edge in $E$ to have different variances.

\subsection{The spectra of $\bar{\BFL}_\BFr$ and $\BFSigma_\BFr^{(\BFn)}$}
\label{sec:parameter_models:spectra}

Establishing asymptotic results for the estimators presented in Section~\ref{sec:inference} relies heavily on understanding the spectrum of $\bar{\BFL}_\BFr$ and of $\BFSigma_\BFr^{(\BFn)}$.
In this section we specify these.

\subsubsection{Decomposition of $\bar{\BFL}_\BFr$}
\label{sec:asymptotics:spectra:laplacian}

\begin{figure}[!ht]
\begin{center}
\includegraphics[height=2in]{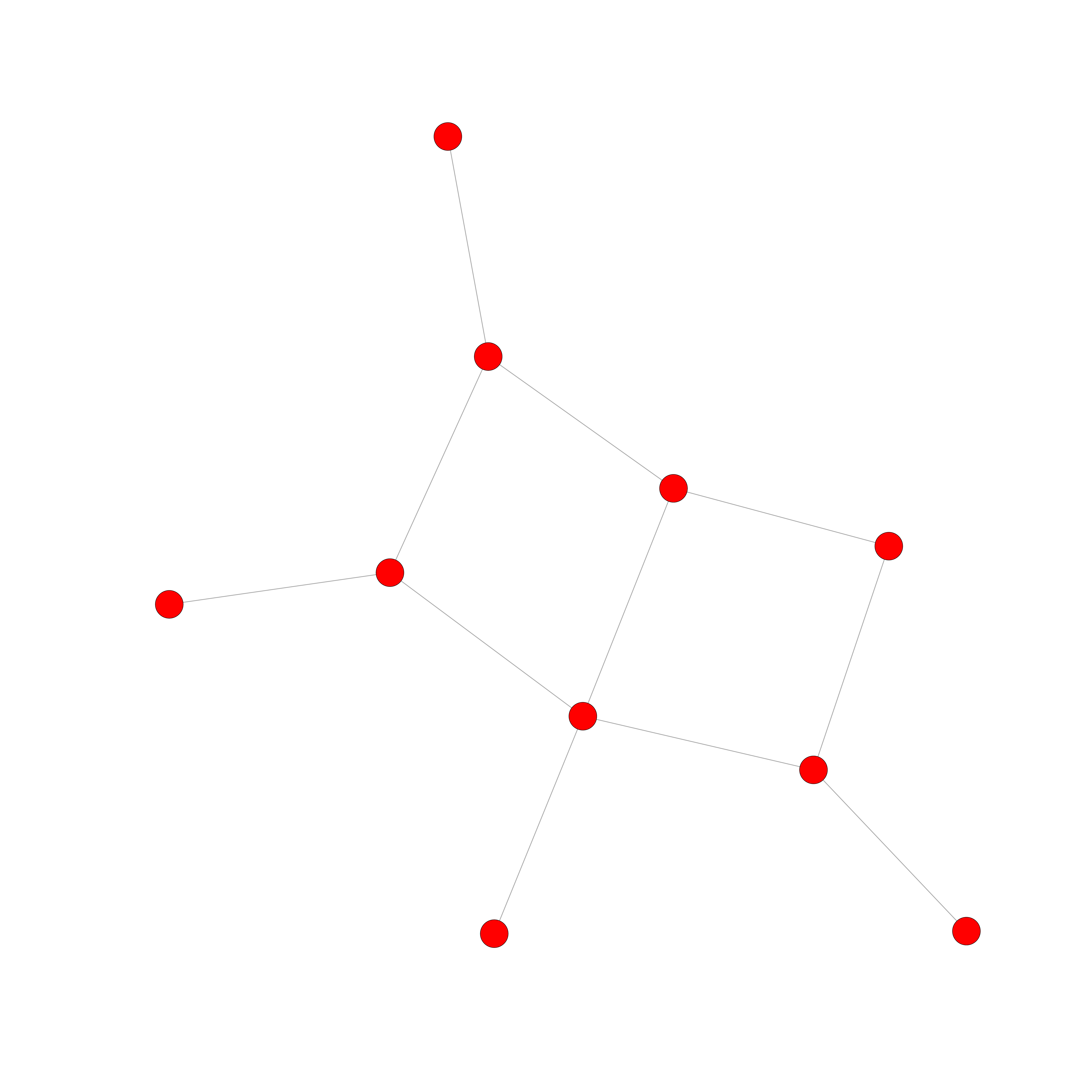}
\includegraphics[height=2in]{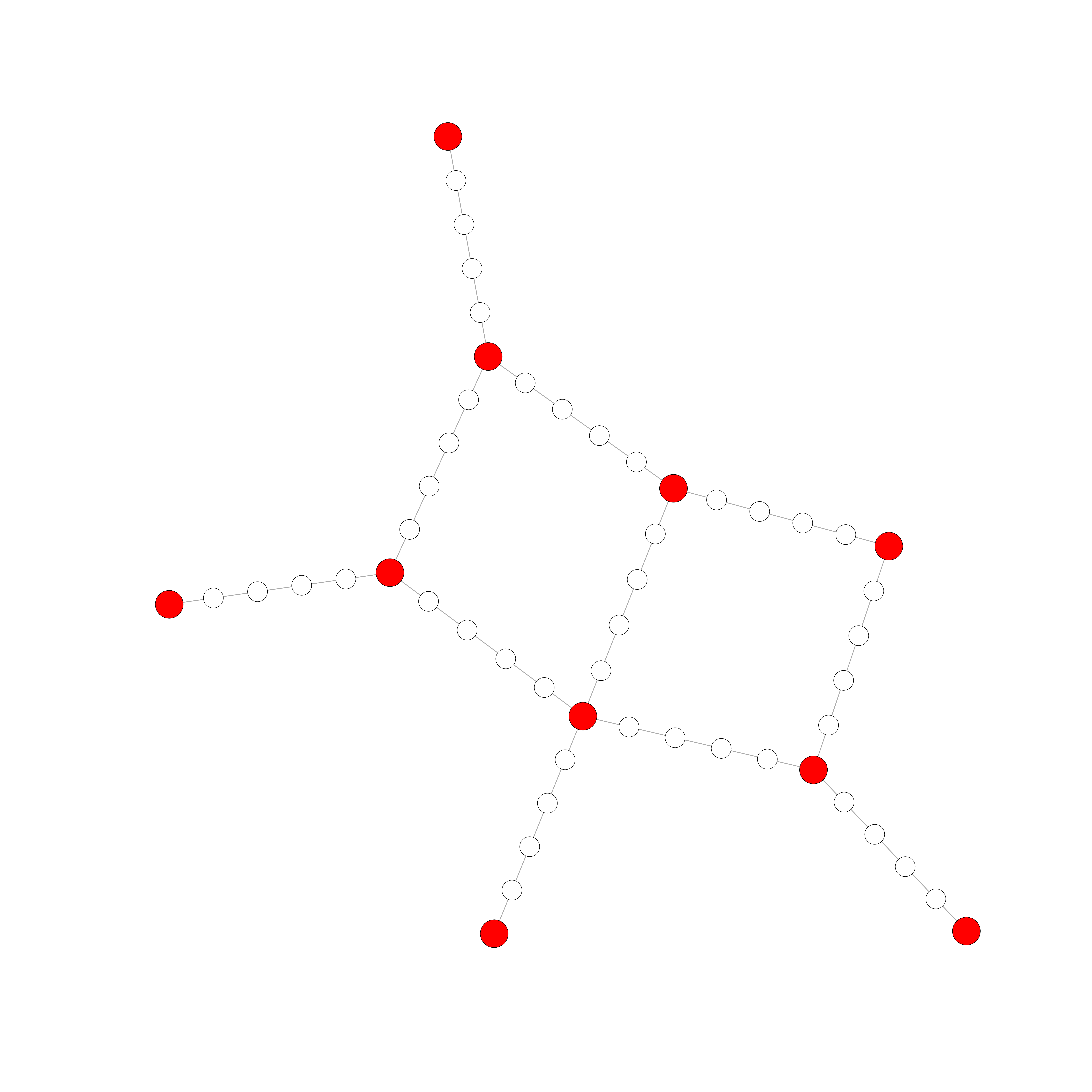}
\includegraphics[height=2in]{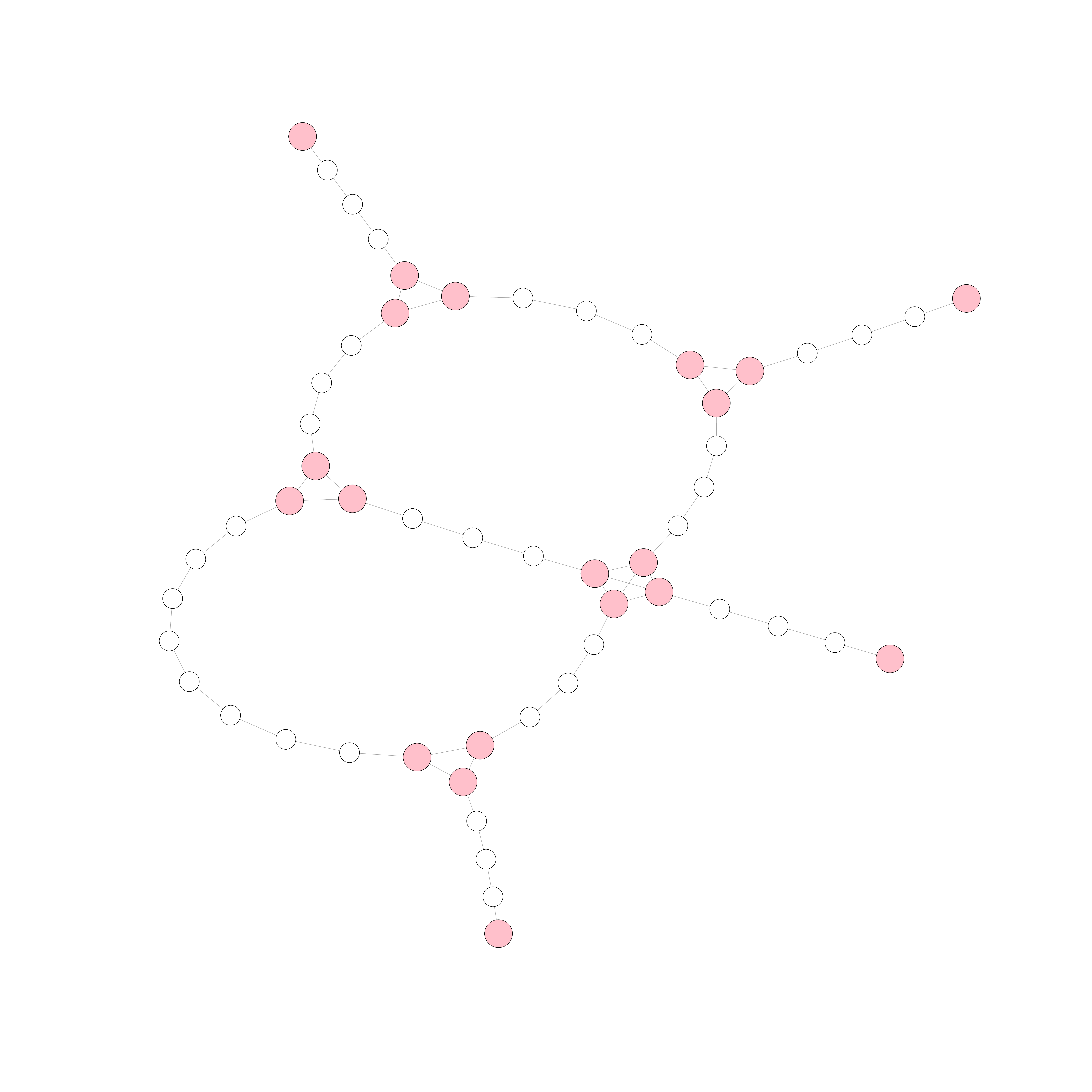}
\caption{The graph $G_\BFr$ (middle), where $r_e=4$ for all $e\in E$, corresponding to the graph $G$ (left) from Figure~\ref{fig:high_res_example}, and the respective line graph $\bar G_\BFr$ (right).
In the middle plot, red vertices correspond to vertices in $G$, and white vertices correspond to vertices that are not in $G$.
In the right plot, pink vertices correspond to edges connecting red and white vertices in $G_\BFr$, and white vertices correspond to edges connecting white vertices in $G_\BFr$.} 
\label{fig:high_res_example_col}
\end{center}
\end{figure}

By construction, in $G_\BFr$ there are two different types of vertices. This is illustrated in the middle plot in Figure~\ref{fig:high_res_example_col}:
red vertices corresponding to vertices of $G$, and 
white vertices corresponding to vertices that are not present in $G$.
In the corresponding line graph $\bar G_\BFr$, which is depicted in the right plot in Figure~\ref{fig:high_res_example_col}, there are two types of vertices:
pink vertices corresponding to edges in $G_\BFr$ that connect a red vertex and a white vertex, and
white vertices that correspond to edges in $G_\BFr$ that connect white vertices.

The structure of the line graph $\bar G_\BFr$ is rather simple:
to each red vertex of degree $d$ in the original graph corresponds a clique of size $d$ in $\bar G_\BFr$, and these cliques are connected via path graphs (the one corresponding to edge $e$ having $r_e-1$ vertices). Recall that $q_\BFr$, the number of edges of $G_\BFr$, can be written as $\sum_{i=1}^q(r_i-1)+2q$.

In the following we make extensive use of the following property.
Since $\bar{\BFL}_\BFr$ is a real, symmetric matrix of dimension $q_\BFr\times q_\BFr$, there exists a matrix $\bar\BFOmega$ of the same dimension such that 
\[
\bar\BFOmega^\top \bar{\BFL}_\BFr \,\bar\BFOmega =
\diag\{\BFell\},
\]
where \[\BFell=(\ell_1, \dots, \ell_{q_\BFr}) = (\ell_{1,1}, \dots, \ell_{1,r_1-1}, \ell_{2,1}, \dots, \ell_{2,r_2-1}, \dots, \ell_{q+1,1}, \dots, \ell_{q+1, 2q}).\]
Proposition~\ref{prop:eigendecomposition} (see Appendix A)
tells us, in particular, that if $r:=\min\{r_1,\ldots,r_q\}$ is large, then $\ell_{i,j}$ is well approximated by
$
4 \sin(\pi (j-1)/(2\, r_i))^2, j = 1, \dots, r_i, i=1,\dots,q,
$
with the remaining $2q$ eigenvalues being bounded by the maximal degree of a vertex in $G$.
In other words, the matrix $\bar\BFOmega^\top \bar{\BFL}_\BFr \,\bar\BFOmega$ and the matrix $\BFOmega^\top \bar{\BFL}_\BFr \,\BFOmega$, with $\BFOmega$ as defined in Proposition~\ref{prop:eigendecomposition}, are close. Informally, this means that the eigenvalues of $\bar{\BFL}_\BFr$ are asymptotically (as $r\to\infty$, that is) going to coincide with eigenvalues specified above. 

\subsubsection{Decomposition of $\BFSigma_\BFr^{(\BFn)}$}
\label{sec:asymptotics:spectra:covariance}

The spectrum of $\BFSigma_\BFr^{(\BFn)}$ is highly dependent on the structure of the variance-covariance matrix $\BFSigma_\BFr^{(\BFn)}$.
Here, we consider the model specified in~\eqref{eq:our_covariance} in Section~\ref{sec:parameter_models:covariance}, where $\BFSigma_\BFr^{(\BFn)} =n^{-1} \diag\{\sigma_1^2\BFone_{r_1}, \dots, \sigma_q^2\BFone_{r_q}\}$; alternative models can be handled as well at the expense of a substantial amount of additional notation and computations but this model should be flexible enough for any situation where we can think of the $n$ travel times being collected at each edge as being collected from a random subset of all particles traversing that edge.

It can immediately be seen that, by definition of the matrices $\BFOmega$ that feature in the proof of Proposition~\ref{prop:eigendecomposition},
\[
\BFOmega^\top \BFSigma_\BFr^{(\BFn)} \,\BFOmega =
\diag\{\BFD_1, \dots, \BFD_q, \BFD_{q+1}\},
\]
where the first blocks are given by
$\BFD_i := n^{-1}\,\sigma_i^2\, \BFI_{r_i-1}$,
for $i=1,\dots, q$, and where the last block is
$\BFD_{q+1} := n^{-1}\, \diag\{\sigma_{s(1)}^2, \dots, \sigma_{s(q)}^2\}$,
with $s(\cdot)$ denoting some given permutation of the edges $\{1,\dots,q\}$.

\subsubsection{Decomposition of $\BFH(\lambda,\BFtheta)$}
\label{sec:asymptotics:spectra:smoother}

Based on the decompositions from the previous two sections, we can now also (approximately) diagonalize our smoother matrix $\BFH(\lambda,\BFtheta)$. More concretely, up to a controllable error, $\BFH(\lambda,\BFtheta)\approx\BFOmega\,\diag\{\BFh\}\,\BFOmega^\top $, where
\[
h_{i,j} := 
h_{i,j}(\lambda,n) = 
\frac1{\displaystyle 1+4\,\frac\lambda n \sigma_i^2 \sin\Big(\frac{\pi(j-1)}{2\,r_i}\Big)^2}, 
\qquad j=1,\dots,r_i, \quad i=1,\dots, q+1,
\]
are the entries of $\BFh$.
We refer to Lemma~\ref{lem:series_sine} (see Appendix A) for the precise statement.

\subsection{Performance of the estimation procedure}
\label{sec:asymptotics:consistency}

We conclude this section by presenting our main result. 
It addresses the consistency of the estimators for the expected travel times and for the variances that we proposed in Section~\ref{sec:inference}. 
The proof of Theorem~\ref{theo:moments} can be found in Appendix~\ref{app:proofs}. 

\begin{theorem}
\label{theo:moments}
Suppose that 
\begin{equation}
\label{eq:Xn}\BFX^{(\BFn)} \sim 
{\mathscr N}\Big(\BFmu_{0,\BFr},\; n^{-1}\diag\big(\sigma_{0,1}^2\BFone_{r_1}, \dots, \sigma_{0,q}^2\BFone_{r_q}\big)\Big)\end{equation}
for some $\BFmu_{0,\BFr}\in\mathscr{M}_{\BFr}(C)$.
Consider then $n\in\mathbb{N}$, and $\BFr\in\mathbb{N}_0^q$ such that $n = o\big(\min_{i=1,\dots,q}r_i\big)^2$ and define the collection
\[
\Lambda_{n,\BFr}: = 
\Big\{\lambda>0: n = o(\lambda), \lambda = o\big(\min_{i=1,\dots,q}r_i\big)^2\Big\}.
\]
Consider also
\[
\hat\lambda := 
\arg\min_{\lambda\in\Lambda_{n,\BFr}} {\rm GCV}(\lambda),
\]
for ${\rm GCV}(\lambda)$ as defined in~\eqref{eq:GCV}, as well as $\hat\BFtheta = (\hat\sigma_1^2, \dots, \hat\sigma_q^2)$ with each $\hat\sigma_i^2 = \hat\sigma_i^2(\hat\lambda)$ defined as in~\eqref{eq:variances_estimators}, and finally, define \[\hat\BFmu_\BFr := \hat\BFmu_\BFr\big(\hat\lambda, \BFSigma_\BFr^{(\BFn)}(\hat\BFtheta)\big).\]
Then, as long as either $n\to\infty$ or $\min_{i=1,\dots,q}r_i\to\infty$,
$\hat\BFmu_\BFr$ is consistent in probability for $\BFmu_{0,\BFr}$, and each
$\hat\sigma_i^2$ is consistent in probability for $\sigma_{0,i}^2$.
\end{theorem}

\section{Numerical validation of estimation procedure}
\label{sec:validation}
In this section we exemplify the performance and some properties of the estimation procedure from Section~\ref{sec:inference}, using a selection of illustrative network instances. Each of our examples aims to assess a specific feature of the estimator. To simplify the interpretation of the results, we consider relatively small networks, but we emphasize that the computational burden of our algorithm, being linear in the data, is low. 

First we outline the data generation mechanism for the numerical experiments. Each example corresponds to 
\begin{itemize}
    \item[$\circ$] a graph $G$,
    \item[$\circ$] a set of resolution parameters $\BFr$,
    \item[$\circ$] a sample size vector $\BFn$,
     \item[$\circ$] the length of each of the edges,
    \item[$\circ$] per edge a velocity function (explained in detail below), determining the per-edge mean travel time,
    \item[$\circ$] the variance of the per-edge travel time.
\end{itemize}

We proceed by explaining the concept of the velocity function. This describes the  expected instantaneous velocity of particles as they traverse the corresponding edge. A velocity function is defined on the closed interval $[0,1]$. It provides the expected instantaneous velocity at
 the relative position $x\in[0,1]$ on the edge under consideration, for each  edge $e\in E$ (starting from $v_a(e)$; note that this function is not necessarily symmetric). These functions allow us to model the expected instantaneous velocities of particles as being non-constant, which is in line with the idea that the level of congestion will generally not be evenly spread throughout an edge. Of course, what we are actually interested in are the travel times, but knowing the lengths of road segments we can easily switch between velocities and travel times. The advantage of working with velocities is that while travel times scale with the length of the edges, instantaneous velocities do not, and are therefore more intuitive.
To be clear, the precise shape of the velocity function is not our target for inference. Instead, 
the velocity functions only serve  the purpose of allowing us to generate ground truths for the examples that follow.

Combining all of the above, we can determine a vector $\BFmu_\BFr$ and a variance-covariance matrix $\BFSigma_\BFr^{(\BFn)}$ that act as the ground truth for that example.
We then generate a sample of travel times $(X_{e,i}: i=1,\dots,n_e, e\in E_\BFr)$ from a given distribution which, when averaged at each edge, lead to a realization of a random vector $\BFX^{(\BFn)}$ with expectation $\BFmu_\BFr$ and variance-covariance matrix $\BFSigma_\BFr^{(\BFn)}$.
The data vector $\BFX^{(\BFn)}$ in turn leads to an estimate $\hat{\BFmu}_\BFr$ of $\BFmu_\BFr$ and, consequently, an estimate $\hat{\BFmu}$ of $\BFmu$.

To compare the estimate with the true ground truth, we use the relative squared error:
\[
\text{R}_e(\hat{\BFmu}_\BFr) := \left(\frac{\hat{\mu}_{\BFr,e}-\mu_e}{\mu_e}\right)^2, \qquad e\in E_{\BFr}.
\]
Since $\hat{\BFmu}_\BFr$ is random, the relative squared error only measures the performance of the estimator for a particular realization of $\hat{\BFmu}_\BFr$. What we are actually interested in is the expected relative squared error of $\hat{\BFmu}_\BFr$:
\[
\text{RSE}_e = \mathbb{E}\left(\frac{\hat{\mu}_{\BFr,e}-\mu_e}{\mu_e}\right)^2, \qquad e\in E_{\BFr}.
\]
Therefore, we carry out the estimation procedure for $M$ independent samples $\BFX_j^{(\BFn)}$ and compute the average of the errors of the respective estimates $\hat{\BFmu}_{\BFr,j}$, $j = 1,\dots,M$. It follows from the law of large numbers that
\[
\bar{\text{R}}_{e,M} = 
\frac{1}{M}\sum_{j=1}^M \text{R}_e(\hat{\BFmu}_{\BFr,j}) 
\xrightarrow{\rm a.s.} 
\text{RSE}_e, \qquad e\in E_{\BFr}.
\]
For $M$ large, we obtain accurate approximations for the true estimation errors for each edge. In the examples that follow we have set $M=100\,000$. The aggregate simulation time, corresponding to all examples appearing in this section, was as low as two hours on an ordinary laptop.

\subsection*{Example 1: A first display of the estimation procedure output}
\label{sec:validation:Ex1}
This example introduces an elementary network and aims to confirm that the estimator performs as intended.

\begin{figure}[!ht]
	\begin{center}
	\includegraphics[width = 0.5\textwidth, trim={2cm 2cm 3cm 1cm}, clip]{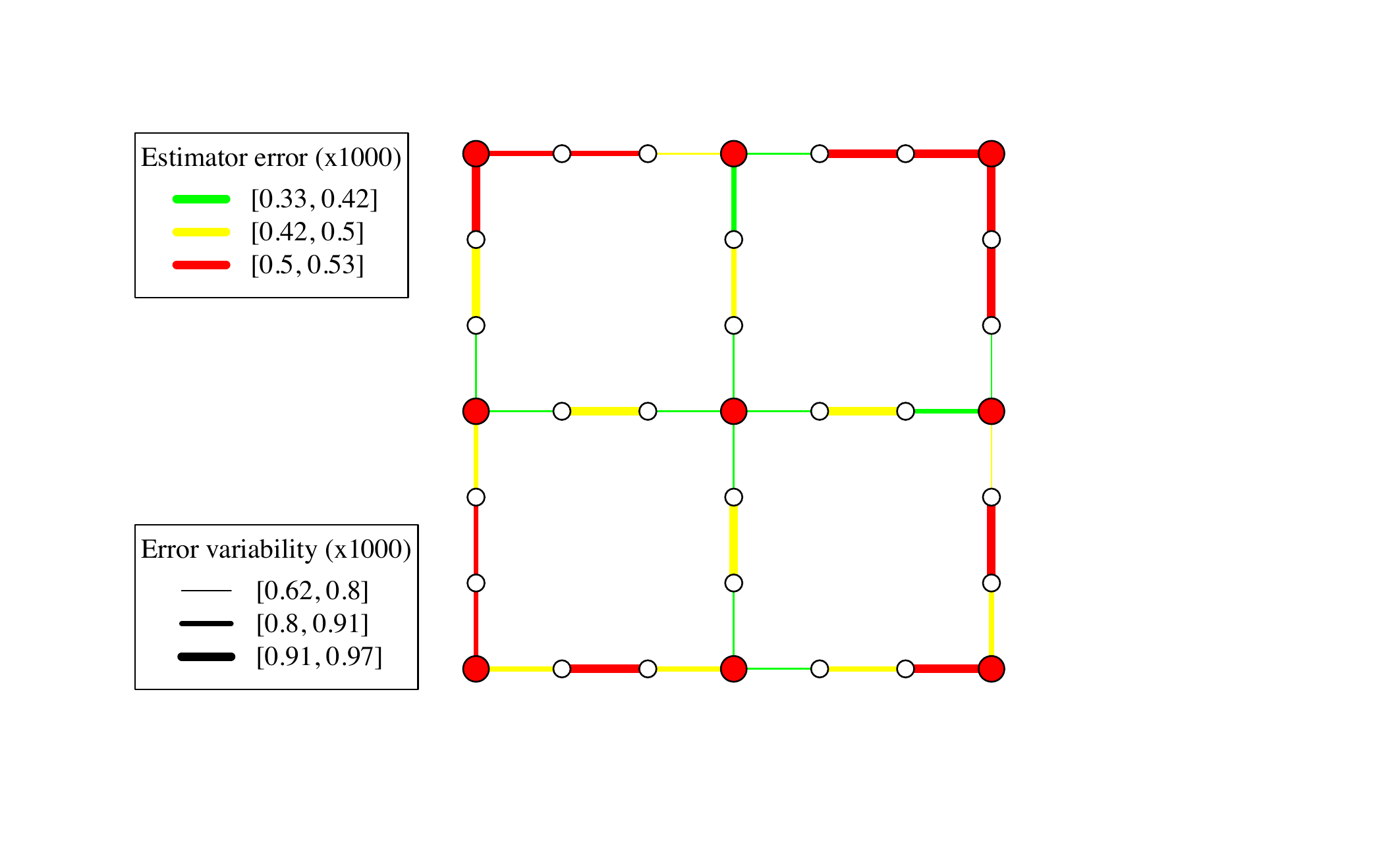}%
	\includegraphics[width = 0.5\textwidth, trim={2cm 2cm 3cm 1cm}, clip]{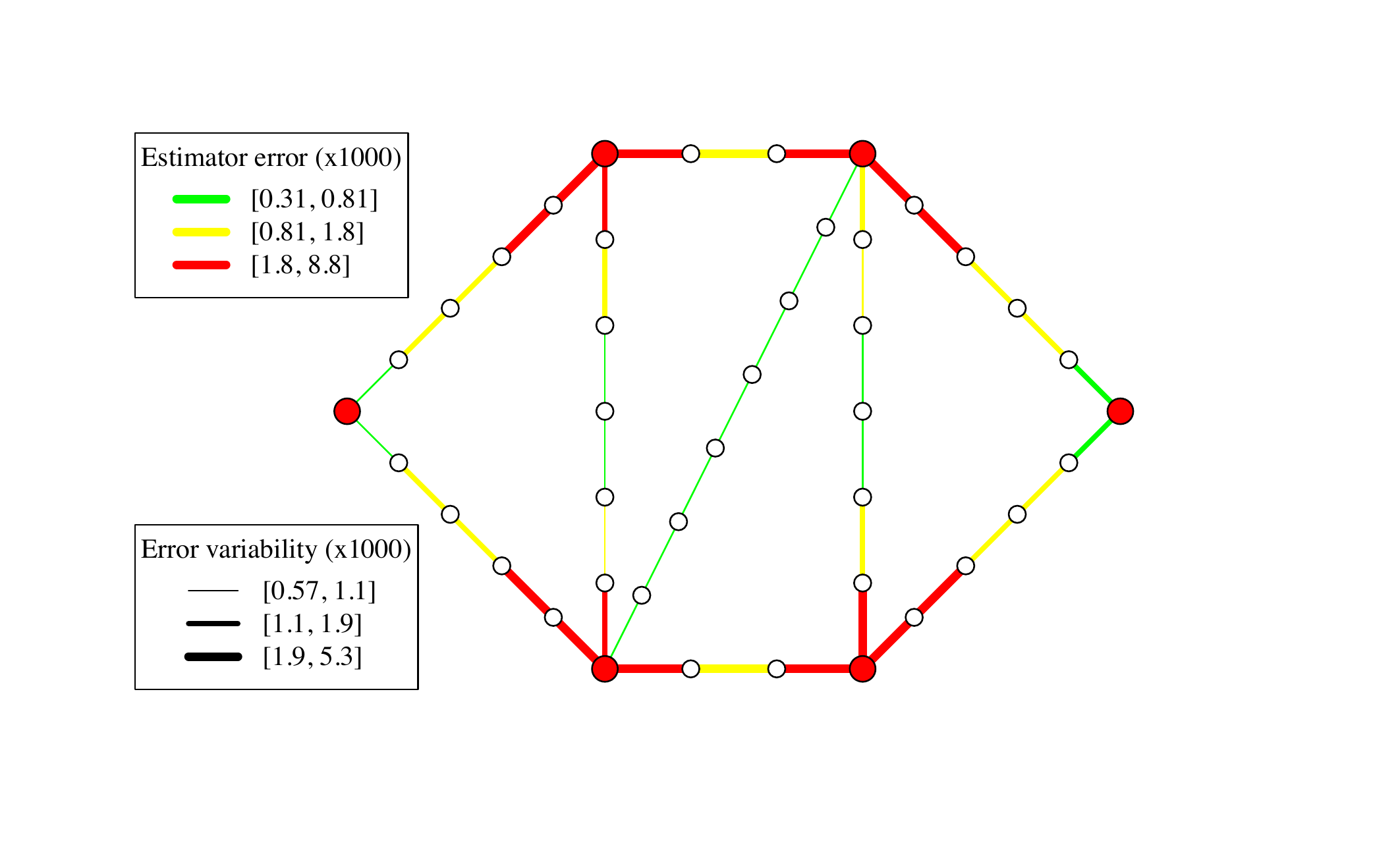}
	\caption{Approximations of the expectation and standard deviation of the relative squared errors corresponding to Example 1.} 
	\label{fig:Ex1}
	\end{center}
\end{figure}

As a first graph we choose a 2 by 2 lattice with edges that each represent a road with a length of 1 kilometer; see the left plot in Figure~\ref{fig:Ex1}. The resolution of each edge is set to $r_e=2$, which means that the higher resolution graph is constructed by replacing every edge in $G$ by a path graph consisting of 3 edges of equal length to obtain $G_\BFr$. The vertices that are added to the higher resolution graph are colored white, whereas the red vertices were already part of the original graph. Besides the lattice, we also consider the graph depicted in the right plot in Figure~\ref{fig:Ex1}. For this second graph, the scale of the figure is such that the length of the two shortest edges of its original graph correspond to 1 kilometer. Since the edges of the original graph do not all have equal length, we choose different resolutions for different edges to ensure that each edge on the graph corresponds to a road segment with approximately the same length. 

In this example, we assume that particles traverse the edges with a constant expected velocity of 30 km/h. The variance of the travel time per kilometer is   $0.01^2$ hour$^2$, at each of the sub-edges. Hence, the sub-edges that arise in the higher resolution graph will also each be traversed with an expected velocity of 30 km/h. We let the travel time variance  corresponding to the sub-edges  be given by $0.01^2$ multiplied by the length of the sub-edge.
Using the lengths of the edges, we translate the expected velocities to expected travel times. Performing a conversion from hours to seconds, each simulation for the left plot in Figure~\ref{fig:Ex1} therefore consists of sampling independent $X_{e,i}\sim{\mathscr N}(40,36^2/3)$, $i=1,\dots,n_e$, $n_e=100$, $e\in E_{\BFr}$, collecting these in $\BFX^{(\BFn)}$, and from this computing the estimate $\hat{\BFmu}_\BFr$ as well as their respective relative squared errors $\text{R}_e(\hat{\BFmu}_\BFr)$ for each edge $e\in E_{\BFr}$. Note that the moments of the samples for the right plot in Figure~\ref{fig:Ex1} are slightly different since the lengths of the sub-edges differ.

This procedure was repeated $M$ times to obtain a sample of $M$ relative squared errors.
These are reported in Figure~\ref{fig:Ex1} where we have summarized the sample mean and sample standard deviation of the relative squared errors at each edge of $G_\BFr$.

Rather than reporting the means $\bar{\text{R}}_{e,M}$ (our proxies for $\text{RSE}_e$) for each edge in $E_\BFr$, we color  the respective edge based on the value of $\bar E_{e,M}$.
Each color corresponds to a range of relative errors, with the break points that define the ranges being three equally spaced quantiles of the sampled $\{\bar{\text{R}}_{e,M}: e\in E_\BFr\}$.
So, for instance, in the left plot in Figure~\ref{fig:Ex1}, $0$, $0.42$, $0.50$, and $0.53$ are respectively, the minimum, $0.33$-, $0.67$-quantile, and maximum of the sampled $\{\bar{\text{R}}_{e,M}: e\in E_\BFr\}$.

In the same spirit, rather than reporting the $|E_\BFr|$ sample standard deviations of the $M$ relative squared errors obtained at each of the edges in $E_\BFr$ in our Monte Carlo simulation, we set the thickness of the respective edge based on the values of the $|E_\BFr|$ sample standard deviations.

\medskip


The left plot in Figure~\ref{fig:Ex1} shows that the edges near the center are colored green and are thin; 
this means that within this graph, we conclude the relative squared error of the edges at the center to have smaller expectation and standard deviation.
In contrast, the edges in the corners are red and thick;
this indicates that within this graph, we observe the relative squared error of the edges at the corners to have larger expectation and standard deviation.
Importantly, the above does not mean that the travel times of the corner edges are poorly estimated, but rather that they have higher relative squared error when compared the edges at the center.
Indeed, the legend reveals that the expected relative squared errors of any of the edges are low: in our simulation, they do not exceed $0.53\times 10^{-3}$ or $0.053\%$ relative squared error.
Also bearing in mind the small error variability, we conclude that in this example the mean travel times of all edges were estimated accurately. The same conclusion holds for the right plot in Figure~\ref{fig:Ex1}, where we see that the expected relative squared errors do not exceed $8.8\times 10^{-3}$.  


A final remark is that the asymptotic result in Theorem~\ref{theo:moments} ensures that the global estimation error for the entire graph (meaning the average error across the entire graph) is low if either all sample sizes $n_e$ are large, or if all resolution parameters are large.
It is however quite instructive to look at the per edge errors as these reveal the effect of the topology of the graph and the local amount of information available.

\subsection*{Example 2: Effect of sample size on estimation error}
\begin{figure}[!ht]
	\begin{center}
		\includegraphics[width = 0.5\textwidth, trim={2cm 2cm 3cm 1cm}, clip]{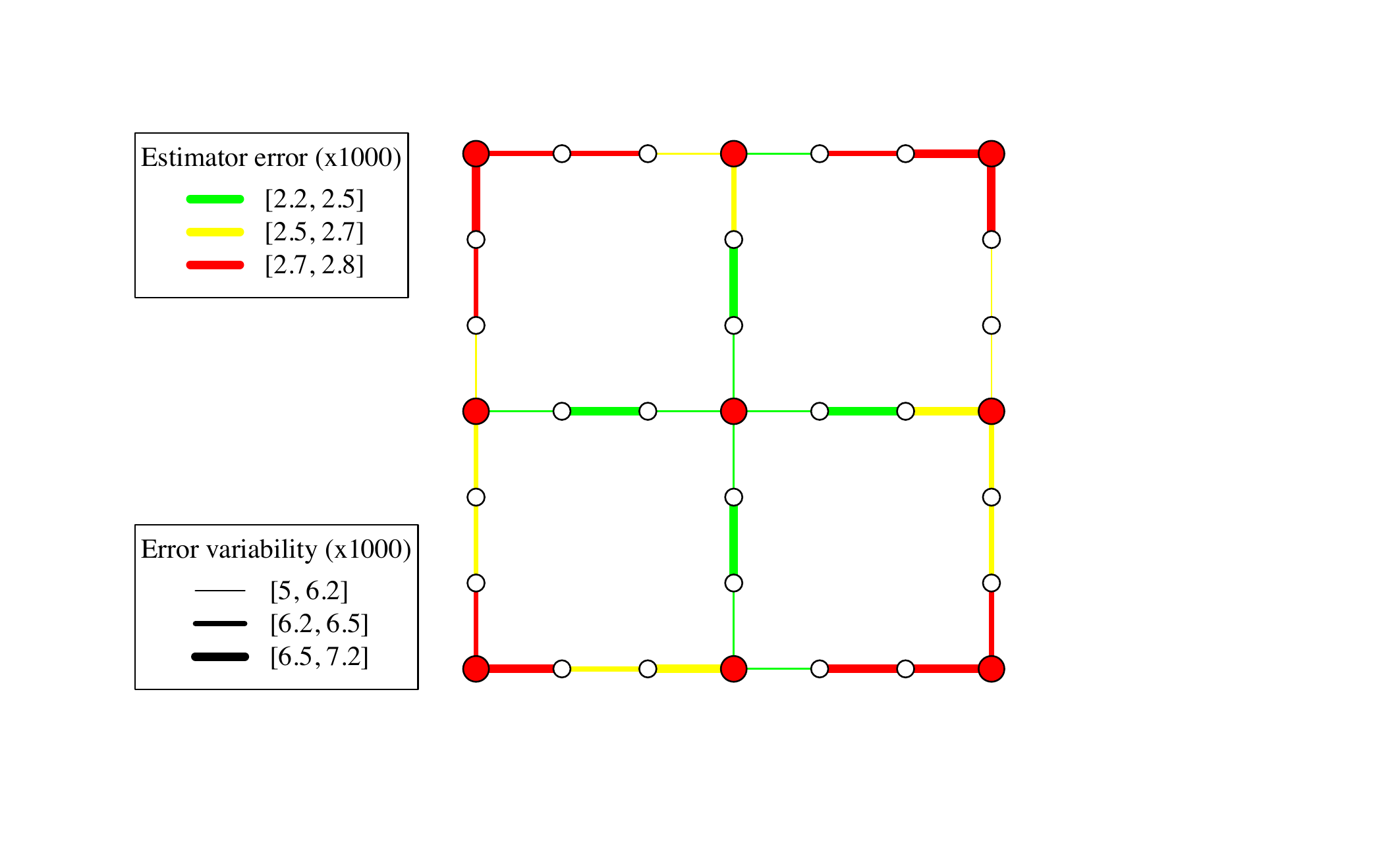}%
		\includegraphics[width = 0.5\textwidth, trim={2cm 2cm 3cm 1cm}, clip]{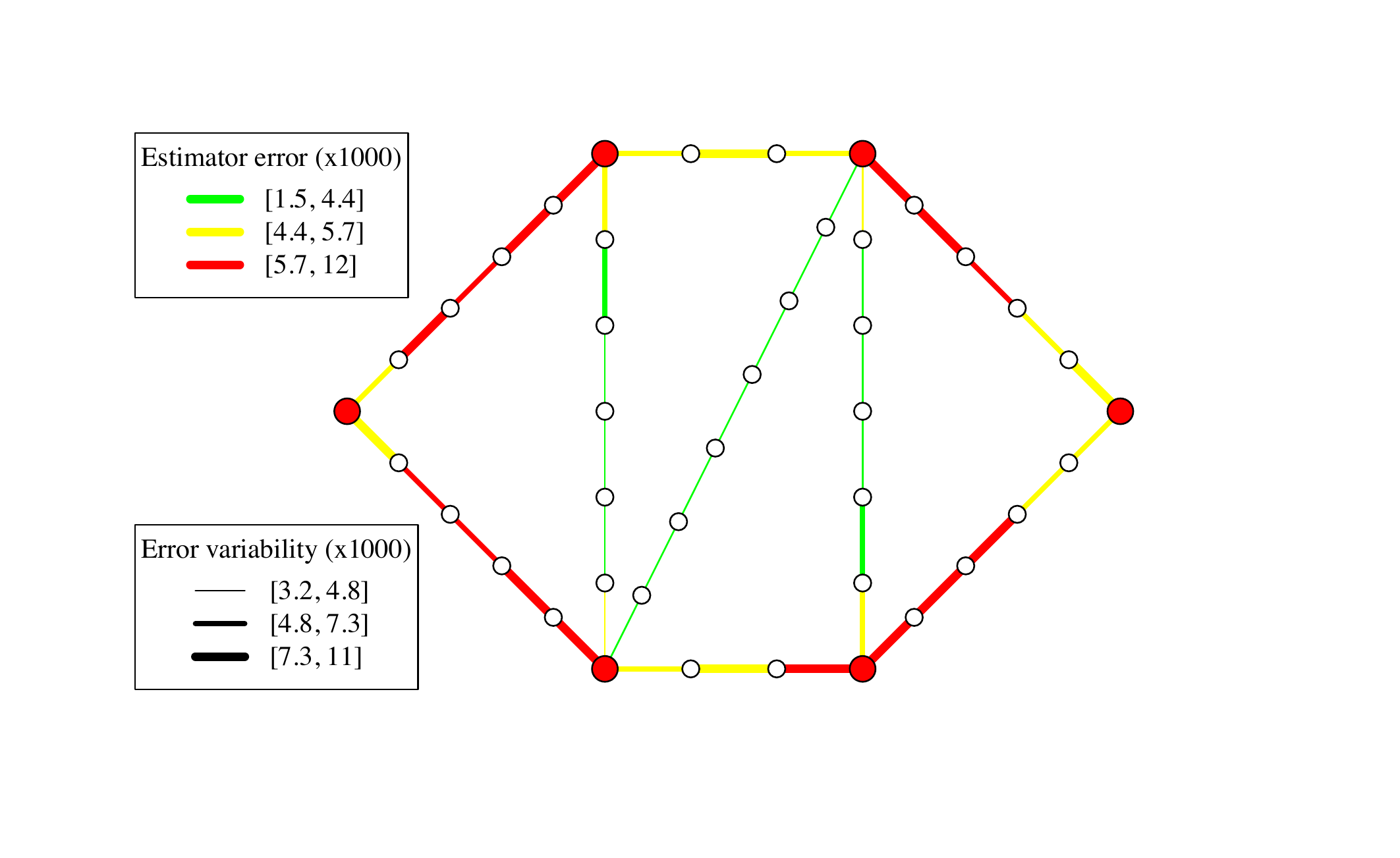}
		\caption{Estimation errors corresponding to the instances discussed in Example 1 but with a smaller sample size of $n_e = 10$ for each $e\in E_{\BFr}$.} 
		\label{fig:Ex2}
	\end{center}
\end{figure}

This example is a continuation of Example 1. It illustrates to what extent the quality of the estimates changes as a function of the sample size. Keeping the setting of Example 1 unaltered, we now choose a smaller sample size of $n_e = 10$ for each $e\in E_{\BFr}$. The results are shown in Figure~\ref{fig:Ex2}.

Comparing the ranges of the quantiles in Figure~\ref{fig:Ex2} with those in Figure~\ref{fig:Ex1}, we see that both the errors and the variability of the errors are substantially lower in the experiment with the higher sample size (i.e., the setting of Figure~\ref{fig:Ex1}). This is in line with our asymptotic results.

We also see that the effect of the topology  on the quality of the per-edge estimates remains similar; only the magnitude of the errors (scale of the error ranges) changes, with the relative magnitude of the errors within each graph (color of the edges) remaining essentially unaltered.
Later examples focus on the effect of the topology  on the estimation error, but we first consider the effect of the amount of smoothing.


\subsection*{Example 3: Effect of smoothing}
\begin{figure}[!ht]
	\begin{center}
		\includegraphics[width = 0.5\textwidth, trim={2cm 2cm 3cm 1cm}, clip]{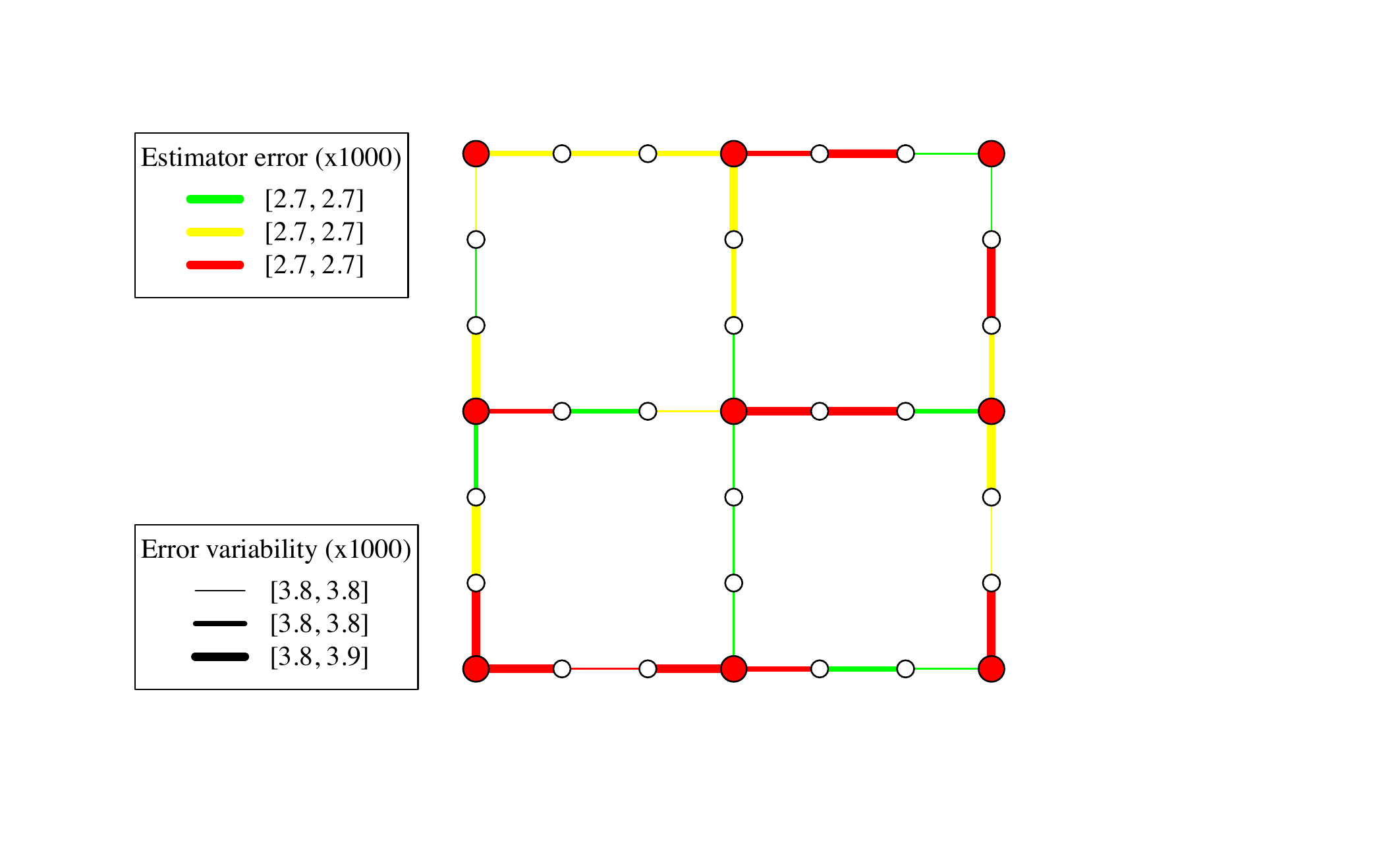}%
		\includegraphics[width = 0.5\textwidth, trim={2cm 2cm 3cm 1cm}, clip]{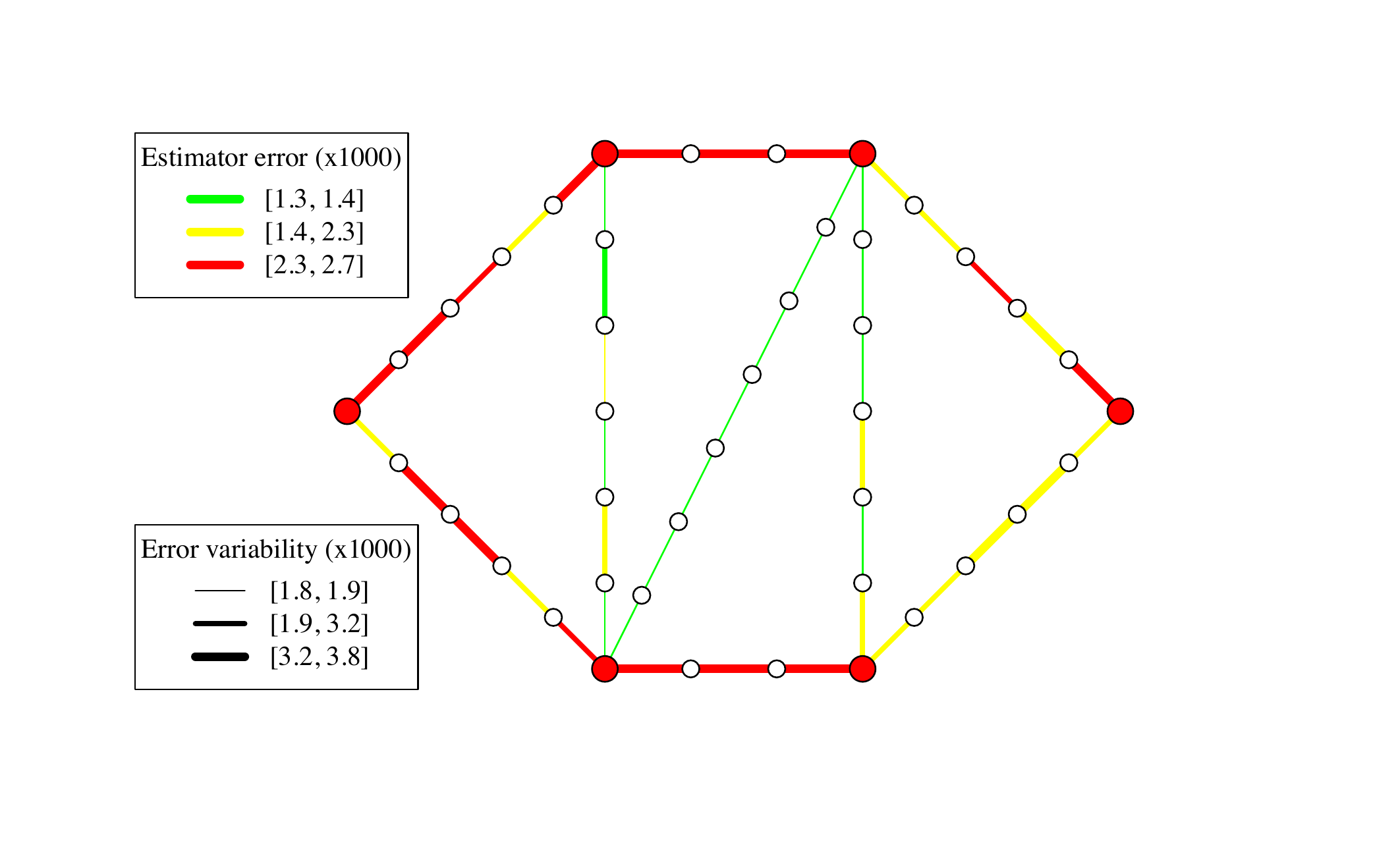}
		\caption{Estimation errors corresponding to the first instance discussed in Example 1 but with no-smoothing.} 
		\label{fig:Ex3.1}
	\end{center}
\end{figure}

This example illustrates the effect of the amount of smoothing on the resulting estimates. 
If we replace the smoother matrix by $\BFI_{q_\BFr}$, then our estimator becomes $\hat{\BFmu}_\BFr = \BFX^{(\BFn)}$.
This can be interpreted as not performing any smoothing at all, as the smoother matrix becomes an identity as $\lambda\to0$; see~(\ref{eq:mean_estimator}).
This means that in this no-smoothing case the expected time to cross each edge is estimated based on information collected at the edge under consideration only, in that there is no sharing of information across neighboring edges. 
Other than the different amount of smoothing, the setting remains the same as that of Example 1.
The output of the estimation procedure is shown in Figure~\ref{fig:Ex3.1}. 

Comparing the ranges of the quantiles of the left plot in Figure~\ref{fig:Ex3.1} with those in Figure~\ref{fig:Ex1}, we see both drastically increased errors and a higher error variability for each quantile. The highest relative squared error, for instance, is more than five times higher in the no-smoothing case. This indicates that the smoother matrix, as expected, enables us to obtain much better estimates.
We also see that the redder edges are more scattered throughout the graph, as in a no-smoothing case the topology of the graph plays no role in the estimation. In the right plot in Figure~\ref{fig:Ex3.1}, the effect of smoothing is less pronounced. Looking at the estimated average squared error across the entire graph, however, the smoothing case ($2.5$ seconds squared) outperforms the no-smoothing case ($2.7$ seconds squared).

\begin{figure}[!ht]
	\begin{center}
		\includegraphics[width = 0.5\textwidth, trim={2cm 2cm 3cm 1cm}, clip]{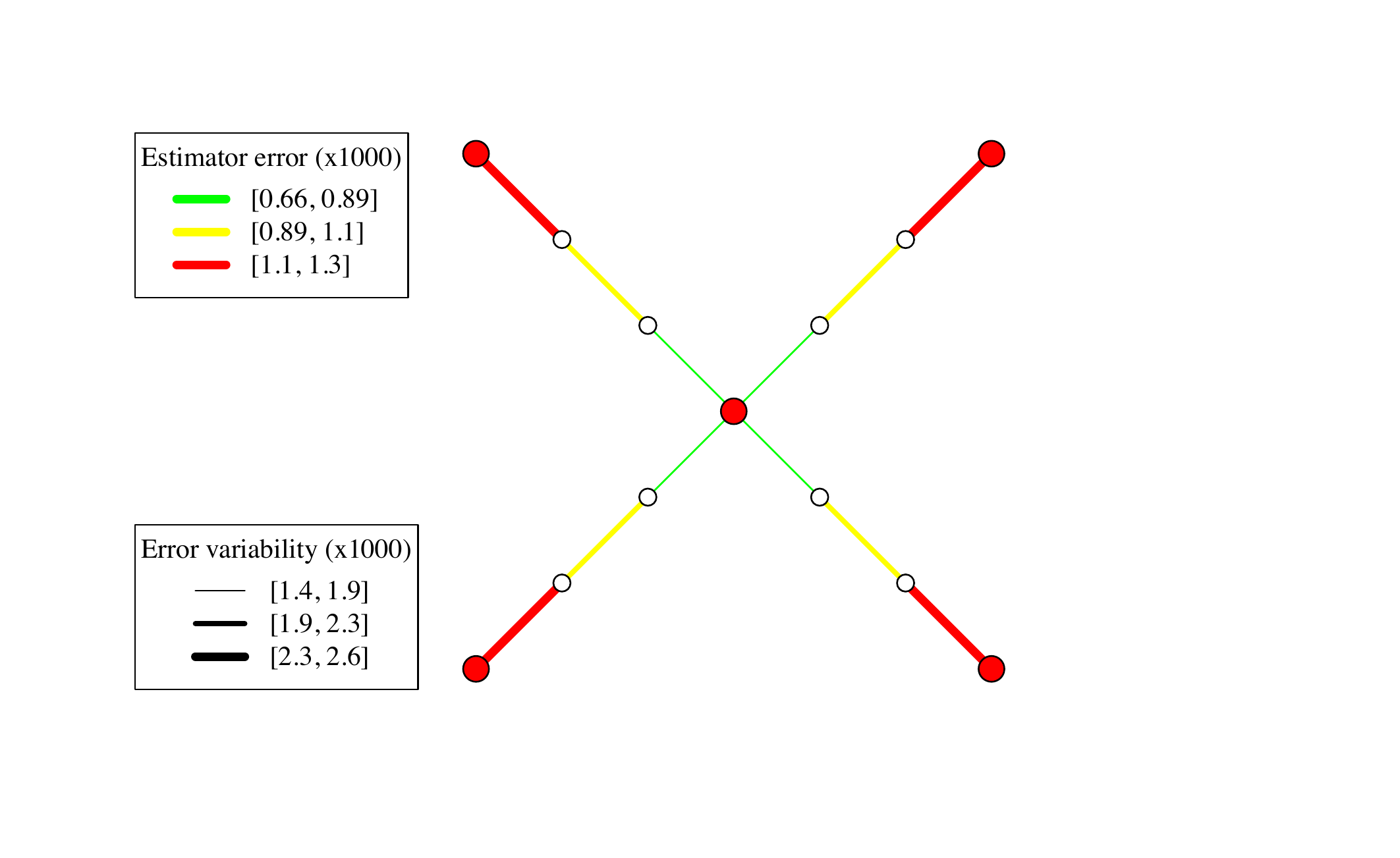}%
		\includegraphics[width = 0.5\textwidth, trim={2cm 2cm 3cm 1cm}, clip]{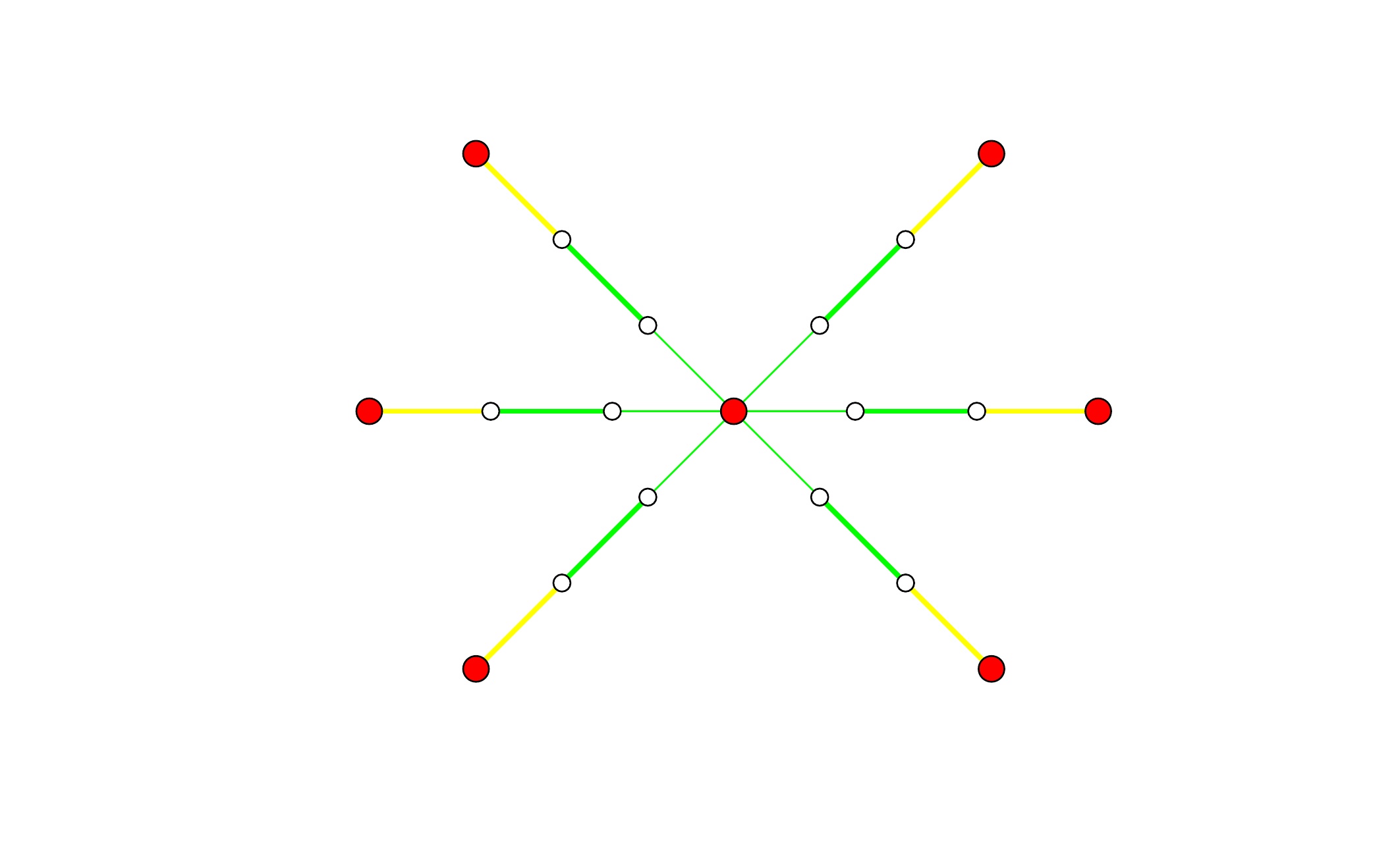}
		\caption{Illustration of the benefits of smoothing for edges with more neighbors.} 
		\label{fig:Ex3.2}
	\end{center}
\end{figure}

The reason why the smoothing parameter plays an important role, specifically in the context of this example, is that each edge has equal constant expected velocity, so that the vector $\BFmu_\BFr$ consists of equal entries whenever edges have equal length, as is the case in the left plot in Figure~\ref{fig:Ex3.1}.
The GCV based procedure that gives us a data-driven choice of $\lambda$ picks up on this, and selects a large value for $\lambda$ which in turn allows the estimate of the expected travel time at each edge to pull more information from neighboring edges.
This results in improved estimates for the expected travel times.
Figure~\ref{fig:Ex3.2} provides another example of the advantage of working with a smoothing parameter.

We use a common scale in both plots in Figure~\ref{fig:Ex3.2} to facilitate easy comparison.
We recall that red does not necessarily indicate a large error but instead an error that is \emph{comparatively} larger;
the values of $\bar{\rm R}_{e,M}$ for each edge $e$ are rather low.

If $\lambda$ is larger, then the estimated travel time of an edge incorporates the observations from more neighboring edges, enabling information to be `carried over' between these edges. The more neighbors an edge has, the more it can benefit from this effect. This effect is demonstrated in Figure~\ref{fig:Ex3.2}, where we see that edges at the center of the graph have better estimates, since these have more neighbors. Moving further away from the center, the edges become more isolated, and as a consequence these edges benefit less from the smoothing effect. 
(This effect was already visible when comparing the errors of the central edges and its less central counterparts in Figure~\ref{fig:Ex3.1} with their counterparts in Figure~\ref{fig:Ex1}.)
Note also that the positive effect of a high number of neighbors in combination with high smoothing benefits not just the edges that share the vertex with the highest degree; the effect extends to their neighbors, neighbors of neighbors, etc.

\subsection*{Example 4: Effect of over-smoothing rough signals}
\begin{figure}[!ht]
	\begin{center}
		\includegraphics[width = 0.5\textwidth, trim={2cm 2cm 3cm 1cm}, clip]{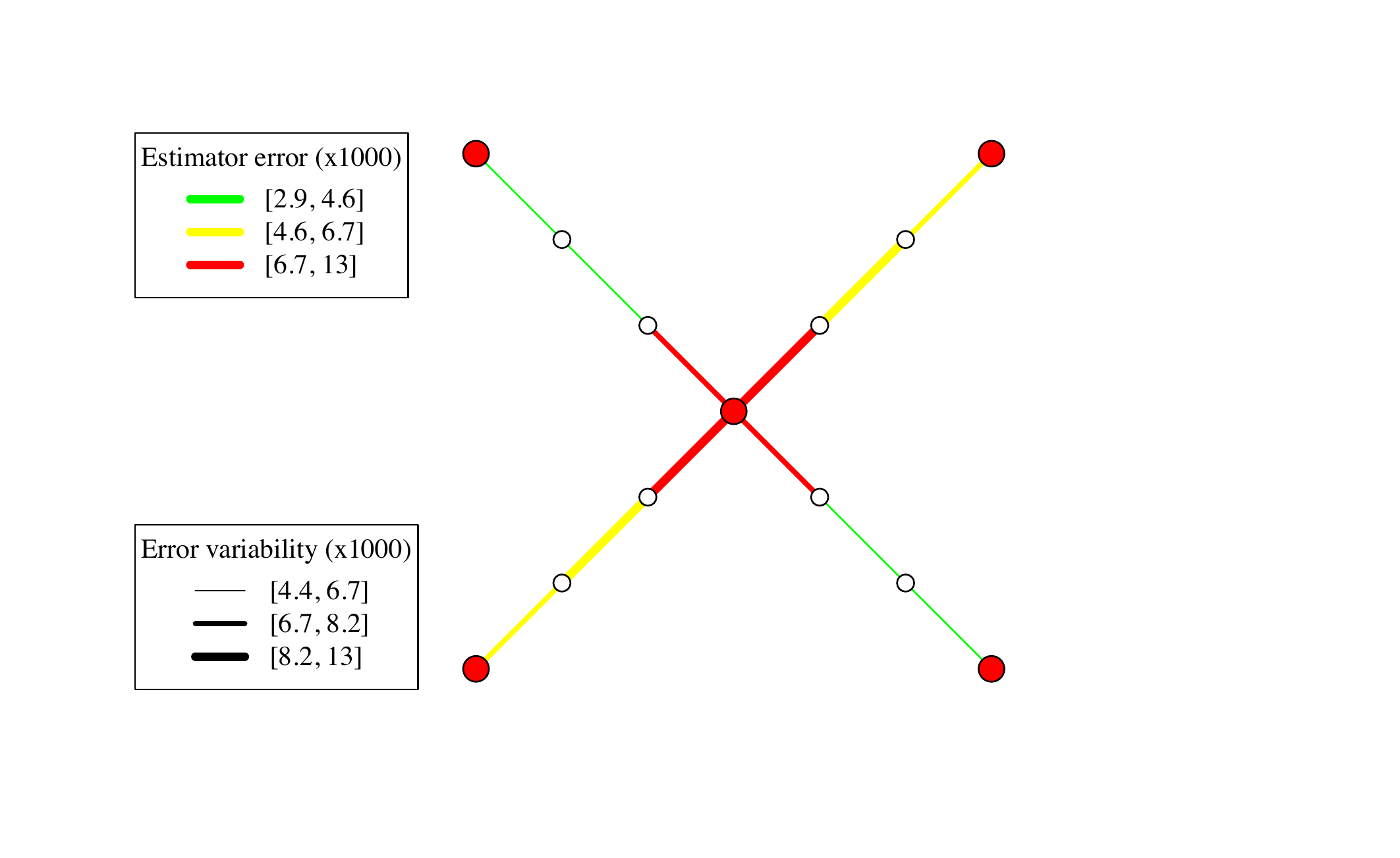}
		\caption{Illustration of the effect of lack of smoothness of the travel times. The average squared error across the entire graph is $7.4$ seconds squared.} 
		\label{fig:Ex4.1}
	\end{center}
\end{figure}

Of course, the carry-over effect that was described in Example 3 will only be beneficial if neighbors of an edge have similar expected travel times. In this example we see what happens if not all edges are traversed with the same expected velocity. Instead of assuming an expected velocity of 30 km/h on all edges (as we did in Examples 1-3), the SW-NE edges are now traversed with an expected velocity of 40 km/h. Hence, at the center of the graph we have adjacent edges traversed at different velocities. 
Figure~\ref{fig:Ex4.1} summarizes the results.

As the estimator carries over observations corresponding to lower expected travel times to edges with higher expected travel times and vice versa, the bias of the estimator at these edges increases. If we move further away from the center, we see that the relative error of the estimator decreases as these edges do have neighbors with the same expected velocity: for these edges it is beneficial to share information.

\begin{figure}[!ht]
	\begin{center}
		\includegraphics[width = 0.5\textwidth, trim={2cm 2cm 3cm 1cm}, clip]{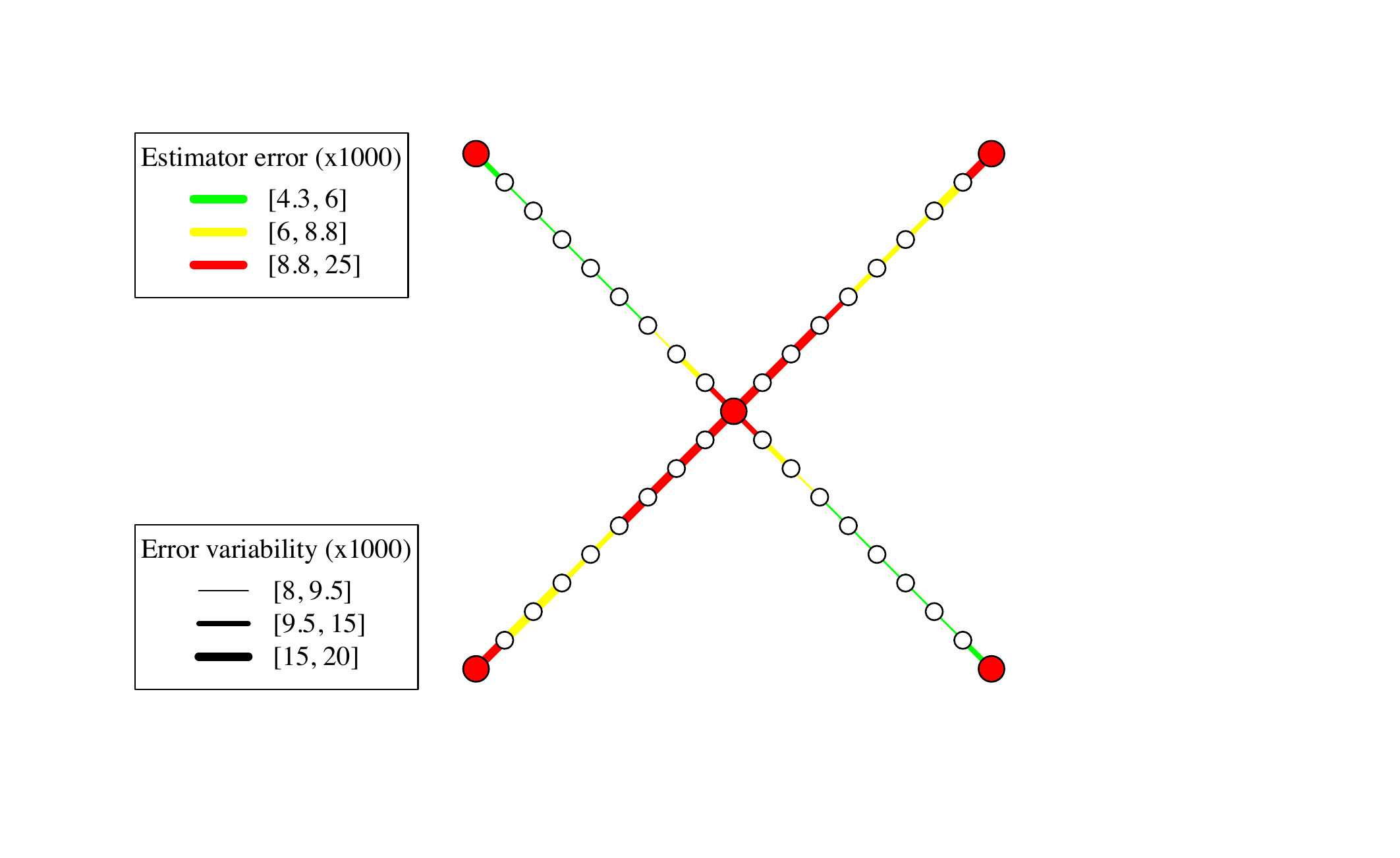}
		\caption{Estimation errors corresponding to the instance discussed in Figure~\ref{fig:Ex4.1} but with a higher resolution and with $n_e = 34$ for each $e\in E_{\BFr}$. The average squared error across the entire graph is $1.1$ seconds squared.} 
		\label{fig:Ex4.2}
	\end{center}
\end{figure}

While the expected travel time is not \emph{smooth} close to the intersection at the middle of the graph, it is so everywhere else.
This explains the large errors in Figure~\ref{fig:Ex4.1}, relative to those in the left plot on Figure~\ref{fig:Ex3.2}.
While lack of smoothness runs contrary to the principle of smoothing, its effect can be mitigated by increasing the resolution of the graph.
Figure~\ref{fig:Ex4.2} shows the effect of increasing the resolution parameters from $r_e=2$ to $r_e=8$. Note the smaller sample size, which makes the average squared error of Figures~\ref{fig:Ex4.1} and ~\ref{fig:Ex4.2} comparable since both graphs have a similar number of observations on each edge of the original graph~$G$.

We see that the increase in resolution does not eliminate the problem -- as is does not (and cannot) do away with the lack of smoothness at the intersection -- but it does isolate the higher errors to just the area close to the intersection which is the area where the expected travel time, as a function on the edges of the graph, lacks smoothness.
So, increasing resolution does not eliminate the problem but it rather  concentrates it.
In the following example we illustrate how this is true for less trivial choices of the expected travel times.

\subsection*{Example 5: Capturing inhomogeneous speeds}
In all of the preceding examples we assumed that particles traverse each of the edges at a constant expected velocity. This is clearly not realistic as drivers reduce their speed as they approach curves or intersections.

\begin{figure}[!ht]
\begin{center}
\includegraphics[height=6cm]{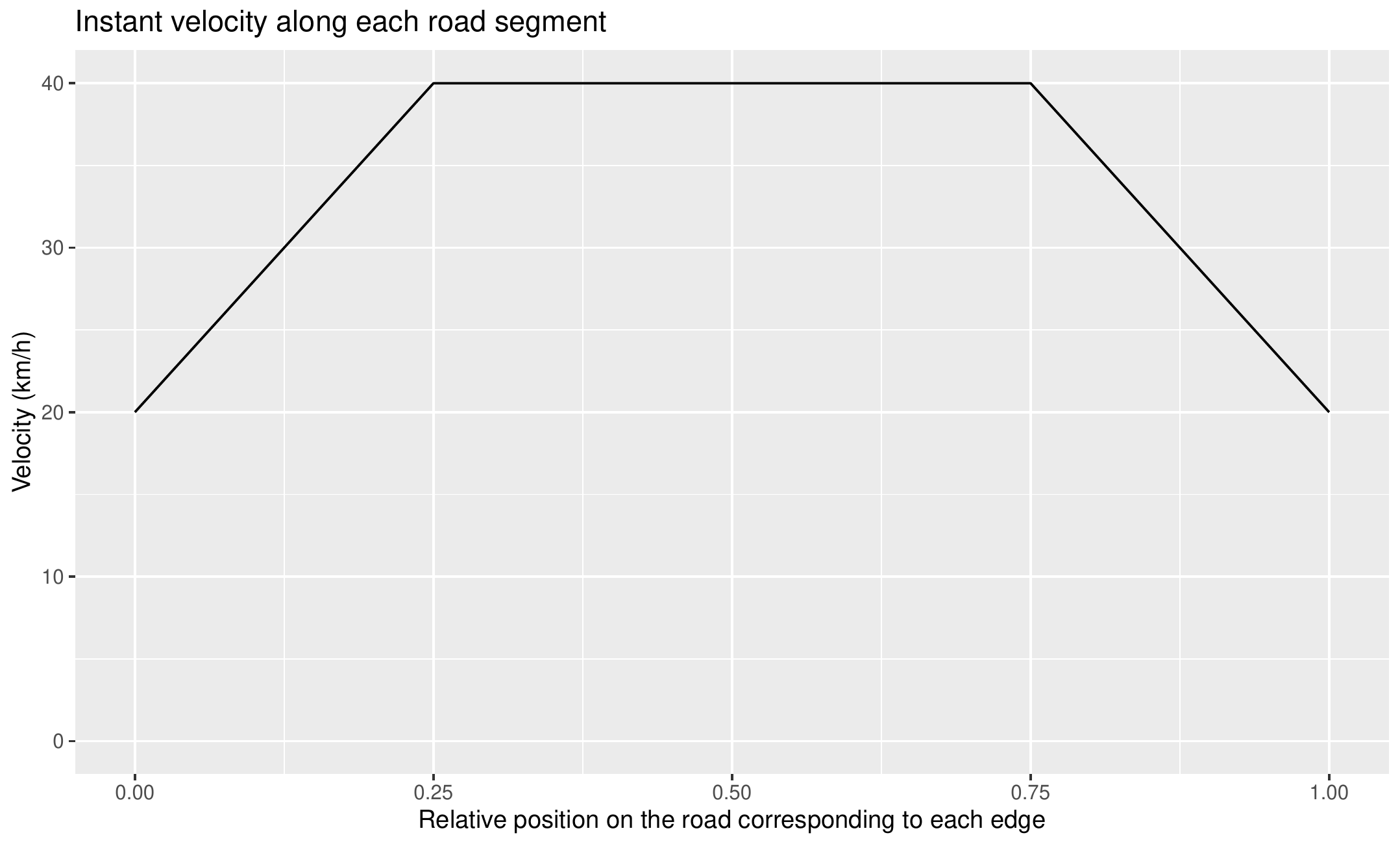}
\caption{The non-constant velocity function used in Example 5. The relative position `0' corresponds to the instantaneous velocity at the start of the road corresponding to the edge, while the relative position `1'  corresponds to the instantaneous velocity at the end of the road segment.} 
\label{fig:Ex5velocity}
\end{center}
\end{figure}

In this example we return to the graph from Example 1. We model each particle to enter each road segment (represented by an edge in the original graph $G$) at a relatively low velocity, accelerate until they reach a higher velocity, and then decelerate again when approaching the end of an edge. Specifically, at each edge we model the expected velocity using the trapezoid function depicted in Figure~\ref{fig:Ex5velocity}.

This function is the same for every edge in $G$, and is used to compute the expected travel times needed to traverse each of the road segments represented by the edges in $G_\BFr$, which are now going to be different from edge to edge.
It is noted that in fact any choice of the expected instantaneous velocities that leads to smoothly varying expected travel times is allowed in our setting.

\begin{figure}[!ht]
	\begin{center}
		\includegraphics[width = 0.5\textwidth, trim={2cm 2cm 3cm 1cm}, clip]{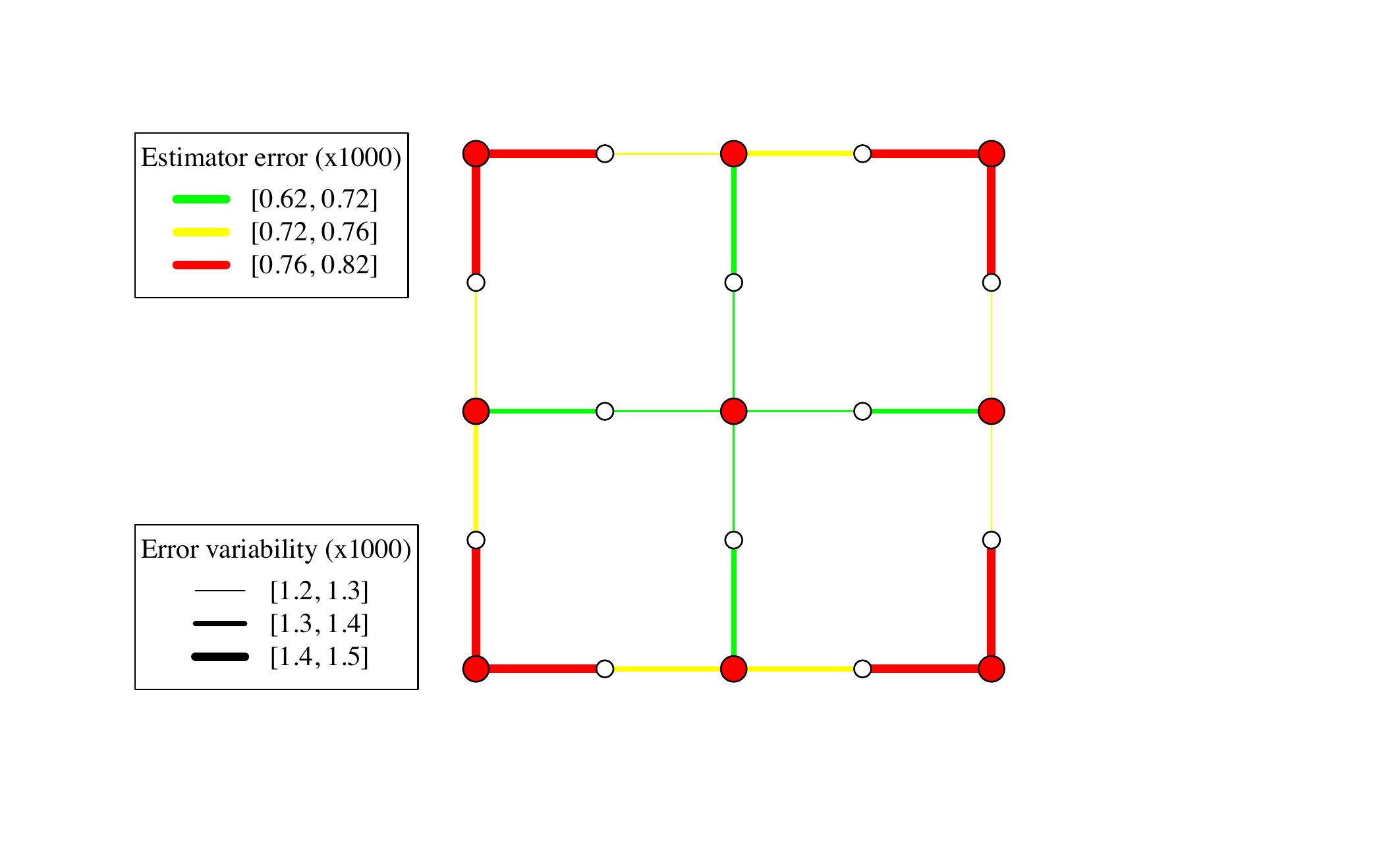}%
		\includegraphics[width = 0.5\textwidth, trim={2cm 2cm 3cm 1cm}, clip]{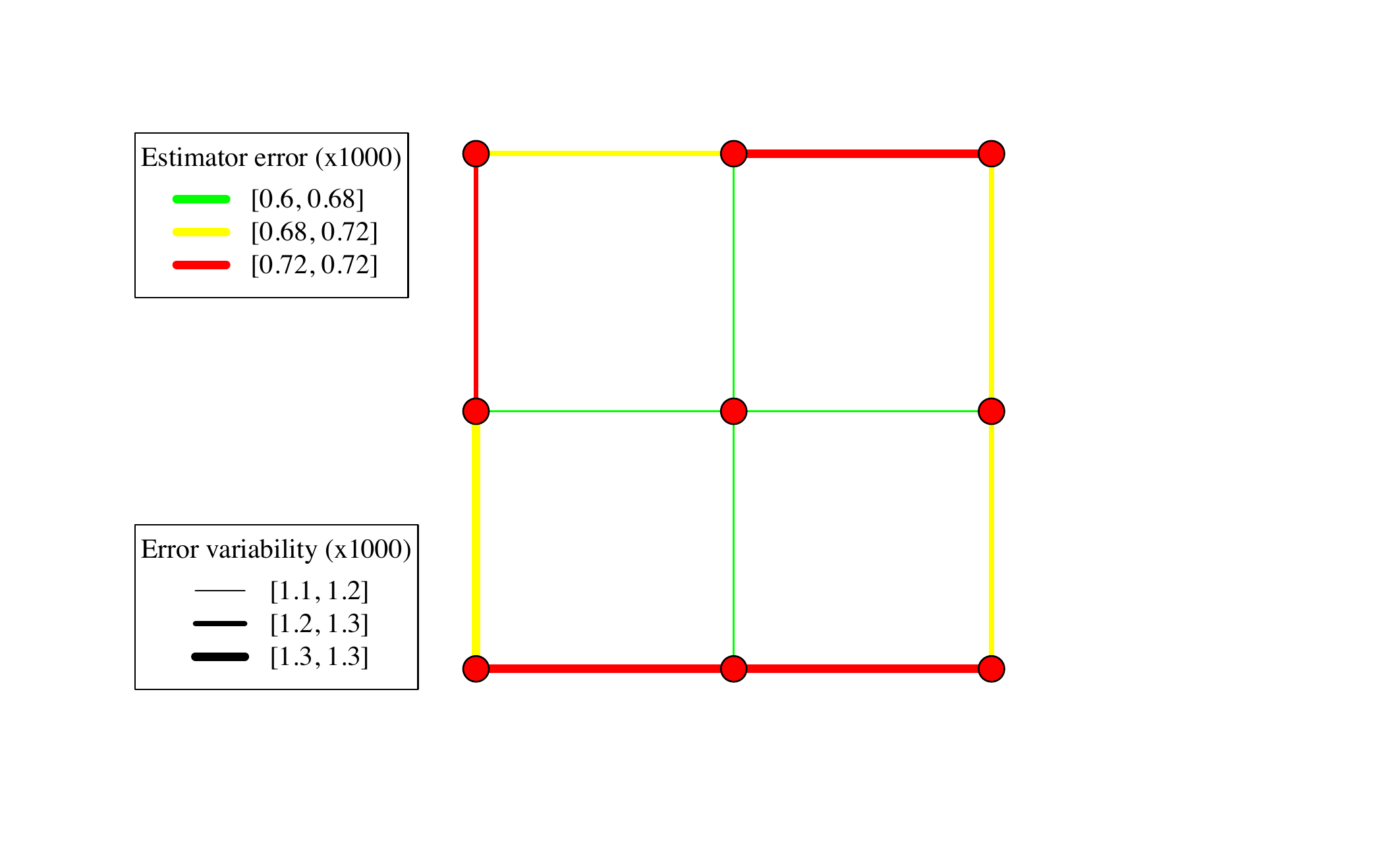}
		\includegraphics[width = 0.5\textwidth, trim={2cm 2cm 3cm 1cm}, clip]{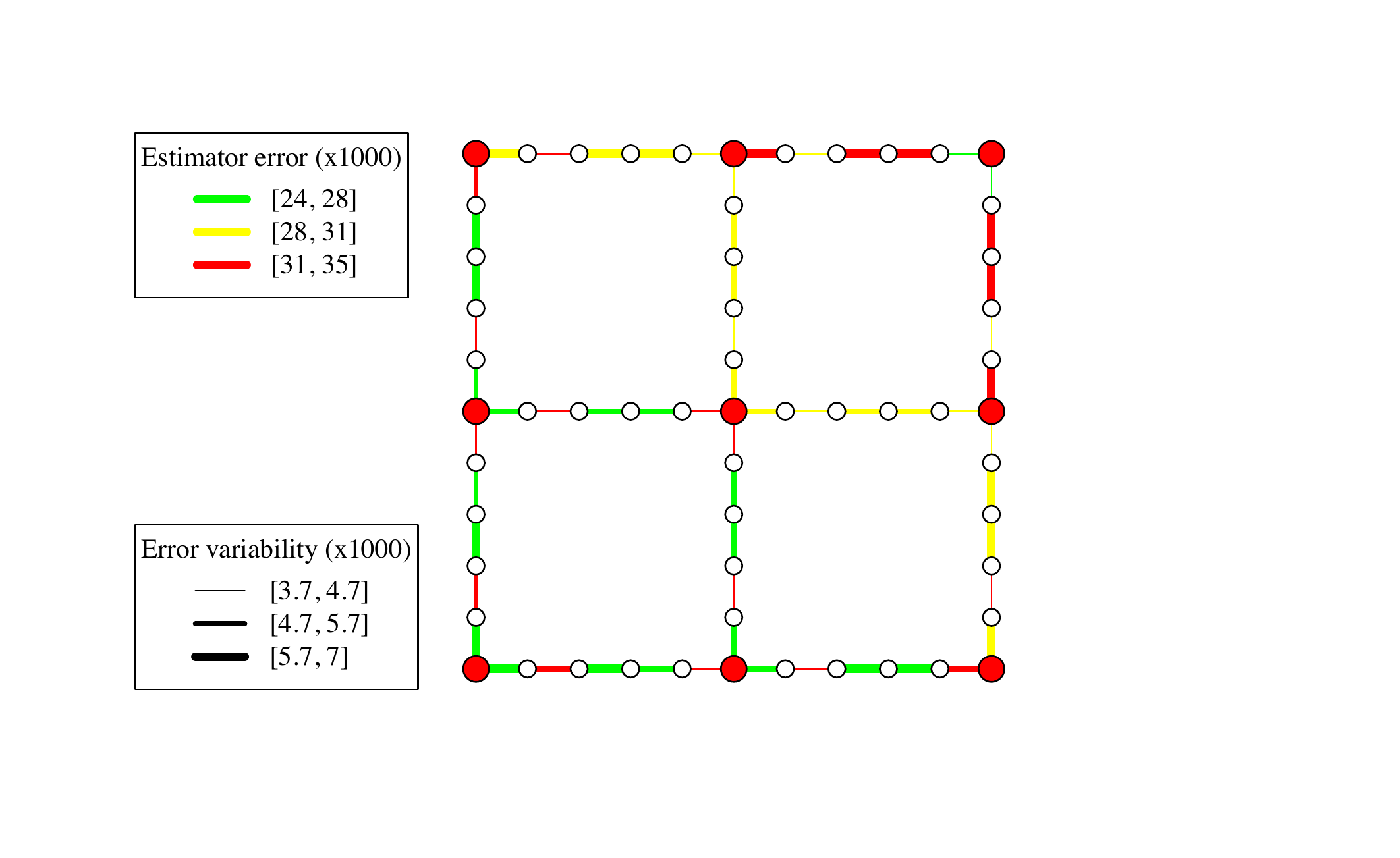}%
		\includegraphics[width = 0.5\textwidth, trim={2cm 2cm 3cm 1cm}, clip]{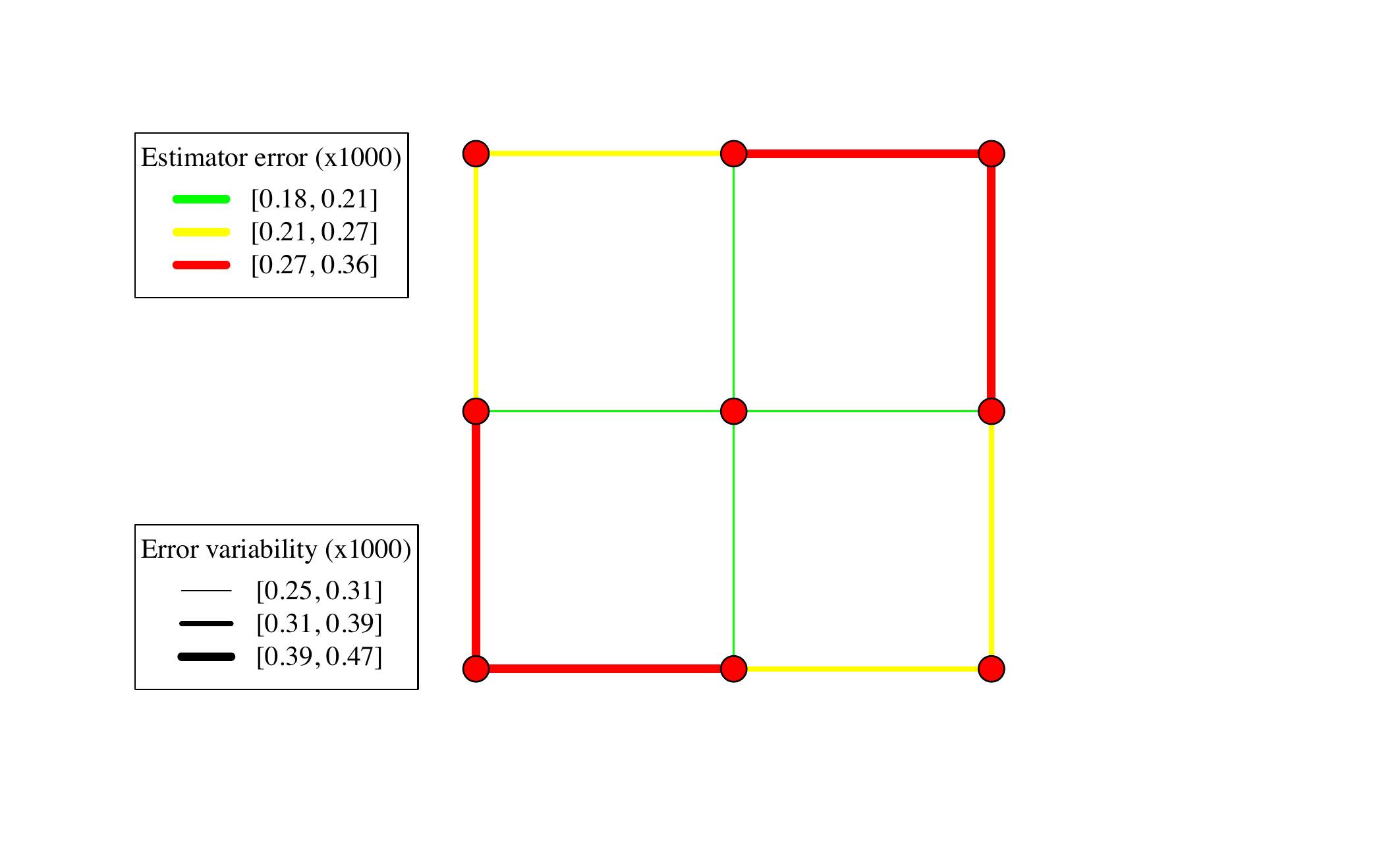}
		\caption{Illustration of the effect of a higher resolution parameter on the estimation errors corresponding to the original graph for trapezoidal-shaped velocity functions.} 
		\label{fig:Ex5}
	\end{center}
\end{figure}

In this example, adjacent edges in the higher resolution graph are, in principle, traversed with different expected velocities. This may introduce bias into the estimator such that it no longer performs well locally, as we have seen in Example 4. However, as also noticed in Example~4, when increasing the resolutions in $\BFr$  the differences in velocities between neighboring edges decrease, due to the fact that the velocity function is (almost everywhere) differentiable. Overall, this results in better estimates of the expected travel times. This effect is shown in Figure~\ref{fig:Ex5}.

In the left plots the estimation procedure is carried out for two choices of the resolution, namely $r_e=1$ (top) and $r_e=4$ (bottom). In the right hand side plots we see the approximation of the expected relative errors for the original graph $G$.
The plots reveal that the estimates for the original graph $G$ improve as the resolution of its corresponding high-resolution graph $G_\BFr$ increases.
(Note the different scales for the errors in each plot.)

\subsection*{Example 6: Non-normal data}
\begin{figure}[!ht]
	\begin{center}
		\includegraphics[width = 0.5\textwidth, trim={2cm 2cm 3cm 1cm}, clip]{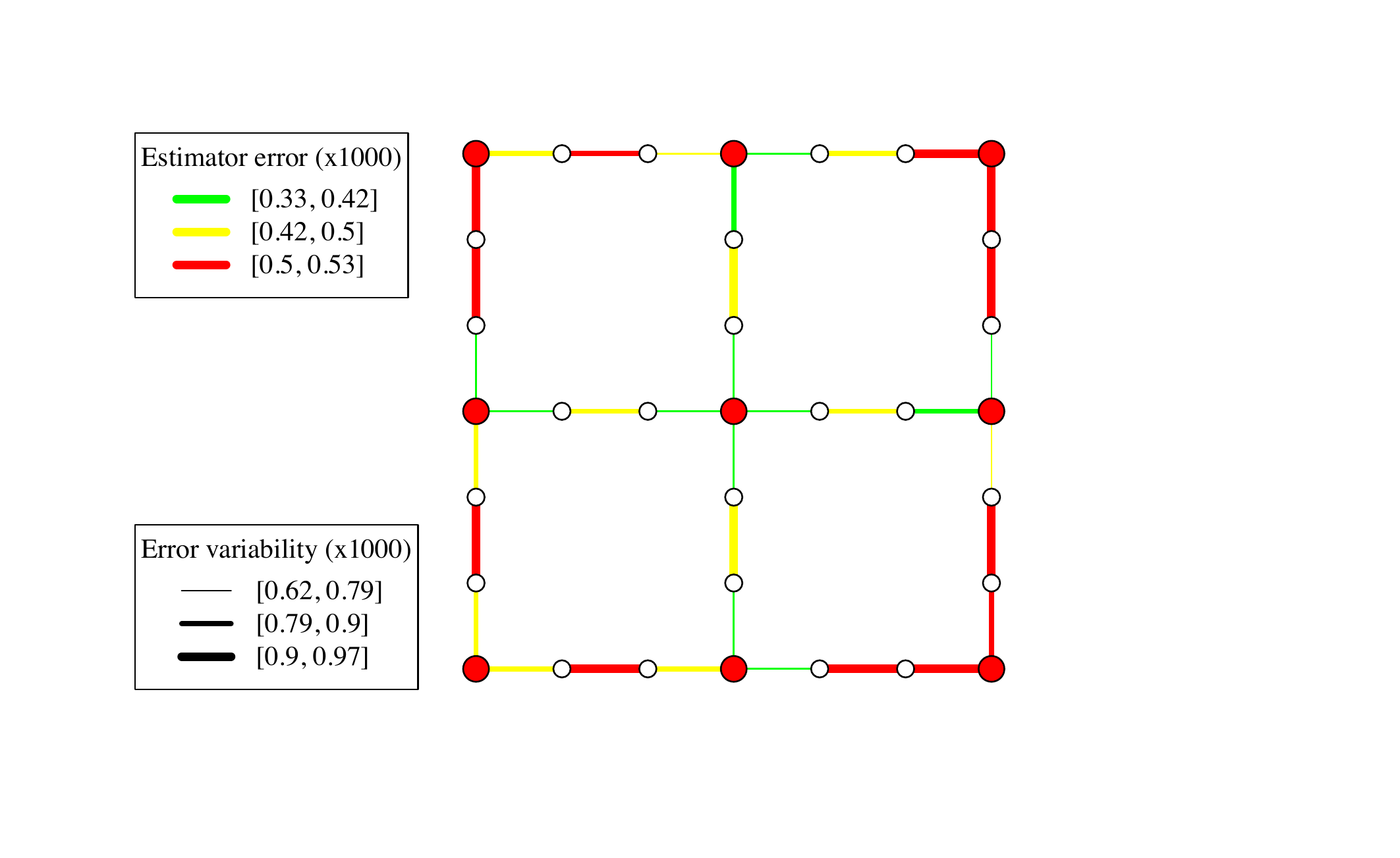}%
		\includegraphics[width = 0.5\textwidth, trim={2cm 2cm 3cm 1cm}, clip]{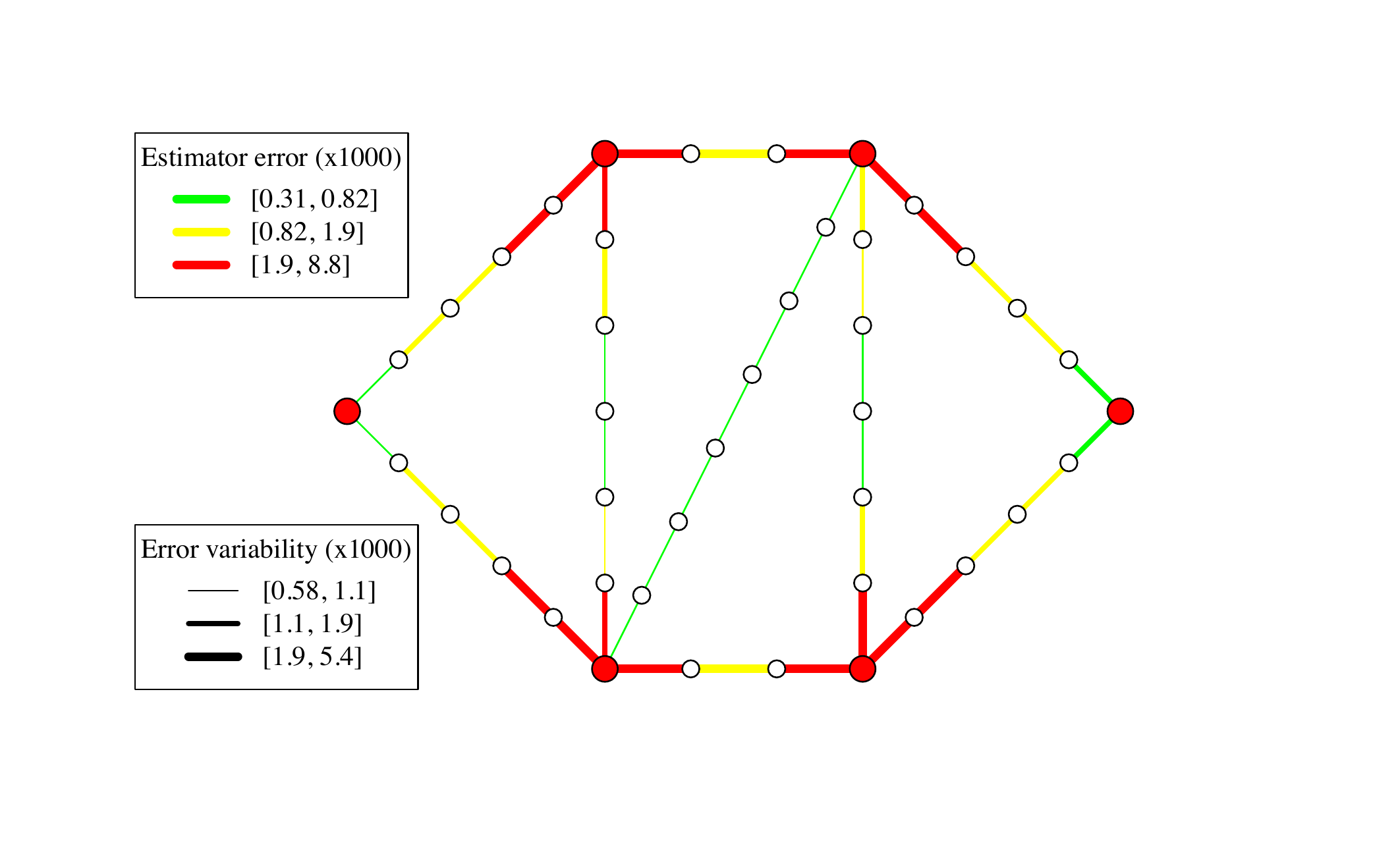}
		\caption{Estimation errors corresponding to the instances discussed in Example 1 but with non-normal data.} 
		\label{fig:Ex6}
	\end{center}
\end{figure}

Our estimator relies on the modeling assumption that the vector of averages $\BFX^{(\BFn)}$ is (approximately) normally distributed so that our normal posterior is a good proxy for the posterior distribution of $\BFmu$; see \eqref{eq:model}.
In all preceding examples we sampled the individual observations $X_{e,i}$ from a normal distribution so that this assumption is fulfilled by default. 
It is important to note that, bearing in mind that the entries of the data vectors are actually averages and appealing to the central limit theorem, this assumption is in practice by approximation fulfilled.
To demonstrate that our estimation procedure does not rely heavily on the normality assumption, in this example we simulate data from a different distribution.
Figure~\ref{fig:Ex6} shows the results.

More specifically, we have sampled the data from the gamma distribution with shape and rate parameters chosen so that the observations match the expectation and variance of the observations from Example 1. 
Comparing with Figure~\ref{fig:Ex1}, we see that the errors of the estimates effectively match with those obtained from normal data, thus corroborating the claim that the normality assumption does not play a crucial role. 

\bigskip
The preceding examples illustrate the performance and some properties of the estimation procedure from Section~\ref{sec:inference}. In particular, we have seen that the estimation procedure provides us with accurate estimates of the mean travel times and that the accuracy increases with the sample sizes and/or resolutions. We compared the estimates to those obtained from the estimation procedure that does not use smoothing, clearly revealing the beneficial effect of the smoothing parameter. Especially the estimates corresponding to edges with many neighbors with similar velocities benefit substantially from the smoothing parameter. In the case that neighboring edges do not have similar velocities, the smoothing parameter results in an increased bias at these edges. However, this can be mitigated by choosing a higher resolution. The resolution parameter also proves to be of great importance for non-constant expected velocity functions, as the differences in velocities between neighboring edges decrease for higher resolutions. Lastly, even though the estimator relies on a normality assumption, we have seen that we also obtain accurate estimates for non-normal data.

\section{Route selection examples}\label{sec:route}
Now that we have developed an estimator and illustrated its performance, we proceed by discussing their use in route selection. As pointed out in the introduction, a route is deemed optimal if it aligns with preferences of the individual driver. In this section, these preferences are expressed in terms of objective functions, which represent the routes' `disutilities'. We consider an elementary test network as well as a more sophisticated network. 

We wish to identify a path $\psi^*_{i,j}$ in a set $\Psi_{i,j}$ of feasible simple paths from vertex $i$ to vertex~$j$. Denoting the disutility of a route by the driver's objective function $f(\cdot)$, our goal is to find the best route among the feasible paths, namely
\begin{equation}
\label{eq:optim_statement}
\psi^*_{i,j} = 
\argmin_{\psi_{i,j}\in\Psi_{i,j}} f(\psi_{i,j}),
\end{equation}
assuming uniqueness of the minimizing path.

In this section, we consider several possible choices for $f$; namely:
\begin{itemize}
\item[$\circ$] expected travel time;
\item[$\circ$] quantile of the posterior expectation (for instance the $0.95$-quantile, i.e., the $95$\%-percentile);
\item[$\circ$] quantile of the distribution of the estimator of the expected travel time;
\item[$\circ$] sum of the squared difference of the expected velocities of the consecutive edges;
\item[$\circ$] mean of the squared difference of the expected velocities of the consecutive edges.
\end{itemize}
Evidently, some of these choices have more practical appeal than others, but our exposition also serves the goal of demonstrating the generality of our approach. 

Clearly, the value that the disutility $f(\psi_{i,j})$  takes for a given path $\psi_{i,j}$ depends on the distribution of the data, and is hence not known a priori to the driver, entailing  that the optimization problem described in~\eqref{eq:optim_statement} is not directly solvable.
We therefore express $f(\psi_{i,j})$ in terms of parameters that we can estimate, thus yielding an optimization problem that we \emph{can} solve:
\begin{equation}
\label{eq:optim_statement_estim}
\hat\psi^*_{i,j} = 
\argmin_{\psi_{i,j}\in\Psi_{i,j}} \hat f(\psi_{i,j}).
\end{equation}
We think of $\hat\psi^*_{i,j}$ as a proxy for $\psi^*_{i,j}$.

In the following examples we illustrate solving the optimization problem in~\eqref{eq:optim_statement_estim} for the disutilities mentioned above.

\medskip

\begin{figure}[!ht]
	\begin{center}
		\includegraphics[width = 0.5\textwidth, trim={5cm 2cm 4cm 2cm}, clip]{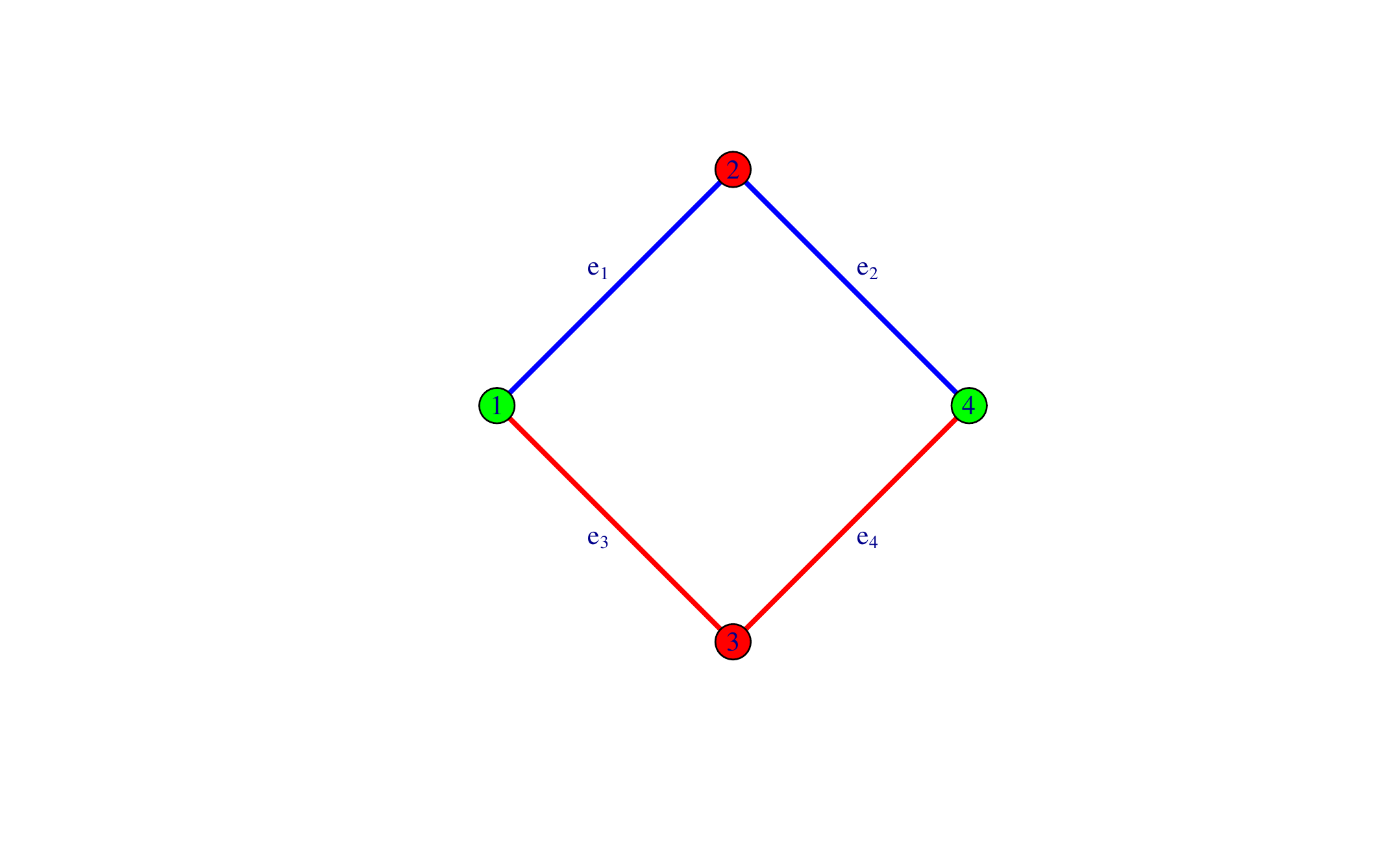}
		\caption{The graph of the elementary traffic network. There are only two possible routes from the origin, vertex 1, to the destination, vertex 4.} 
		\label{fig:opt1}
	\end{center}
\end{figure}

We will consider an elementary traffic network with a graph as in Figure~\ref{fig:opt1}, where each of the four edges has a length of 1 kilometer. Suppose a driver wants to travel from vertex~$1$ to vertex $4$. Clearly, $\Psi_{1,4}$ consists of only two possible paths, viz.\ the red route and the blue route: $\Psi_{1,4}=\{(e_1,e_2),(e_3,e_4)\}$. The optimal route is determined by the output of the estimation procedure and the driver's objective function $f(\cdot)$.

The following paragraphs show the effect of making different choices for the disutility $f$ on the route that gets selected.
The graph in Figure~\ref{fig:opt1} is rather simple so that we obtain more direct insight into the effect of objective functions on the selected route.

\subsection*{Expected travel time}
In case the driver wishes to minimize her expected travel time, the optimal route clearly satisfies
\[
	\psi^*_{1,4}=
	\begin{cases}
		(e_1,e_2) \quad \text{ if } {\mu}_{e_1}+{\mu}_{e_2}<{\mu}_{e_3}+{\mu}_{e_4};\\
		(e_3,e_4) \quad \text{ otherwise.}
	\end{cases}
\]
Below we refer to $(e_1,e_2)$ as Route 1, and to $(e_3,e_4)$ as Route 2.

\begin{figure}[!ht]
\begin{center}
\includegraphics[height=3in]{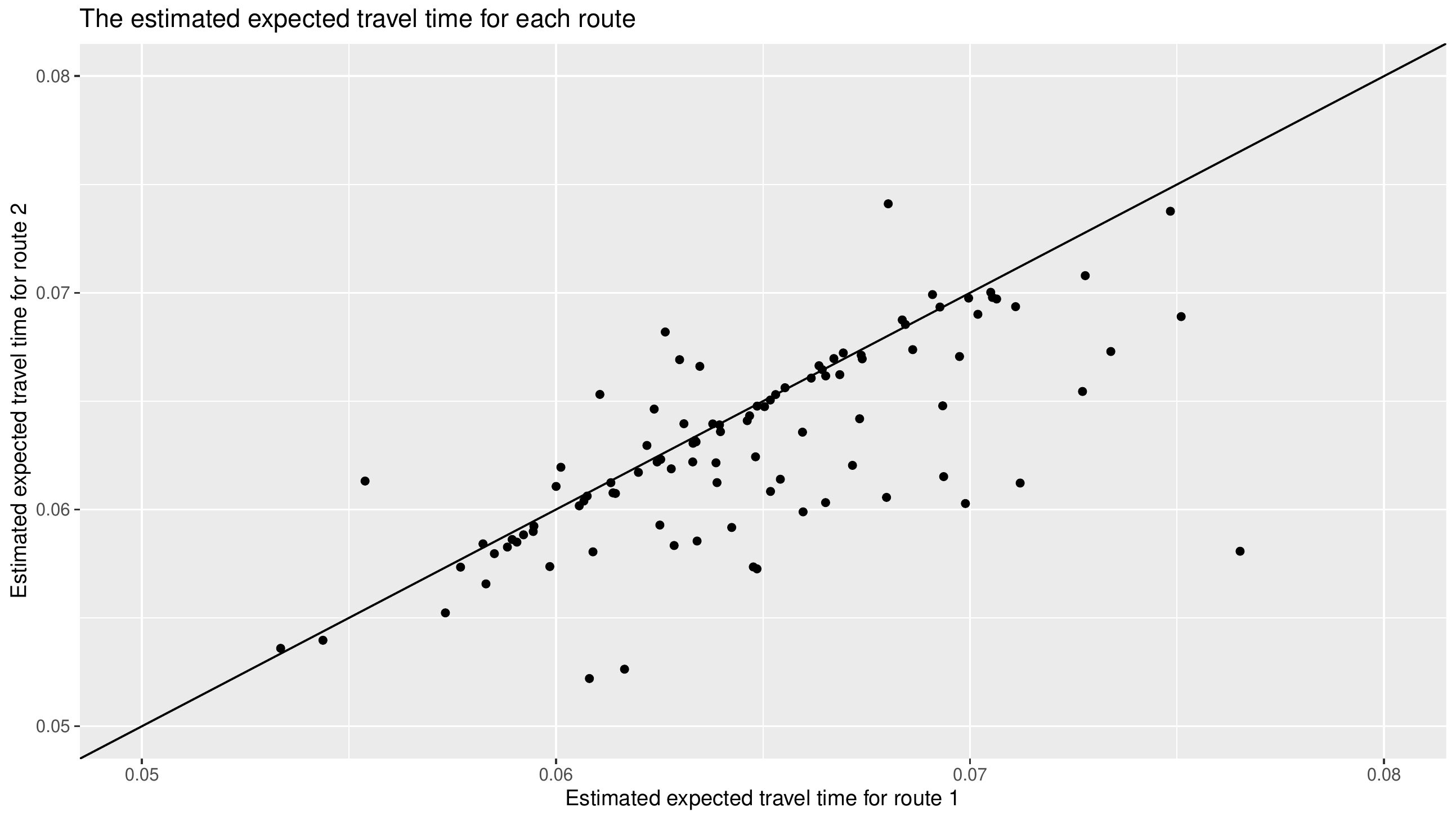}
\caption{The estimated expected travel time for both routes for each of the 100 simulations. Route 1 corresponds to path $(e_1,e_2)$ and Route 2 to path $(e_3,e_4)$. The diagonal line is added to identify the best route according to each simulation more easily; points below this line correspond to the simulations for which Route 2 has the lowest estimated expected travel time.} 
\label{fig:scatter1}
\end{center}
\end{figure}

In this numerical experiment, we sample from a gamma distribution such that $\mathbb{E}X_{e_j}=120$ and $\mathbb{V} X_{e_j}=36^2$ for $j=1,2$, whereas $\mathbb{E}X_{e_j} \approx 109.1$ 
and $\mathbb{V} X_{e_j}=36^2$ for $j=3,4$. These expectations correspond to the travel time in seconds for edges with expected velocities of 30~km/h (for edges $e_1$ and $e_2$) and 33~km/h (for edges $e_3$ and $e_4$).
For each edge we sampled~10 observations so that $n_e=10$ for each $e\in E$. 
We performed the estimation procedure 100 times and plotted the value of the objective function for the pair of routes in each run; the results can be seen in Figure~\ref{fig:scatter1}.

Each point in Figure~\ref{fig:scatter1} represents the estimated expected travel time for the two routes.
We see that Route 2 minimizes the objective function in 76 of the 100 experiments, which is consistent with the fact that Route 2 indeed has a lower expected travel time.
This means that if we were to use our procedure to select a route, then the majority of the time the correct route would be selected.
Importantly, in this example the sample sizes are just $n_e=10$, and that increasing these sample sizes per edge (or having a smaller variance for the data at each edge) would result in Route 2 being selected even more often.

\subsection*{Quantile of the posterior expectation}
\begin{figure}[!ht]
	\begin{center}
		\includegraphics[height=3in]{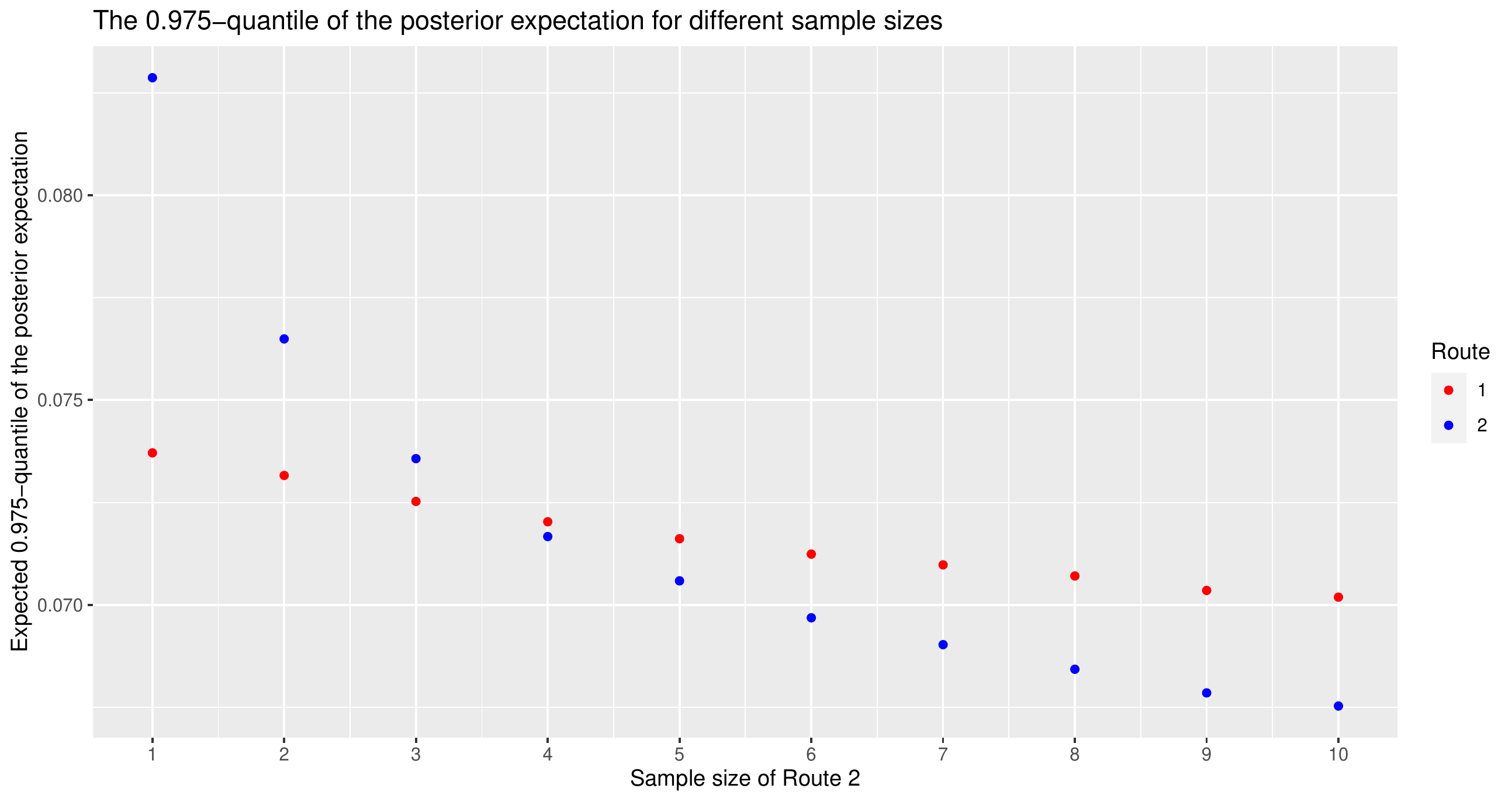}
		\caption{The expected objective value for both routes and for different sample sizes for Route 2, while keeping the number of observations for Route 1 constant at 10.} 
		\label{fig:mean}
	\end{center}
\end{figure}

The objective function in the previous example relies solely on (estimates of) the expected travel times. Our estimator for $\BFmu$ in \eqref{eq:posterior_mu} also provides us with an uncertainty quantification of these estimates. This fact is particularly useful in a situation in which drivers are reluctant to traverse routes for which the estimated mean travel times are low but carry large uncertainty. In this situation it may be more appropriate to minimize a certain (relatively high) quantile of the posterior distribution of the expected travel time of each of the routes. Specifically, in our numerics we choose the $0.975$-quantile of the posterior distribution of the travel time for a route as the objective function. 
In Figure~\ref{fig:mean} we report the average utility obtained across 10\,000 replications (which is an approximation of the expected utility) when the sample sizes for the edges in Route 1 are fixed at $n_{e_1}=n_{e_2}=10$, and try different sample sizes for the edges in Route 2.

From Figure~\ref{fig:mean} we can see that, even though the sample sizes for the edges that make up Route 1 are kept fixed, the estimates for the $0.975$-quantile of the posterior distribution of the expected travel time for Route 1 are changing as we change the number of observations collected at each of the edges that make up Route 2.
This again illustrates how the statistical procedure borrows information from edges that are close the each other.

The uncertainty of the estimated expected travel time of an edge heavily relies on the number of observations for that edge. Therefore, the routing criterion that minimizes a quantile of the posterior distribution of the expected travel time favors routes that are well-explored. We have already seen that Route 2 is most likely to minimize the expected travel time and is therefore expected to also minimize a quantile of the posterior mean in case of equal sample sizes, since the variances of the observations are constant across the graph. However, this is not necessarily the case if Route 2 is not as well-explored as Route 1, as illustrated by Figure~\ref{fig:mean}.

We see that if we have fewer than 4 observations for the edges that make up Route 2, while still having 10 observations for the edges that make up Route 1, the objective function is (on average) minimized by Route 1.
Route 2 may very well have a lower expected travel time, but a user who cares to optimize a quantile of the posterior distribution for the expected travel time for their route may still prefer Route 1 if selecting this route carries less uncertainty.

\subsection*{Quantile of the distribution of the estimator of the expected travel time}
Instead of just focusing on expected travel times only, drivers may also want to incorporate the variance of the travel time into their decision criterion. A higher quantile indicates that a driver is more reluctant to traverse edges with a high travel time variance; we say that this driver is more risk-averse. Conversely, if we consider the $0.5$-quantile, this objective function essentially reduces to the objective function of the shortest expected travel time (at least in a setting in which the median and mean are close, which will be the case in the central limit type of regime that we consider).

We again sample from a gamma distribution such that $\mathbb{E}X_{e_j}=120$ and $\mathbb{V} X_{e_j}=36^2$ for $j=1,2$, but this time let $\mathbb{E}X_{e_j}\approx 109.1$ and $\mathbb{V} X_{e_j}=72^2$ for $j=3,4$. Now, Route 2 has smaller expected travel time but higher variance than Route 1. Again, for each edge we obtain 10 observations so that $n_e=10$ for each $e\in E$. 

The risk-averseness of the driver will determine which route is preferred.

\begin{figure}[!ht]
	\begin{center}
		\includegraphics[width = 0.33\textwidth]{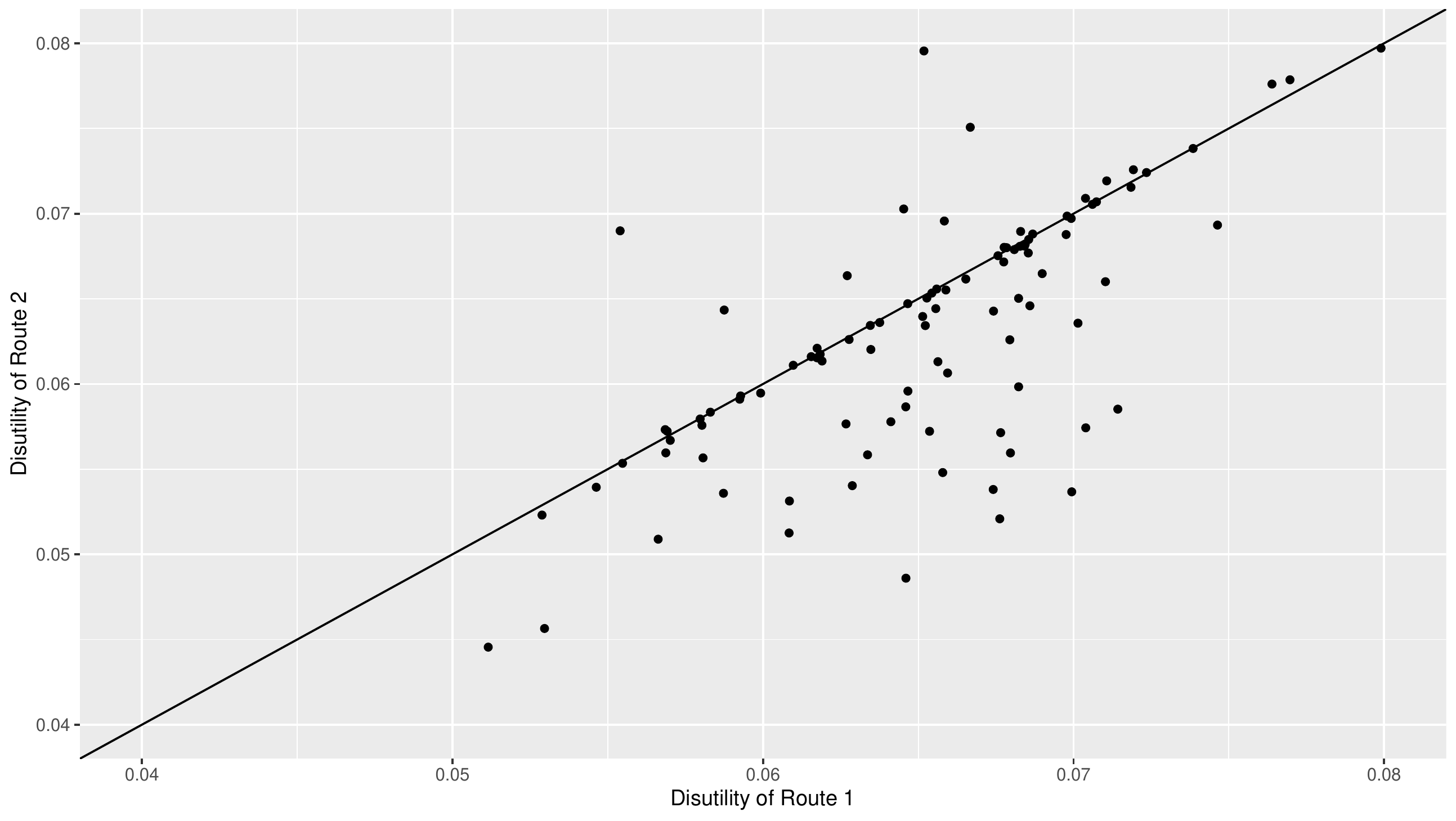}%
		\includegraphics[width = 0.33\textwidth]{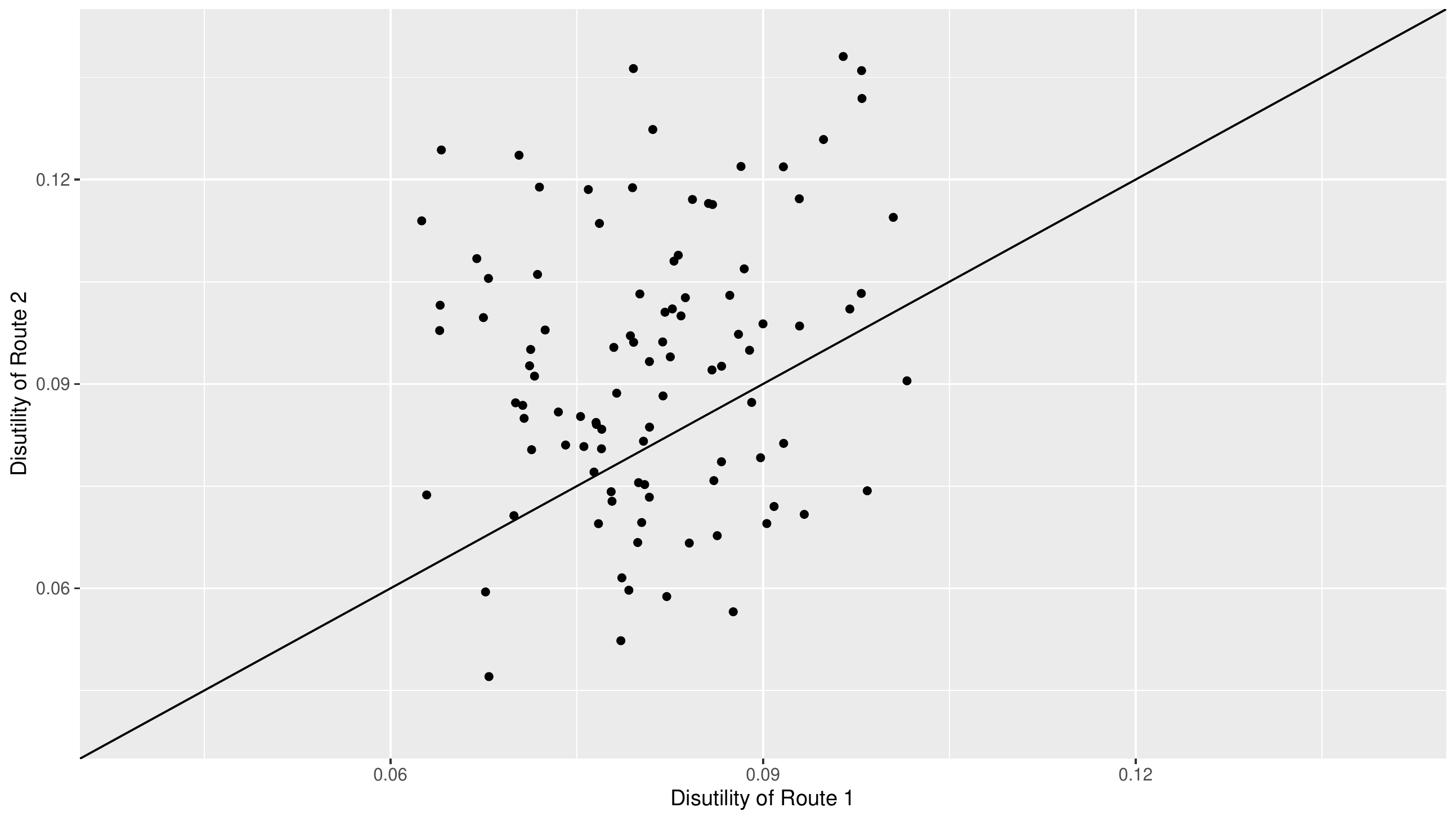}%
		\includegraphics[width = 0.33\textwidth]{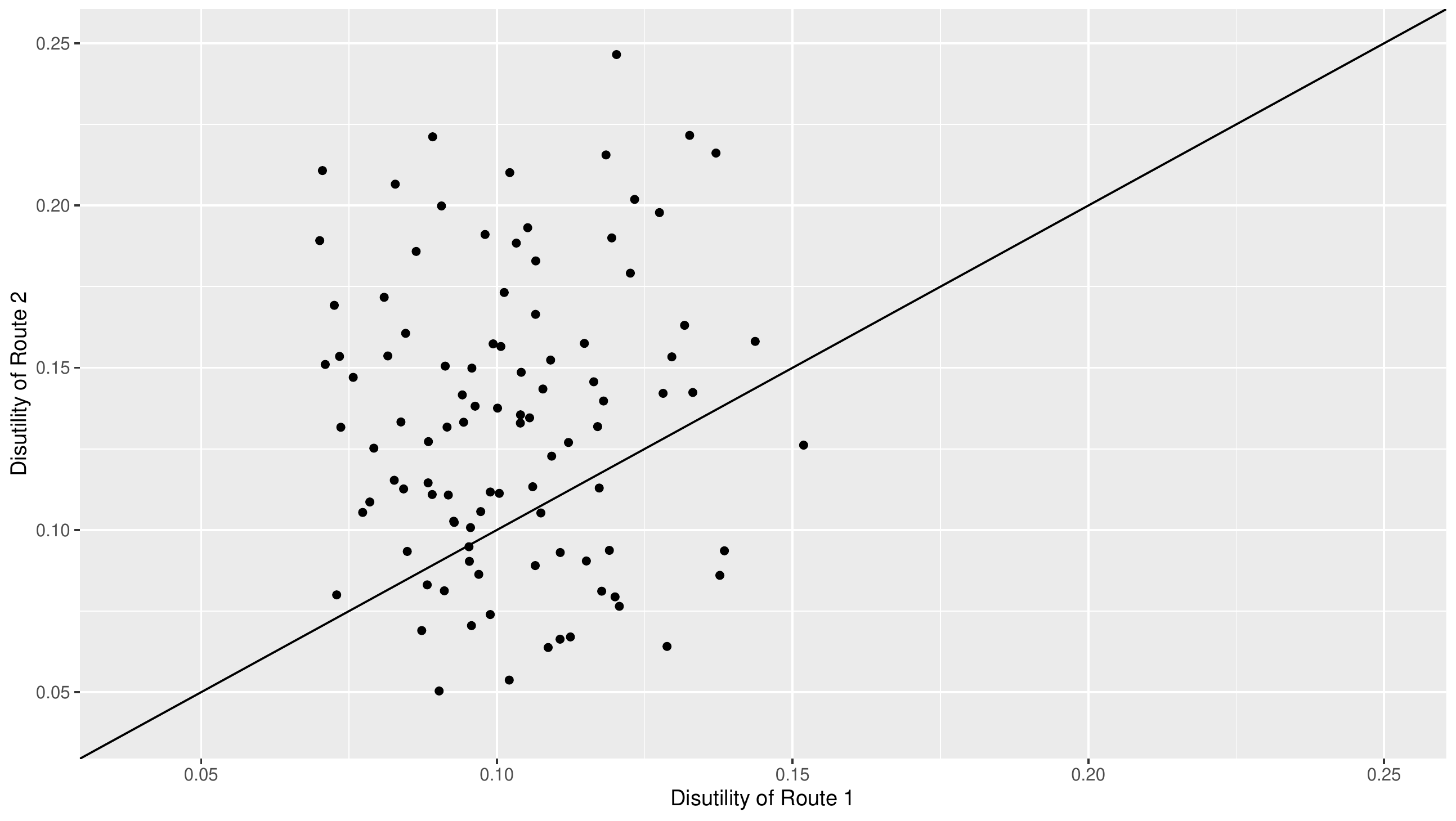}
		\caption{The objective function for the 0.5- (left), 0.8- (middle), and 0.975- (right) quantile of the distribution of the estimator of the expected travel time.} 
		\label{fig:quantiles}
	\end{center}
\end{figure}

In Figure~\ref{fig:quantiles}, we plotted 100 realizations of the objective function of both routes for the different quantiles of the distribution of the estimator of the expected travel time. The quantiles are  easily computed using the fact that the estimator of the travel time of a route is approximately normal. 

We see that for the median (left-most plot), Route 2 is often (correctly) selected as the optimal route. This quantile corresponds to drivers that are risk-neutral and therefore prefer the route that is expected to be faster without regarding the variability of the travel time for that route. As we consider higher quantiles, Route 1 becomes more attractive for increasingly risk-averse drivers. This route has 
a travel time with 
only a slightly higher mean, but it has a smaller variance.

\subsection*{Route selection in a larger traffic network}

\begin{figure}[!ht]
	\begin{center}
		\includegraphics[width = 0.5\textwidth, trim={4cm 1cm 4cm 2cm}, clip]{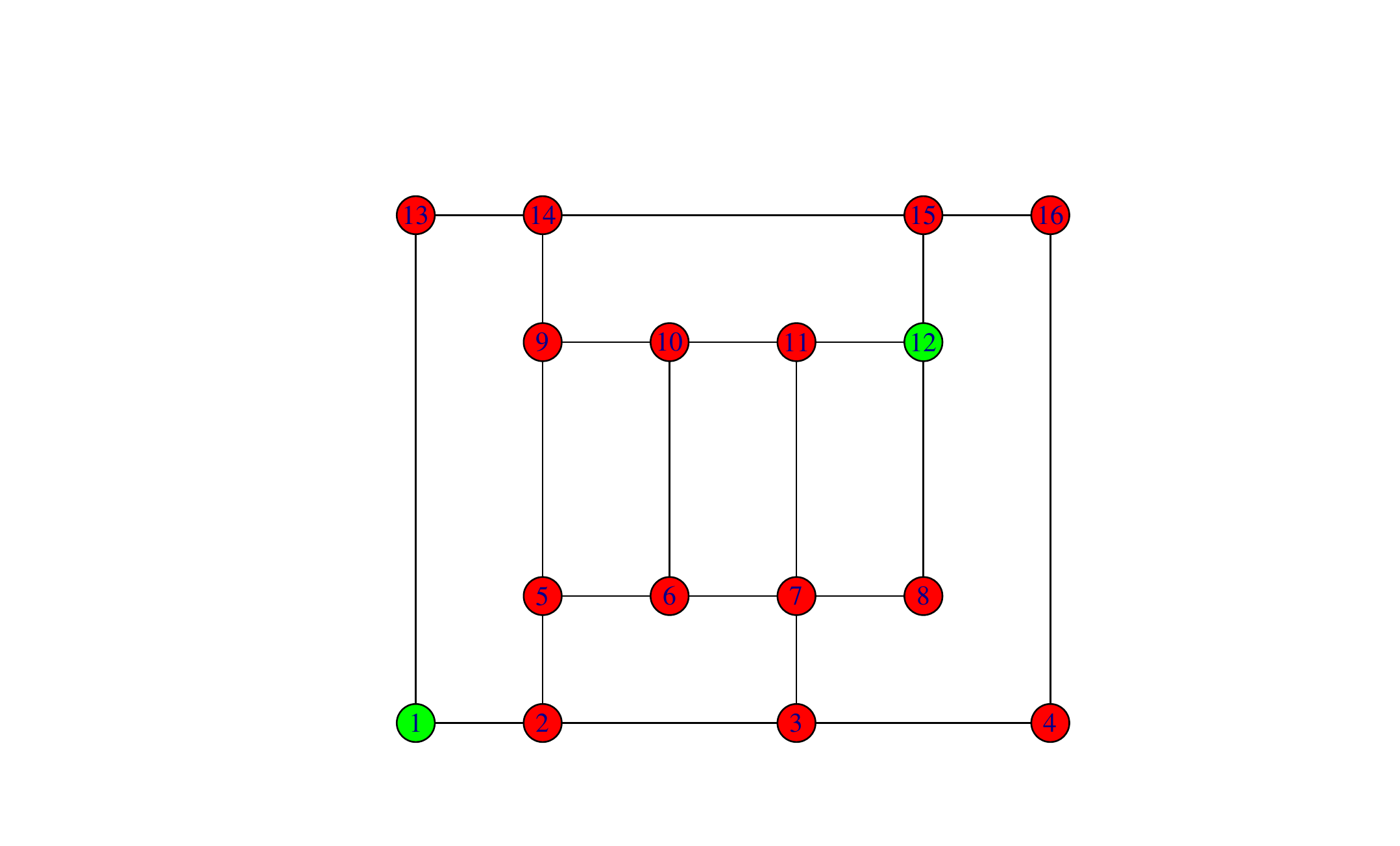}
		\caption{The graph of a larger traffic network. There are many possible routes from the origin, vertex 1, to the destination, vertex 12.} 
		\label{fig:big_network1}
	\end{center}
\end{figure}

We used the elementary traffic network whose graph is depicted in Figure~\ref{fig:opt1}, with just two routes, to discuss the effect the different disutilities have on route selection.
We proceed by studying the larger network depicted in Figure~\ref{fig:big_network1}.
This example uses all disutilities $f$ given in the list at the beginning of this section.

For this traffic model, we let particles traverse the outer edges of the graph at a constant expected velocity of 60 km/h. For the inner edges, we assume the trapezoid shaped expected velocities as in Figure~\ref{fig:Ex5velocity}. Moreover, the intersection at vertex 7 has the property that drivers are not required to decelerate when approaching this intersection (i.e., they keep driving at 40 km/h).
We also assume that the travel times have a standard deviation of $72$ seconds per kilometer (i.e., 0.02 hour per kilometer) for each edge, except for the edges $(2,3)$, $(3,4)$, $(4,16)$, $(12,15)$, and $(15,16)$, which have a standard deviation of $36$ seconds per kilometer (i.e., 0.01 hour per kilometer).
As for sample sizes, at the \emph{outer edges} -- i.e., 
$(1,2)$, $(1,13)$, $(2,3)$, $(3,4)$, $(4,16)$, $(13,14)$, $(14,15)$, and $(15,16)$ --
we collected $100$ observations, while at the remaining edges, the \emph{inner edges}, we collected $10$ observations.
The scale of the figure is such that the length of the edge $(1,2)$ corresponds to 1 kilometer.

\begin{figure}[!ht]
	\begin{center}
		\includegraphics[width=0.5\textwidth, trim=1cm 0cm 4cm 1cm, clip]{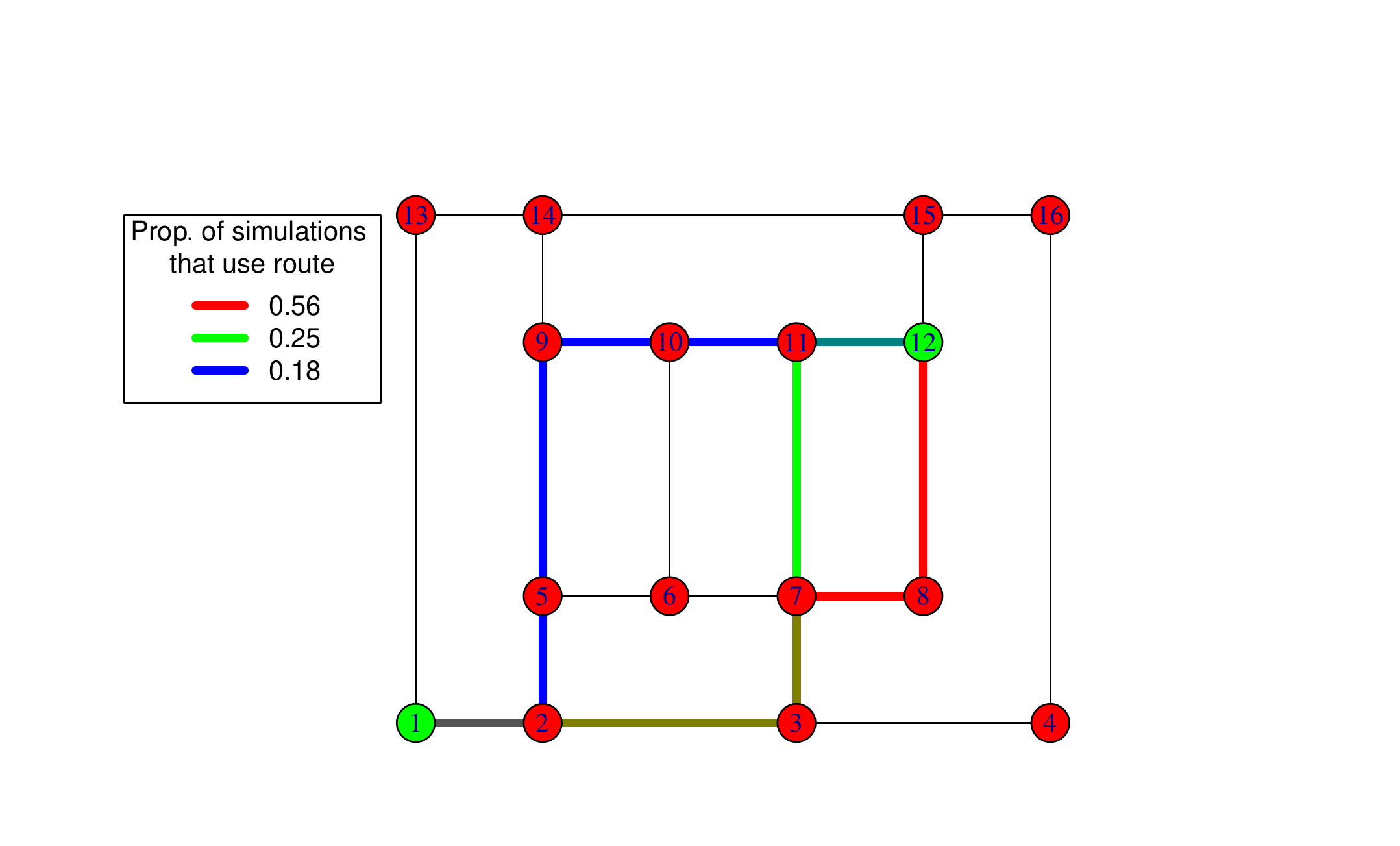}%
		\includegraphics[width=0.5\textwidth, trim=1cm 0cm 4cm 1cm, clip]{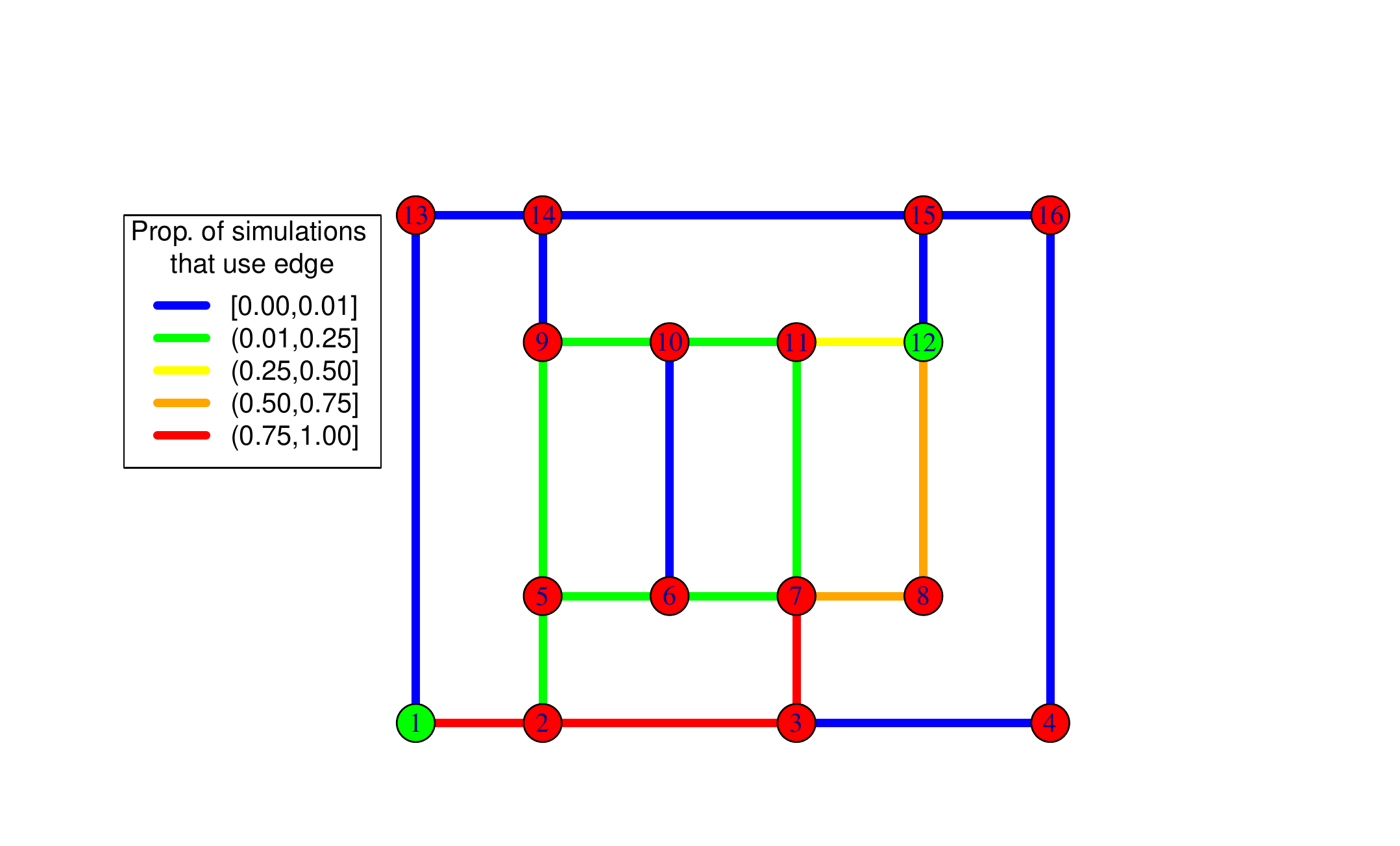}
		\caption{Using the expected travel time as objective function. The left picture shows the three routes that minimize the objective function most often, whereas the right picture shows the frequency with which each edge is part of the optimal route.}
		\label{fig:util1}
	\end{center}
\end{figure}

We wish to find the routes from vertex $1$ to vertex $12$ that minimize each of the objective functions listed above. We performed $M=$10\,000 simulations and identified the route that minimizes the objective function for each simulation and for each objective function. 
For each objective function we report a figure with two plots; see, e.g., Figure~\ref{fig:util1} for reference.
In the left plot we report the three routes that were most frequently selected as the optimal route across the $M$ simulations.
These are colored red (most often selected), green (second most often selected), and blue (third most often selected); in the legend we report the fraction of the $M$ simulations in which each route was selected as optimal route.
Note that some edges in such plots are part of multiple routes, so we color an edge olive when the red and green routes overlap, purple when red and blue overlap, teal when green and blue overlap, and gray when all three routes overlap.
In the right plot we report a heat-map that depicts the frequency with which each edge is part of the optimal route.
This can be thought of as indicative of the congestion of the network, if all vehicles were to select which route to take based on the same objective function.




\medskip

First, we perform this simulation with the objective function set to the expected travel time; cf.~Figure~\ref{fig:util1} for the results. 
The right plot of Figure~\ref{fig:util1}  shows that the edges $(1,2)$ and $(2,3)$ are the ones that feature more often in the optimal route. This is not surprising, as the expected velocities of these edges are relatively high, while these edges are also part of the routes with the shortest travel distance. On the other hand, the remaining outer edges are (almost) never used. Indeed, in this particular network, the routes consisting of these remaining outer edges have longer expected travel times due to the longer travel distance, and the estimator picks up on this successfully. This is in line with the left plot of Figure~\ref{fig:util1}, where we see the routes with the smallest expected travel time. Observe that the two routes that minimize the expected travel time most frequently, utilize vertex 7. Recall that the intersection at this vertex does not require the drivers to decelerate, which indeed contributes to a lower expected travel time. 

\begin{figure}[!hb]
	\begin{center}
		\includegraphics[width=0.5\textwidth, trim=1cm 0cm 4cm 1cm, clip]{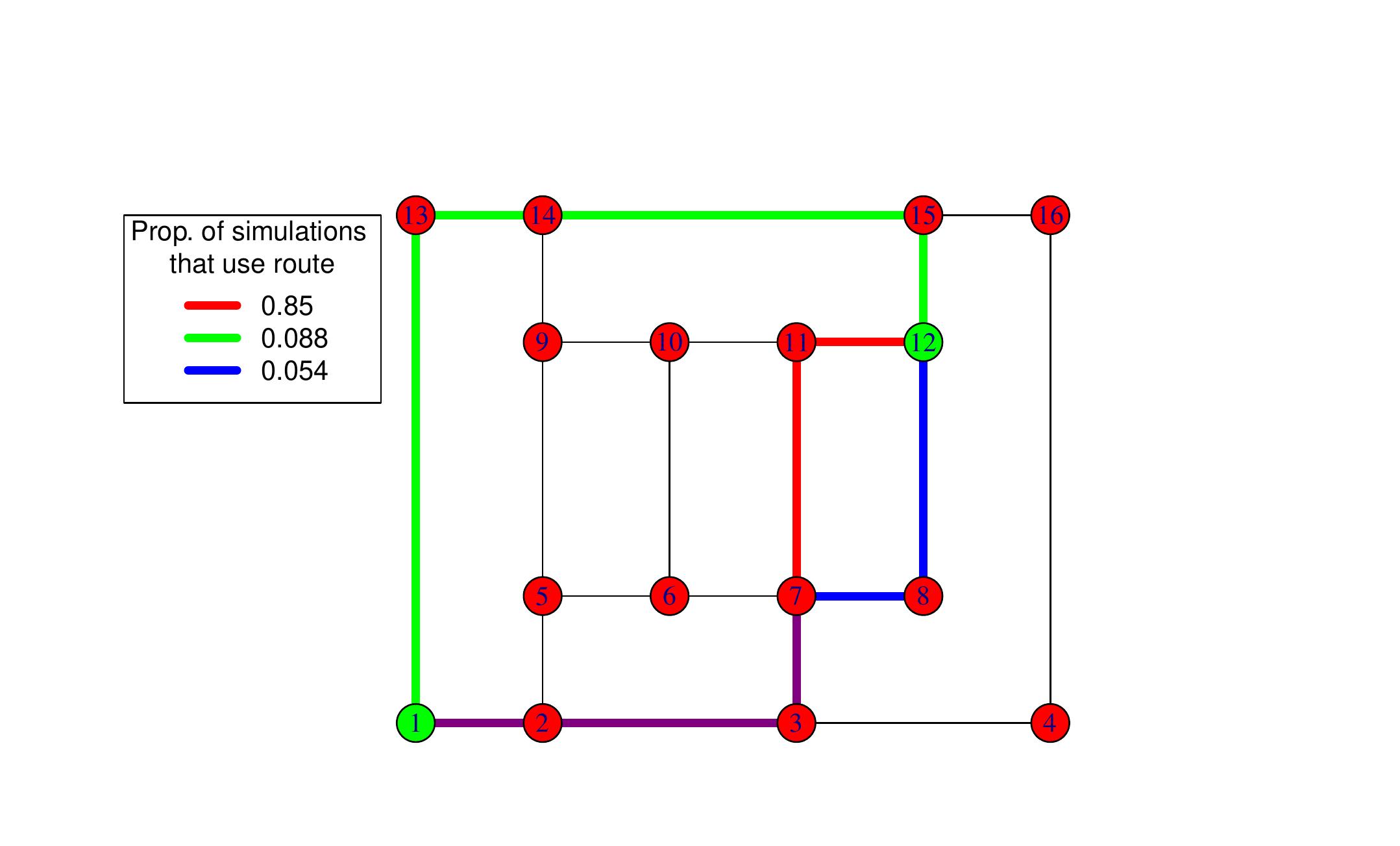}%
		\includegraphics[width=0.5\textwidth, trim=1cm 0cm 4cm 1cm, clip]{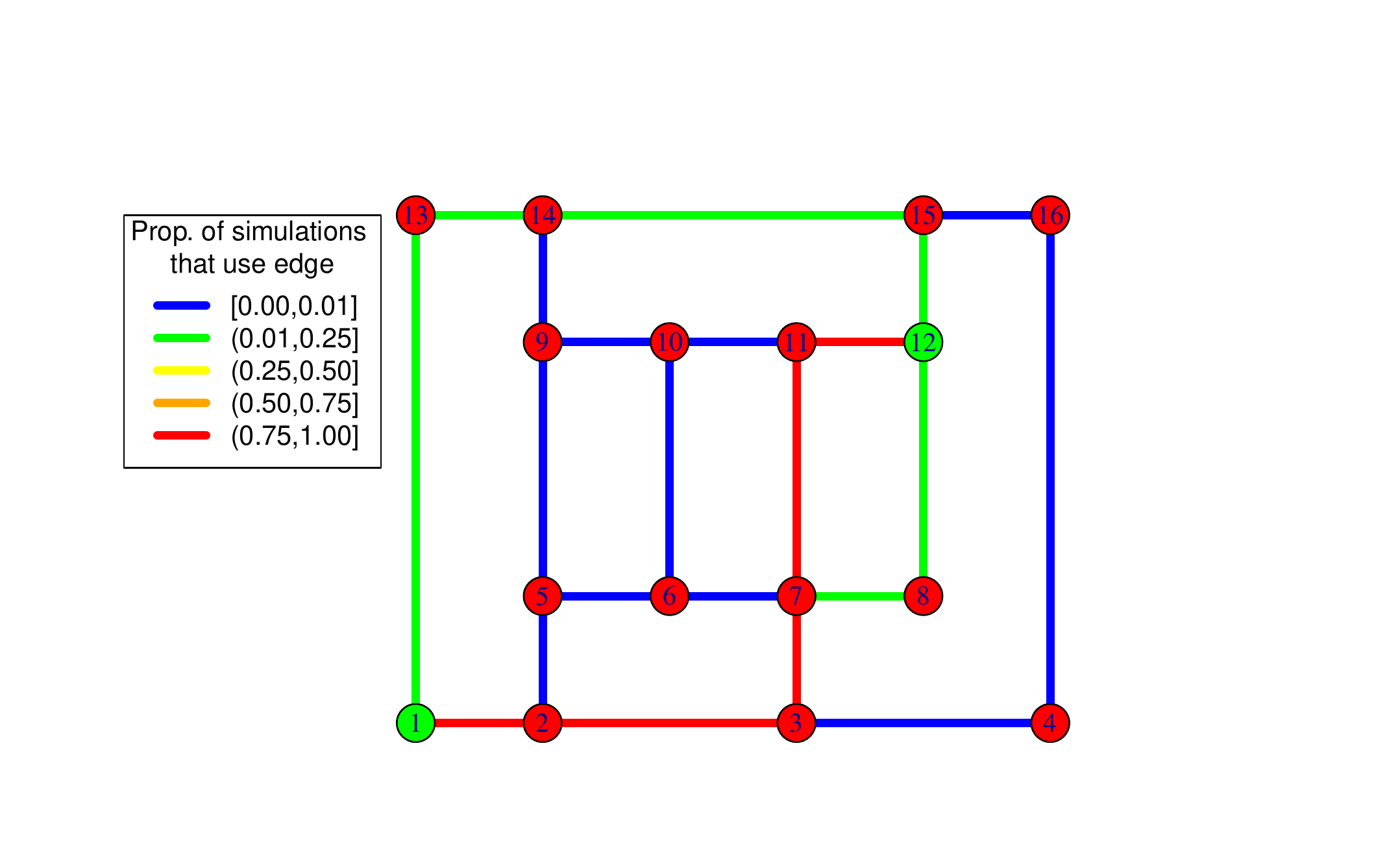}
		\caption{Using the 0.975-quantile of the posterior distribution as objective function. The left picture shows the three routes that minimize the objective function most often, whereas the right picture shows the frequency with which each edge is part of the optimal route.}
		\label{fig:util2}
	\end{center}
\end{figure}

Next, in Figure~\ref{fig:util2}, we consider the objective to minimize the 0.975-quantile of the posterior expectation. Comparing Figure~\ref{fig:util2} (right) to Figure~\ref{fig:util1} (right), we see that the outer edges $(1,13)$, $(13,14)$, and $(14,15)$ are more frequently part of the optimal route when we minimize the $0.975$-quantile of the posterior expectation rather than just the expected travel time. We have seen earlier in this section that the sample sizes for the edges play a prominent role here. As mentioned before, in our simulation we generated $100$ observations for each outer edge of the graph, whereas we only generated $10$ observations for the other edges. In other words, traversing the outer edges carries less uncertainty and, consequently, the $0.975$-quantile of the marginal posterior expectation of the routes that utilize the outer edges will be more concentrated around their expectation, resulting in smaller $0.975$-quantiles. In the left plot of Figure~\ref{fig:util2} we see that the route that utilizes the outer edges $(1,13)$, $(13,14)$, and $(14,15)$ is the second most selected route when minimizing the $0.975$-quantile of the posterior expectation. Note that this route was not yet visible in Figure~\ref{fig:util1} (left), indicating that the lower uncertainty of this route compensates for the somewhat higher expected travel time. The other routes in Figure~\ref{fig:util2} (left), however, coincide with routes in Figure~\ref{fig:util1} (left); observe that a lower expected travel time also contributes to a lower $0.975$-quantile of the marginal posterior expectation.

\begin{figure}[!ht]
	\begin{center}
		\includegraphics[width=0.5\textwidth, trim=1cm 0cm 4cm 1cm, clip]{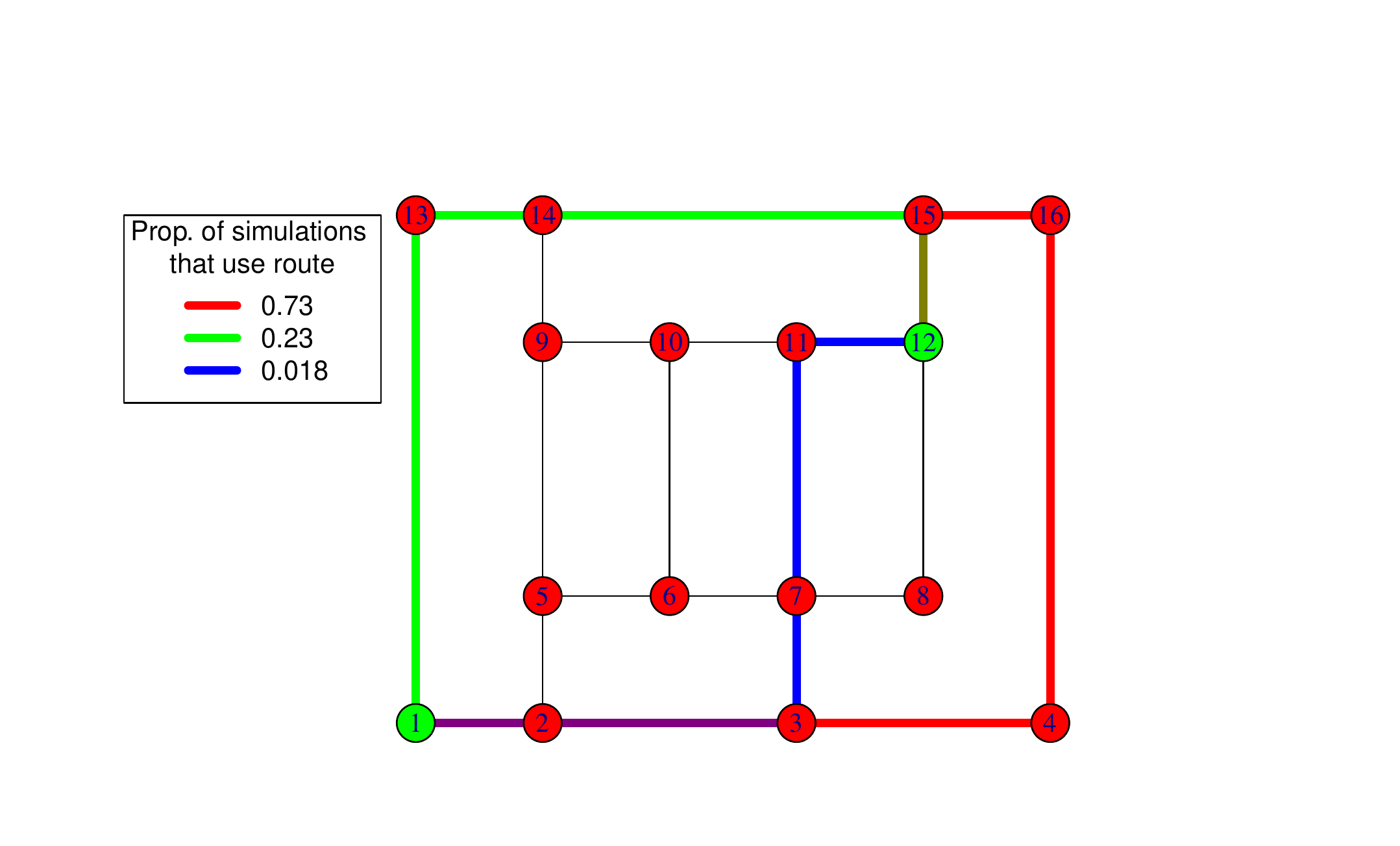}%
		\includegraphics[width=0.5\textwidth, trim=1cm 0cm 4cm 1cm, clip]{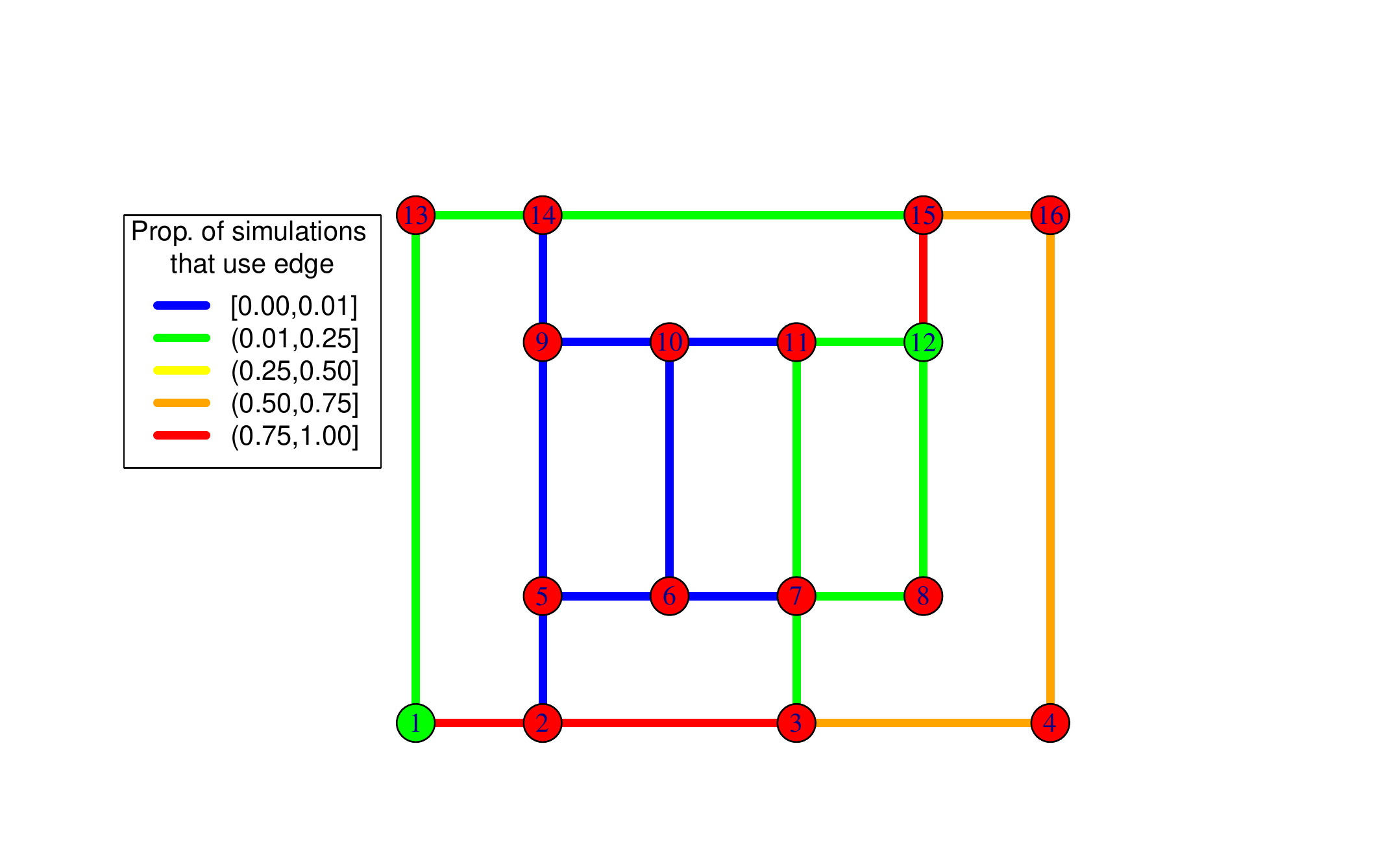}
		\caption{Using the 0.975-quantile of the distribution of the estimator of the expected travel time as objective function. The left picture shows the three routes that minimize the objective function most often, whereas the right picture shows the frequency with which each edge is part of the optimal route.}
		\label{fig:util3}
	\end{center}
\end{figure}

Another quantity that we may want to minimize is the $0.975$-quantile of the distribution of the estimator of the expected travel time, see Figure~\ref{fig:util3}. Now, the variance of the travel times of the edges has become an important factor. Recall that we assumed that the travel times have a standard deviation of $72$ seconds per kilometer for each edge, except for edges $(2,3)$, $(3,4)$, $(4,16)$, $(12,15)$, and $(15, 16)$, which have a standard deviation of $36$ seconds per kilometer. Therefore, if the objective is to minimize the $0.975$-quantile, the routes consisting of edges with a lower variance now also become attractive for the minimization. In the left plot of Figure~\ref{fig:util3} we indeed see that the route that utilizes the previously mentioned edges is now most frequently minimizing the $0.975$-quantile, while this route was never optimal for the previously studied utilities. This indicates that the lower variance compensates the somewhat higher travel time. 

\begin{figure}[!ht]
	\begin{center}
		\includegraphics[width=0.5\textwidth, trim=1cm 0cm 4cm 1cm, clip]{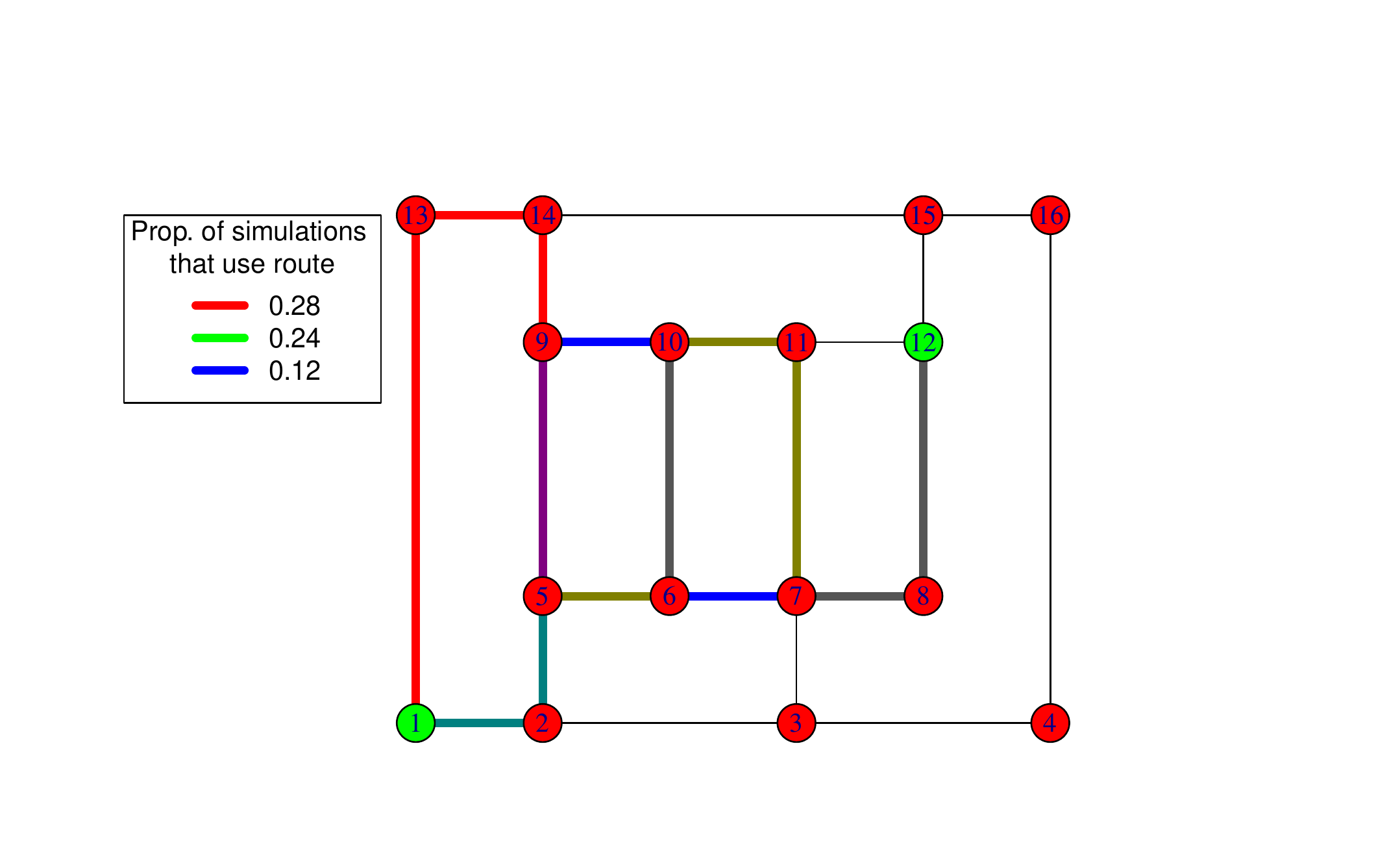}%
		\includegraphics[width=0.5\textwidth, trim=1cm 0cm 4cm 1cm, clip]{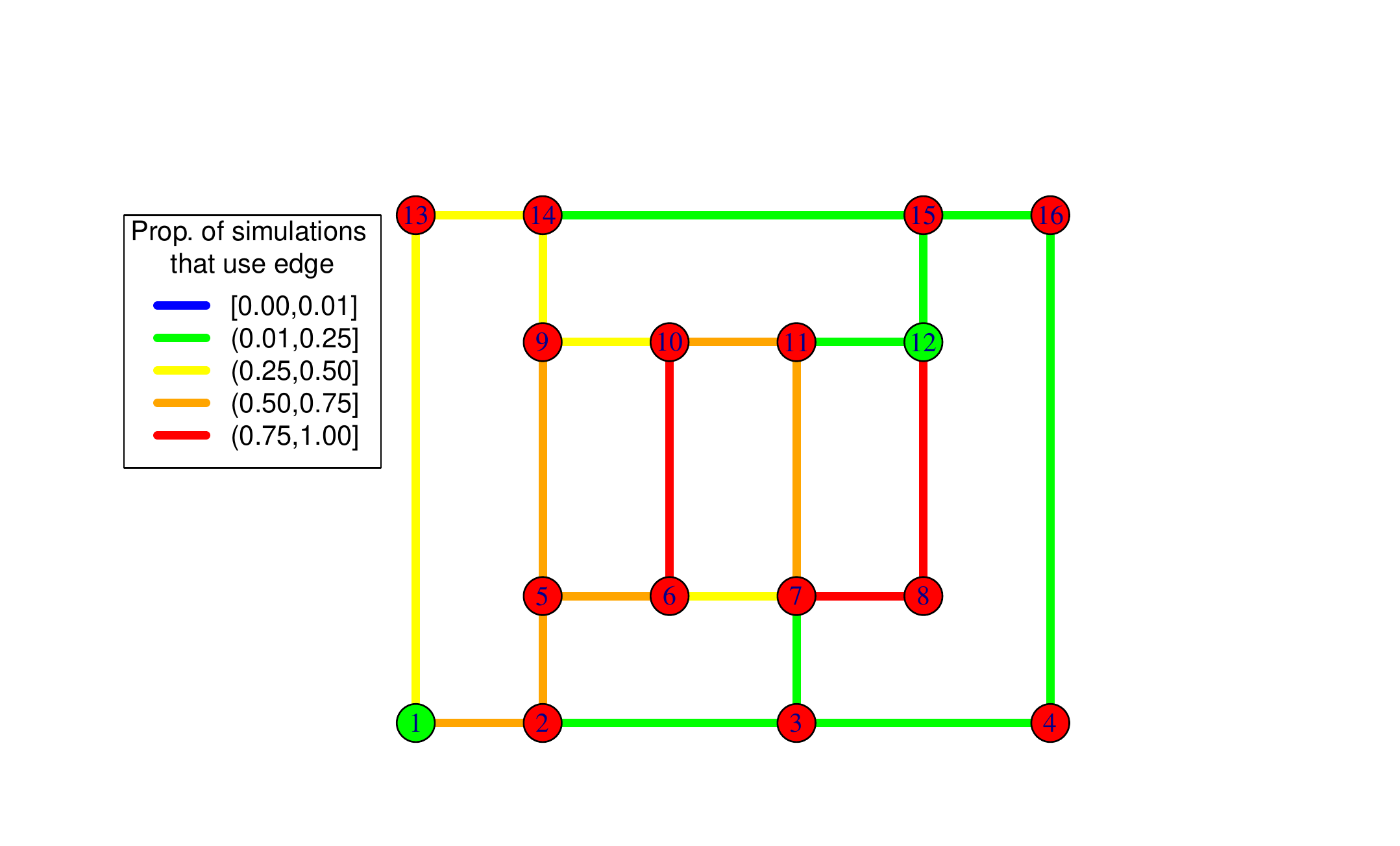}
		\caption{Using the mean of the squared difference of the expected travel time of the consecutive edges as objective function. The left picture shows the three routes that minimize the objective function most often, whereas the right picture shows the frequency with which each edge is part of the optimal route.}
		\label{fig:util5}
	\end{center}
\end{figure}

The last two objective functions whose use we illustrate here are respectively the mean- and the sum of the squared difference of the expected travel time at consecutive edges of a path.
These squared differences quantify the variation of speed across two neighboring road segments.
Minimizing the mean of the squared difference of the expected travel times at consecutive edges on a path effectively amounts to keeping the velocity as constant as possible.
Figure~\ref{fig:util5} illustrates the results corresponding to this utility. 
From the left plot we see that the driver favors routes that make use of the inner part of the network.
Judging by how often each of these routes is selected, there is no clearly preferred route.
Because using the mean of the squared differences does not penalize the number of road segments that are used, we see that no optimal routes (among the top three) continue from vertex 10 to vertex 12 directly, but rather take longer routes (that apparently minimize speed variability) to reach the destination vertex.
Looking at the right plot in Figure~\ref{fig:util5}, we see that the heat-map is less concentrated than for other objective functions.
Again, this is not surprising since this objective function is \emph{distance-indifferent}, so that there is no concentration around routes with lower expected travel time.

\begin{figure}[!ht]
	\begin{center}
		\includegraphics[width=0.5\textwidth, trim=1cm 0cm 4cm 1cm, clip]{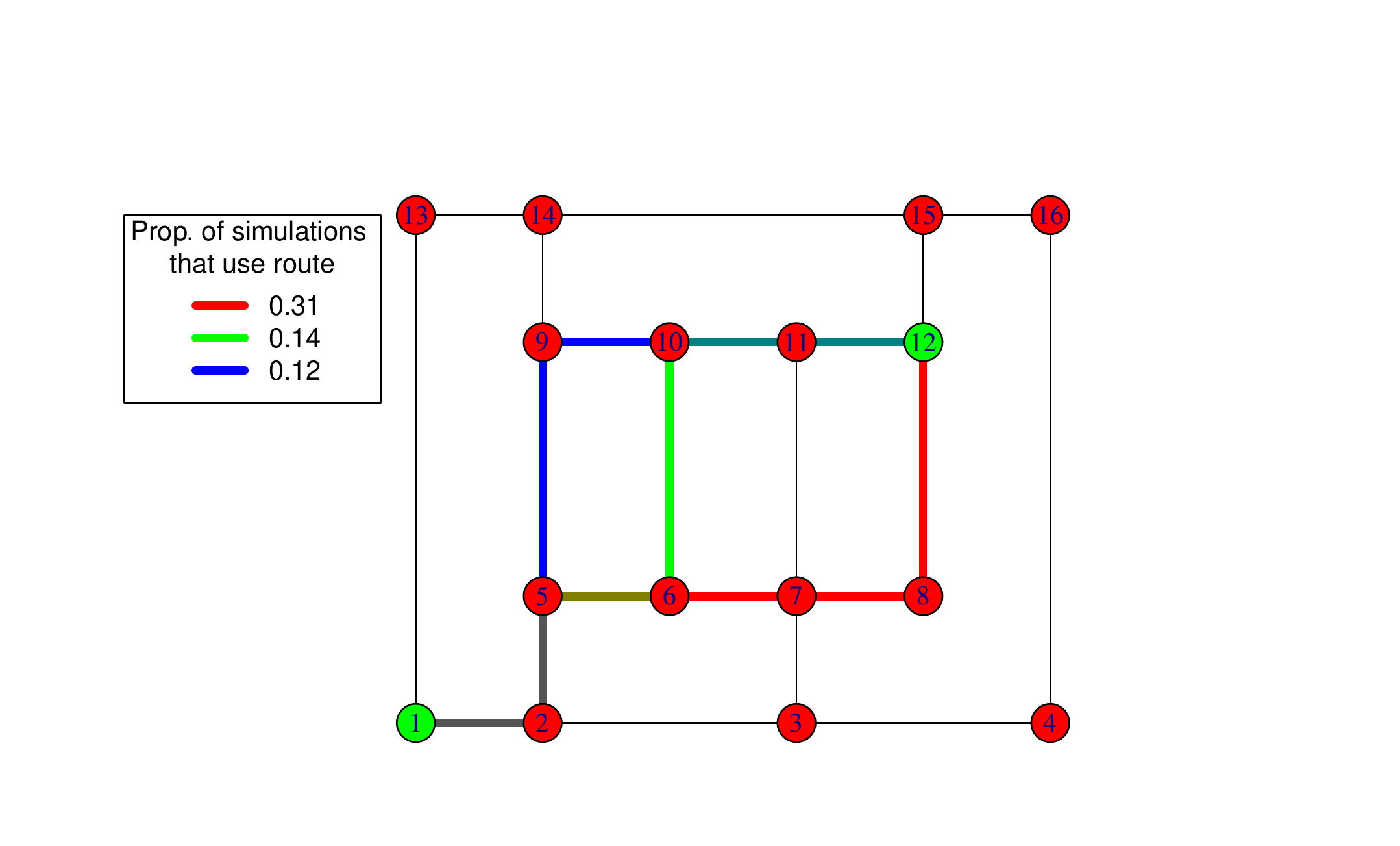}%
		\includegraphics[width=0.5\textwidth, trim=1cm 0cm 4cm 1cm, clip]{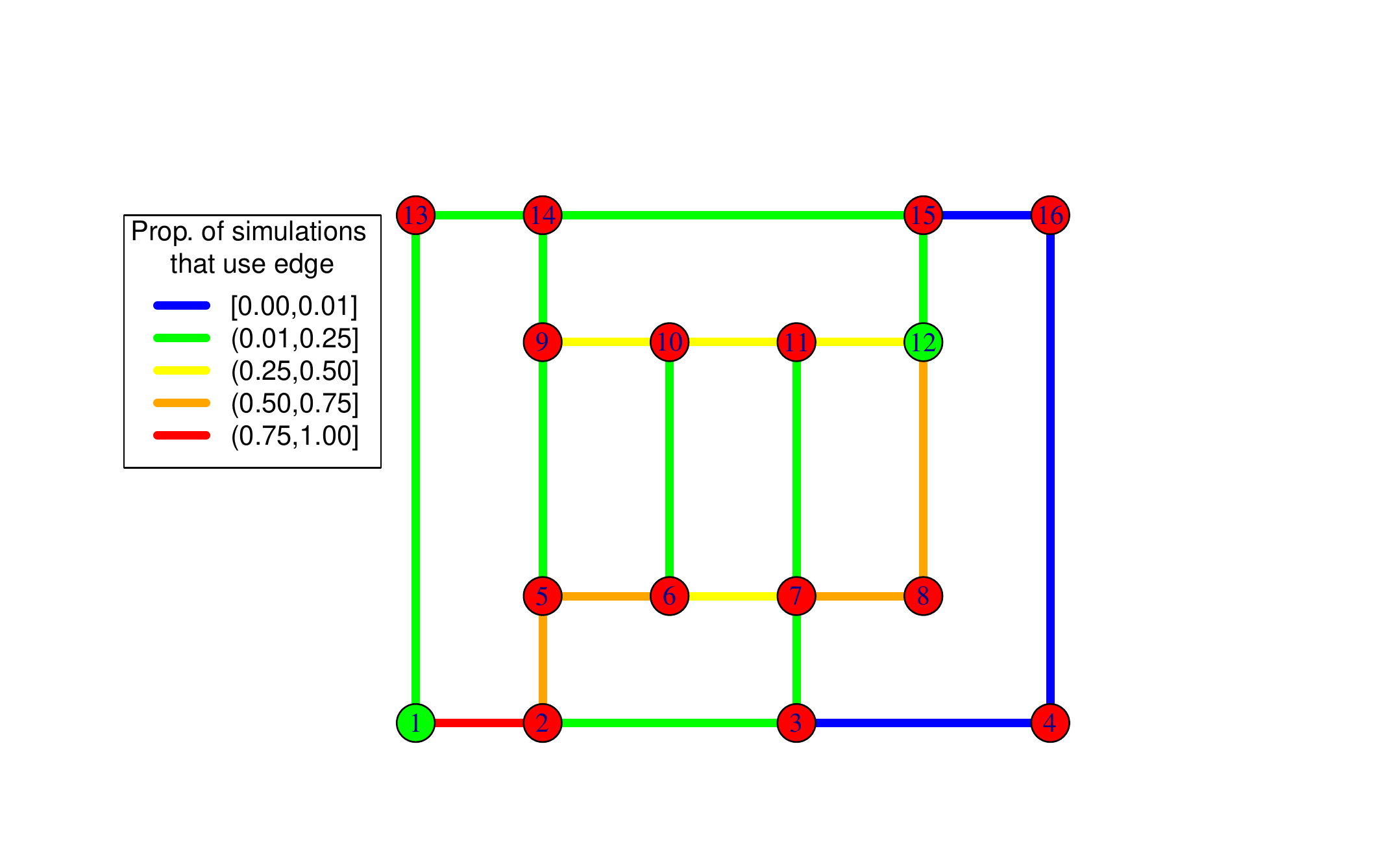}
		\caption{Using the sum of the squared difference of the expected travel time of consecutive edges as objective function. The left picture shows the three routes that minimize the objective function most often, whereas the right picture shows the frequency with which each edge is part of the optimal route.}
		\label{fig:util4}
	\end{center}
\end{figure}

When aiming at minimizing the {\it sum} of the squared differences of the expected travel times at consecutive edges, there is now a downside to taking longer routes (since the sum will include more terms) as well as routes through which it is more difficult to keep a constant velocity.
As a consequence of this the right plot in Figure~\ref{fig:util4} shows more concentration. 
In fact, the left plot shows a clear preference for taking the red route, or otherwise traveling to vertex 5 and then either making use of the path graph $5-6-7-8$ or of the path graph $9-10-11-12$ to reach the destination.

\bigskip

In this section we modeled the preferences of the drivers in terms of various objective functions. The experiments showed that different objective functions may lead to very different optimal routes. Whereas the expected travel time solely relies on the point estimate of this quantity, the quantile of the posterior expectation also takes the estimation uncertainty into account and favors well-explored routes. Besides the expected travel time, drivers may also want to incorporate the travel time variability into their objective function. This demand can be fulfilled by basing the decision on the quantile of the distribution of the estimator of the expected travel time. The precise quantile to be considered reflect the risk-averseness of the driver. One is of course free to come up with alternative objective functions. For example, drivers who strongly dislike velocity fluctuations (for driving comfort, or for reducing fuel consumption) may want to minimize the sum- or mean of the squared difference of the expected velocities of the consecutive edges.

\section{Discussion and concluding remarks}\label{sec:disc}


In this paper we have developed a framework for estimating travel times in a road traffic network. We have followed a Bayesian approach that produces estimates under minimal distributional assumptions. In a series of experiments we have assessed several aspects of the resulting estimation procedure. In addition we have argued how our framework can be used to support route selection. 

In operations research the classical paradigm is to separate the estimation phase from the decision making phase: it is common practice to work with stochastic models assuming that the underlying distributions and parameters are known. In this paper we depart from this approach, in that we advocate taking into account estimation uncertainty when selecting a route. As such, our work can be seen as part of the branch of research in which learning and optimization are integrated; see e.g., \cite{CESA,KESK} for other examples. Our data-driven approach is facilitated by the abundant travel time measurements that are available nowadays due to GPS-based technologies. Another example of a domain in which decisions are made by explicitly taking into account the estimation error, is that of measurement-based admission control \cite{DUFF, GANE,GROS}: based on estimates of the bandwidth consumption of traffic streams that are currently present, it is decided whether newly arriving streams can be accommodated.

We use a Gaussian model to represent the joint distribution of the mean travel time at each of the different edges of a traffic network.
This model provides a good approximation for the distribution of the data under a large variety of sampling regimes, e.g., measuring travel times of randomly chosen particles traversing each edge.
The expectation is endowed with an appropriate prior leading to a (joint) posterior distribution on the expected travel times for the entire network.

The posterior variance-covariance matrix is estimated using the empirical Bayes approach, while the posterior expected travel time for each edge can be estimated from respective (empirical) posterior distribution.
The empirical posterior provides not just estimates for the expected travel times but also quantifies the uncertainty in said estimates.
Furthermore, from this posterior one can also explicitly compute the joint posterior distribution for the expected travel time of a collection of paths on the graph.
Our approach therefore provides estimates and uncertainty quantification for expected travel times for arbitrary paths on the network, as well as other functionals of the model parameters.

The introduction of a higher resolution version of the network allows consistent inference as long as either the resolution of the network (which determines at what spatial resolution data can be collected, e.g., via GPS signal) or the number of observation collected at each edge increases.
Use of higher resolutions makes the approach robust to how traffic flows through the network as particles are not required to keep constant velocity while traversing edges.
This notion of resolution is useful to capture features such as slowdowns at intersections or curves.

We ran multiple simulations to illustrate the performance of our approach.
We also explored the possibility of utilizing the posterior distribution on expected travel times on the network as a statistical tool to support finding optimal paths on the network, with the utility of the path determined by the user.

\medskip 

Several directions for follow-up research can be thought of --- we here provide three possible themes. In the first place, one could focus on operationalizing the approach presented in this paper in a practical context. Secondly, one could aim at developing systematic procedures for route selection based on the estimates produced by our estimation framework (possibly with an interface by which a driver can, explicitly or implicitly, reveal her utility curve). Finally, one could try to explicitly incorporate specific features of the network at hand (such as: including the precise locations of pedestrian crossings, knowledge of the algorithm used by traffic lights at intersections, speed limits that are imposed on specific individual segments, etc.).



\clearpage
\appendix
\section{Auxiliary Results}
\label{app:auxiliary}

In this appendix we present a series results that have appeared in the main text, and provide their respective proofs.


\begin{proposition}
\label{prop:posterior}
Under the modeling assumption in~\eqref{eq:model}, and supposing that $\BFmu_\BFr\mid\big(\lambda, \BFSigma_\BFr^{(\BFn)}\big)$ is endowed with the prior in~\eqref{eq:prior},
the corresponding posterior distribution is
\begin{equation}
\BFmu_\BFr\mid\big(\lambda, \BFSigma_\BFr^{(\BFn)}, \BFX^{(\BFn)}\big)\sim
{\mathscr N}\Big(\hat{\BFmu}_\BFr(\lambda, \BFSigma_\BFr^{(\BFn)}),\; \big(\{\BFSigma_\BFr^{(\BFn)}\}^{-1} + \lambda \bar{\BFL}_\BFr \big)^{-1}\Big), 
\qquad \lambda> 0.
\tag{\ref{eq:posterior}}
\end{equation}
\end{proposition}
\proof{}
Note that $\big(\BFmu_\BFr,\BFX^{(\BFn)}\big)\mid\big(\lambda, \BFSigma_\BFr^{(\BFn)}\big)$ has a joint Gaussian distribution since the prior was chosen independently of $\BFX^{(\BFn)}$. 
Using the law of total expectation and the law of total variance, it follows by standard algebra that
\[
\big(\BFmu_\BFr,\BFX^{(\BFn)}\big)\mid\big(\lambda, \BFSigma_\BFr^{(\BFn)}\big)\sim
{\mathscr N}\left(
\begin{pmatrix}
\BFzero \\
\BFzero
\end{pmatrix}
,\;
\begin{pmatrix}
\frac1\lambda \bar{\BFL}_\BFr^- & \frac1\lambda \bar{\BFL}_\BFr^-\\
\frac1\lambda \bar{\BFL}_\BFr^- & \BFSigma_\BFr^{(\BFn)} + \frac1\lambda \bar{\BFL}_\BFr^-
\end{pmatrix}
\right).
\]
Using the well known expression for the conditional distribution of Gaussian random vectors, we have that
\[
\BFmu_\BFr\mid\big(\lambda, \BFSigma_\BFr^{(\BFn)}, \BFX^{(\BFn)}\big)\sim
{\mathscr N}\Bigg(\frac1\lambda\bar{\BFL}_\BFr^-\Big(\BFSigma_\BFr^{(\BFn)}+\frac1\lambda \bar{\BFL}_\BFr^-\Big)^{-1}\BFX^{(\BFn)},\; \frac1\lambda \bar{\BFL}_\BFr^- - \frac1\lambda \bar{\BFL}_\BFr^-\Big(\BFSigma_\BFr^{(\BFn)}+\frac1\lambda \bar{\BFL}_\BFr^-\Big)^{-1}\frac1\lambda \bar{\BFL}_\BFr^-\Bigg).
\]
It is clear that
\[
\frac1\lambda\bar{\BFL}_\BFr^-\Big(\BFSigma_\BFr^{(\BFn)}+\frac1\lambda \bar{\BFL}_\BFr^-\Big)^{-1}\BFX^{(\BFn)} =
\Big(\BFI_{q_\BFr} + \lambda \BFSigma_\BFr^{(\BFn)}\bar{\BFL}_\BFr \Big)^{-1}\BFX^{(\BFn)} =
\hat{\BFmu}_\BFr(\lambda, \BFSigma_\BFr^{(\BFn)}).
\]
By applying the matrix inversion lemma \cite[Section 0.7.4]{HORN} to the posterior variance, the proof is complete.
\endproof

\begin{proposition}
\label{prop:covariance_estimate}
Consider the variance-covariance matrix given in~\eqref{eq:eigen-decomposition}, and suppose that 
\begin{equation}
\label{eq:def_theta_hat}
\hat{\BFtheta} =
\arg\min_{\BFtheta}\, 
\big\{\BFX^{(\BFn)}\big\}^\top \Big(\BFSigma_\BFr^{(\BFn)}(\BFtheta) + \frac1\lambda \bar{\BFL}_\BFr^-\Big)^{-1}\BFX^{(\BFn)}
 + \ln \Big|\BFSigma_\BFr^{(\BFn)}(\BFtheta) + \frac1\lambda \bar{\BFL}_\BFr^-\Big|.
\tag{\ref{eq:empirical_bayes_estimator}}
\end{equation}
The matrix $\hat{\BFSigma}_\BFr=\BFSigma_\BFr^{(\BFn)}(\hat{\BFtheta})$ satisfies the relation
\begin{equation}
\label{eq:covariance_estimate}
\BFSigma_\BFr^{(\BFn)}(\hat\BFtheta) =
\frac{
\big(\BFX^{(\BFn)} - \BFH(\lambda,\hat\BFtheta)\BFX^{(\BFn)}\big)\big(\BFX^{(\BFn)} - \BFH(\lambda,\hat\BFtheta)\BFX^{(\BFn)}\big)^\top }{\tr\big(\BFI_{q_\BFr} - \BFH(\lambda,\hat\BFtheta)^\top \big)}.
\end{equation}
\end{proposition}
\proof{}
The $i$-th component of a solution $\BFtheta$ in~\eqref{eq:empirical_bayes_estimator} should satisfy the first-order condition
\[
\big\{\BFX^{(\BFn)}\big\}^\top \frac{\partial}{\partial\theta_i}\Big(\BFSigma_\BFr^{(\BFn)}(\BFtheta) + \frac1\lambda \bar{\BFL}_\BFr^-\Big)^{-1}\BFX^{(\BFn)}
 + \frac{\partial}{\partial\theta_i}\ln \Big|\BFSigma_\BFr^{(\BFn)}(\BFtheta) + \frac1\lambda \bar{\BFL}_\BFr^-\Big| = 0.
\]
Standard results from matrix calculus are that
\[
\frac{\partial \BFM(\BFtheta)^{-1}}{\partial\theta_i} =
-\BFM(\BFtheta)^{-1} \frac{\partial \BFM(\BFtheta)}{\partial\theta_i} \BFM(\BFtheta)^{-1},
\qquad\text{} \frac{\partial \ln |\BFM(\BFtheta)|}{\partial\theta_i} =
\tr\left(\BFM(\BFtheta)^{-1}\frac{\partial \BFM(\BFtheta)}{\partial\theta_i}\right),
\]
so that, relying on the eigen-decomposition~\eqref{eq:eigen-decomposition} of $\BFSigma_\BFr^{(\BFn)}(\BFtheta)$, we conclude that for each $i=1, \dots, q_\BFr$,
\begin{align*}
\tr\Big(\big(\BFSigma_\BFr^{(\BFn)}(\BFtheta)& + \frac1\lambda \bar{\BFL}_\BFr^-\big)^{-1} \BFE_i \Big) \\&=
\big\{\BFX^{(\BFn)}\big\}^\top \Big(\BFSigma_\BFr^{(\BFn)}(\BFtheta) + \frac1\lambda \bar{\BFL}_\BFr^-\Big)^{-1} \BFE_i 
\Big(\BFSigma_\BFr^{(\BFn)}(\BFtheta) + \frac1\lambda \bar{\BFL}_\BFr^-\Big)^{-1}\,\BFX^{(\BFn)}.
\end{align*}
Multiplying both sides of the $i$-th equation with $\theta_i$ and add the resulting $q_\BFr$ equations we get that the solution should satisfy the relation
\begin{align*}
\tr\Big(\big(\BFSigma_\BFr^{(\BFn)}(\BFtheta)& + \frac1\lambda \bar{\BFL}_\BFr^-\big)^{-1} \BFSigma_\BFr^{(\BFn)}(\BFtheta)\Big)\\& =
\big\{\BFX^{(\BFn)}\big\}^\top \Big(\BFSigma_\BFr^{(\BFn)}(\BFtheta) + \frac1\lambda \bar{\BFL}_\BFr^-\Big)^{-1} \BFSigma_\BFr^{(\BFn)}(\BFtheta) \Big(\BFSigma_\BFr^{(\BFn)}(\BFtheta) + \frac1\lambda \bar{\BFL}_\BFr^-\Big)^{-1} \BFX^{(\BFn)}.
\end{align*}

It is straightforward to check that
\begin{equation}
\label{eq:smoother_identity}
\Big(\BFSigma_\BFr^{(\BFn)}(\BFtheta) + \frac1\lambda \bar{\BFL}_\BFr^-\Big)^{-1} \BFSigma_\BFr^{(\BFn)}(\BFtheta) =
\BFI_{q_\BFr} - \BFH(\lambda,\BFtheta)^\top ,
\end{equation}
where we parametrize the smoother matrix $\BFH$ in terms of $\BFtheta$ so that
\begin{equation}
\tag{\ref{eq:smoother}}
\BFH(\lambda,\BFtheta) =
\Big(\BFSigma_\BFr^{(\BFn)}(\BFtheta)^{-1} + \lambda \bar{\BFL}_\BFr \Big)^{-1}\BFSigma_\BFr^{(\BFn)}(\BFtheta)^{-1}.
\end{equation}
We conclude that the solution must satisfy
\[
\tr\big(\BFI_{q_\BFr} - \BFH(\lambda,\BFtheta)^\top \big) =
\big(\BFX^{(\BFn)} - \BFH(\lambda,\BFtheta)\BFX^{(\BFn)}\big)^\top \BFSigma_\BFr^{(\BFn)}(\BFtheta)^{-1} \big(\BFX^{(\BFn)} - \BFH(\lambda,\BFtheta)\BFX^{(\BFn)}\big) 
\Big).
\]
which, using the invariance under cyclical permutations of the trace, can also be written as
\begin{equation}
\label{eq:solution_condition}
\tr\big(\BFI_{q_\BFr} - \BFH(\lambda,\BFtheta)^\top \big) =
\tr\big(
\BFSigma_\BFr^{(\BFn)}(\BFtheta)^{-1} \big(\BFX^{(\BFn)} - \BFH(\lambda,\BFtheta)\BFX^{(\BFn)}\big)\big(\BFX^{(\BFn)} - \BFH(\lambda,\BFtheta)\BFX^{(\BFn)}\big)^\top \big).
\end{equation}
We finish the proof by noting that~\eqref{eq:covariance_estimate} solves the above.
\endproof

\begin{corollary}
\label{cor:covariance_estimate}
Denote by $\BFH_i(\lambda,\BFtheta)$ the $(r_i+1)\times(r_i+1)$ sub-matrix of $\BFH(\lambda,\BFtheta)$ corresponding to the sub-edges of edge $i$, and by $\BFX^{(\BFn)}_i$ the averages collected at the sub-edges of edge $i$.
Assume that the variance-covariance matrix of the data satisfies~\eqref{eq:our_covariance}.
Then empirical Bayes estimators of the $\sigma_i^2$, $i=1,\dots,q$, are, for $\lambda>0$, given by
\begin{equation}
\label{eq:variances_estimators}
\hat\sigma_i^2 =
\frac{\big\{\BFX_i^{(\BFn)}\big\}^\top \big(\BFI_{r_i+1}-\BFH_i(\lambda,\BFone)\big)^\top \big(\BFI_{r_i+1}-\BFH_i(\lambda,\BFone)\big) \BFX_i^{(\BFn)}}{\tr\big(\BFI_{r_i+1}-\BFH_i(\lambda,\BFone)\big)/n}.
\end{equation}
\end{corollary}
\proof{}
Note that the variance-covariance matrix in~\eqref{eq:our_covariance}, is of a parametric form $\BFSigma_\BFr^{(\BFn)}(\BFtheta)$ where $\BFtheta=(\sigma_1^2, \dots, \sigma_q^2)^\top $.
As such, following the same argument as in the proof of Proposition~\ref{prop:covariance_estimate}, the estimate of each $\sigma_i^2$, $i=1,\dots,q$, must satisfy
\[
\tr\Big(\big(\BFSigma_\BFr^{(\BFn)}(\BFtheta) + \frac1\lambda \bar{\BFL}_\BFr^-\big)^{-1} \frac{\BFE_i}n \Big) =
\big\{\BFX^{(\BFn)}\big\}^\top \Big(\BFSigma_\BFr^{(\BFn)}(\BFtheta) + \frac1\lambda \bar{\BFL}_\BFr^-\Big)^{-1} \frac{\BFE_i}n 
\Big(\BFSigma_\BFr^{(\BFn)}(\BFtheta) + \frac1\lambda \bar{\BFL}_\BFr^-\Big)^{-1} \BFX^{(\BFn)},
\]
since ${\partial \BFSigma_{\BFr}^{(\BFn)}(\BFtheta)}/{\partial\sigma^2_i} = n^{-1}\BFE_i$, where $\BFE_i:=\diag\{\BFzero,\dots,\BFzero,\BFI_{r_i}, \BFzero,\dots,\BFzero\}$.
Recalling the definition of $\BFH_i(\lambda,\BFtheta)$ and Eqn.~\eqref{eq:smoother_identity}, we solve, for each $i=1,\dots,q$,
\begin{align*}
\frac n{\sigma_i^2} \tr\big(\BFI_{r_i+1}-&\BFH_i(\lambda/\sigma_i^2,\BFone)\big) \\&=
\frac{n^2}{\sigma_i^4} \big\{\BFX_i^{(\BFn)}\big\}^\top \Big(\BFI_{r_i+1}-\BFH_i(\lambda/\sigma_i^2,\BFone)\Big)^\top \Big(\BFI_{r_i+1}-\BFH_i(\lambda/\sigma_i^2,\BFone)\Big) \BFX_i^{(\BFn)},
\end{align*}
which, since the smoothing parameter $\lambda$ is arbitrary, gives estimators of the form \eqref{eq:variances_estimators}.
\endproof

\begin{proposition}
\label{prop:eigendecomposition}
Consider $\bar{\BFL}_\BFr$, the Laplacian of the line graph of $G_\BFr$ of dimension $q_\BFr\times q_\BFr$, for $R=(r_1, \dots, r_q)$.
There exists an orthonormal matrix $\BFOmega$ of the form
\footnote{\tiny Explicit expressions for the entries of $\BFOmega$ can be found in the proof.}, with
$r_{q+1}=2q$,
\begin{equation}\label{eq:omega}
\left(
\frac{\BFomega_{1,1}}{\|\BFomega_{1,1}\|}, \dots, \frac{\BFomega_{1,r_1-1}}{\|\BFomega_{1,r_1-1}\|}, \dots,
\frac{\BFomega_{q,1}}{\|\BFomega_{q,1}\|}, \dots, \frac{\BFomega_{q,r_q-1}}{\|\BFomega_{1,r_q-1}\|} ,
\frac{\BFomega_{q+1,1}}{\|\BFomega_{q+1,1}\|}, \dots, \frac{\BFomega_{q+1,r_{q+1}}}{\|\BFomega_{q+1,r_{q+1}}\|}
\right)^\top,
\end{equation}
such that for all $j=1, \dots, r_i-1$, $i=1,\dots,q+1$ and $j'=1, \dots, r_{i'}-1$, $i'=1,\dots,q+1$, 
\[
\frac{\BFomega_{i,j}^\top}{\|\BFomega_{i,j}\|} \bar{\BFL}_\BFr \frac{\BFomega_{i',j'}}{\|\BFomega_{i',j'}\|} =
\tilde\ell_{i,j} \, \delta_{i,i'} \delta_{j,j'} + 
\frac{\BFomega_{i,j}^\top}{\|\BFomega_{i,j}\|} \BFDelta \frac{\BFomega_{i',j'}}{\|\BFomega_{i',j'}\|},
\]
where for $\delta_{i,j}:=\ind\{i=j\}$ and any $r\leqslant \min_{i=1,\dots,q}r_i$,
\[
\left|\frac{\BFomega_{i,j}^\top}{\|\BFomega_{i,j}\|} \BFDelta \frac{\BFomega_{i',j'}}{\|\BFomega_{i',j'}\|}\right| \le
\frac 4r\delta_{i,i'}\big(1-\delta_{i,q+1}\big) + 
\sqrt{\frac 2r}\big(1-\delta_{i,i'}\big)\big(1-\delta_{i,q+1}\delta_{i',q+1}\big) + 
\delta_{j,j'}\delta_{i,q+1}\delta_{i',q+1},
\]
and where
\[
\tilde\ell_{i,j} :=
4 \sin\left(\frac{\pi (j-1)}{2\, r_i}\right)^2, 
\qquad j = 1, \dots, r_i,\; i=1,\dots,q,
\]
and where $\tilde\ell_{q+1,j}$, $j=1,\dots,2q$, are at most the highest degree of a vertex in $G$.
In addition, the (symmetric) matrix $\BFDelta$ has at least $q_\BFr-6q$ rows (and therefore eigenvalues) equal to $0$, at most $2q$ eigenvalues with norm at most $1$, and at most $4q$ eigenvalues with norm at most~$2$.
\end{proposition}

\proof{}
Label the $q_\BFr$ vertices in the line graph $\bar G_\BFr$ of $G_\BFr$ according to the following ordering of the edges of $G_\BFr$: 
$e_{1,2}$, $\ldots$, $e_{1,r_1}$, $e_{2,2}$, $\cdots$, $e_{2,r_2}$, $\cdots$, $e_{q,2}$, $\ldots$, $e_{q,r_q}$, and label the remaining $2q$ vertices arbitrarily.
The first set of edges are those in $G_\BFr$ that are not incident to vertices of $G$, while the remaining $2q$ edges are those that are.

Consider $\bar{\BFL}_\BFr$, the graph Laplacian of the line graph of $G_\BFr$ and define the $q_\BFr\times q_\BFr$ matrix $\tilde{\BFL}_\BFr$ by the block-diagonal matrix
\begin{equation}
\label{eq:approx_def}
\tilde{\BFL}_\BFr = \diag\big\{\BFL_{r_1-1}, \dots, \BFL_{r_q-1}, \BFL_\star\big\},
\end{equation}
where 
$\BFL_\star$ is a $2q\times 2q$ diagonal matrix containing the degrees of the vertices in $\bar G_\BFr$ corresponding to the last $2q$ edges in $G_\BFr$, and 
$\BFL_r$ is the $r\times r$ Laplacian matrix of a path graph with $\BFr$ vertices,
\[
\BFL_r =
\begin{pmatrix}
1 & 0 & 0 & \cdots & 0\\
0 & 2 & 0 & \ddots & \vdots\\
0 & \ddots & \ddots & \ddots & 0\\
\vdots & \ddots & 0 & 2 & 0\\
0 &\cdots & 0 & 0 & 1
\end{pmatrix} - \begin{pmatrix}
0 & -1 & 0 & \cdots & 0\\
-1 & 0 & -1 & \ddots & \vdots\\
0 & \ddots & \ddots & \ddots & 0\\
\vdots & \ddots & -1 & 0 & -1\\
0 &\cdots & 0 & -1 & 0
\end{pmatrix} =
\begin{pmatrix}
1 & -1 & 0 & \cdots & 0\\
-1 & 2 & -1 & \ddots & \vdots\\
0 & \ddots & \ddots & \ddots & 0\\
\vdots & \ddots & -1 & 2 & -1\\
0 &\cdots & 0 & -1 & 1
\end{pmatrix}.
\]
Following~\cite[Section 1.4.4]{BROU}, the eigenvalues of $\BFL_r$ are
\[
2 - 2 \cos\left(\frac{\pi k}{r}\right) =
4 \sin\left(\frac{\pi k}{2 r}\right)^2, 
\qquad k = 0, \dots, r-1,
\]
so that in particular the eigenvalues belong to the interval $[0,4)$;
the eigenvector corresponding to the eigenvalue $2-\zeta-\zeta^{-1}$ is the $\BFr$-dimensional vector
\[
\BFomega_\zeta =
\big(1+\zeta^{2r-1},\, \dots,\, \zeta^{i-1} + \zeta^{2r-i},\, \dots,\, \zeta^{r-1} + \zeta^r\big)^\top,
\]
where $\zeta^{2r}=1$ so that $\zeta\neq0$.
Note that the norm of this eigenvector is
\[
\BFomega_\zeta^\top\BFomega_\zeta = 
\sum_{i=0}^{r-1}\zeta^{2j} + 2 r \zeta^{2r-1} + \sum_{j=0}^{r-1}\zeta^{-2(j+1)} =
\frac{1-\zeta^{2r}}{1-\zeta^2} + 2r\zeta^{-1} + \zeta^{-2}\frac{1-\zeta^{-2r}}{1-\zeta^{-2}} = 
2r\zeta^{-1},
\]
where we use the fact that $\zeta^{2r}=1$.
It is also straightforward to see that since the $\zeta\in\mathbb{C}$ corresponding to the eigenvalue $2-\zeta-\zeta^{-1}$ must belong to $[0,4)$, then the modulus of $\zeta$ is at most $1$.

Let $\tilde\ell_{i,j}$, $j=1,\dots, r_i-1$, $i = 1, \dots, q$, represent the $j$-th eigenvalue of $\BFL_{r_i-1}$, and define, for $j=1,\dots, r_i-1$ and $i = 1, \dots, q$,
\[
\BFomega_{i,j} = 
\big(\BFzero_{r_1-1}^\top,\, \dots,\, \BFzero_{r_{i-1}-1}^\top,\, 
\BFomega_{\zeta_{i,j}}^\top,\, 
\BFzero_{r_{i+1}-1}^\top,\, \dots,\, \BFzero_{r_q-1}^\top,\BFzero_{2q}^\top\big)^\top 
\in \mathbb{R}^{q_\BFr},
\]
where $\BFomega_{\zeta_{i,j}}$ is the eigenvector corresponding to the eigenvalue $\tilde\ell_{i,j}$ of $\BFL_{r_i-1}$, and $\BFzero_i\in\mathbb{R}^i$ is a the zero vector.
Define also $\tilde\ell_{q+1,1}, \dots, \tilde\ell_{q+1,2q}$ to be the diagonal elements of $\BFL_\star$, and
\[
\BFomega_{q+1,j} = 
\big(\BFzero_{r_1-1}^\top,\, \dots,\, \BFzero_{r_q-1}^\top,
0,\dots,0,1,0,\dots,0\big)^\top, \qquad
j=1,\dots, 2q,
\]
where the $1$ is on position $q_\BFr - 2q + j$.
We conclude that, by construction, for $j=1,\dots,r_i$, $i = 1, \dots, q+1$, and $j'=1,\dots,r_{i'}$, $i'=1,\dots,q+1$ where we set $r_{q+1}:=2q$,
\[
\BFomega_{i',j'}^\top\tilde{\BFL}_R\BFomega_{i,j} =
\BFomega_{\zeta_{i,j}}^\top\BFL_{r_i-1}\BFomega_{\zeta_{i',j'}} = 
\tilde\ell_{i,j}\BFomega_{\zeta_{i,j}}^\top\BFomega_{\zeta_{i',j'}} = 
\tilde\ell_{i,j}\BFomega_{i',j'}^\top\BFomega_{i,j}.
\]
Note further that
\begin{equation}
\label{eq:norm}
\BFomega_{i',j'}^\top\BFomega_{i,j} = 
\BFomega_{\zeta_{i',j'}}^\top\BFomega_{\zeta_{i,j}} =
2 r_i\, \zeta_{i,j}^{-1}\, \delta_{i,i'} \delta_{j,j'},
\end{equation}
so that the $q_\BFr\times q_\BFr$ matrix $\Omega$, as given through \eqref{eq:omega}, is an orthonormal matrix such that $\BFOmega^\top\BFOmega=\BFI_{q_\BFr} = \BFOmega\,\BFOmega^\top$. In addition, $\Omega$ diagonalizes $\tilde{\BFL}_\BFr$, in that
\[
\BFOmega^\top \tilde{\BFL}_\BFr \BFOmega = 
\diag\big\{
\tilde\ell_{1,1}, \dots, \tilde\ell_{1,r_1-1},
\, \dots,\, 
\tilde\ell_{q+1,1}, \dots, \tilde\ell_{q+1,2q}
\big\}.
\]
Note also that the matrix $\BFOmega$ is block diagonal.

Denoting $\BFDelta := \bar{\BFL}_\BFr - \tilde{\,\BFL}_\BFr$, this matrix has the following structure:\[
\BFDelta =
\left(\begin{array}{cccc|c}
\BFD_{1,1}  & \BFzero  & \cdots & \BFzero  & \BFD_{1,q+1}  \\
\BFzero  & \BFD_{2,2}  & \ddots & \vdots & \vdots  \\
\vdots & \ddots& \ddots & \BFzero  & \BFD_{q-1,q+1} \\
\BFzero  & \cdots & \BFzero  & \BFD_{q,q}  & \BFD_{q,q+1} \\
&&&&\vspace{-3.2mm}
 \\ \hline
&&&&\vspace{-3.2mm}\\
\BFD_{1,q+1}^\top  & \cdots  & \BFD_{q-1,q+1}^\top & \BFD_{q,q+1}^\top  & \BFD_{\star} 
\end{array}\right)
\]
where the block structure is the same as in $\tilde{\,\BFL}_\BFr$.
As such, each of the (symmetric) matrices $\BFD_{i,i}$, $i = 1, \dots, q$ is of dimension $(r_i-1)\times(r_i-1)$, $\BFD_\star$ is of dimension $2q\times 2q$ and symmetric, and the matrices $\BFD_{i,q+1}$, $i=1,\dots, q$ are of dimension $(r_i-1)\times 2q$.
It is clear that $\BFD_{i,i}=\diag\{1,0,\dots,0,1\}$, $i = 1, \dots, q$, because of the ordering that we picked for the edges in $G_\BFr$ and since all of the vertices on which those edges are incident have degree 2 in $\bar G_\BFr$.
The matrix $\BFD_\star$ has zeroes on its diagonal and outside the diagonals it has $-1$'s; in each row, the number of 
$-1$'s is at most the maximal degree on $G$.
Finally, each matrix $\BFD_{i,q+1}$, $i=1,\dots, q$, has exactly two entries equal to $-1$ (on different rows and columns), and all other entries are equal to $0$; the two $-1$ entries correspond to how each path subgraph connects to some edge in the original graph $G$, with the exact locations of the $-1$'s depending on the ordering of the path subgraphs.

Based on the description above it follows immediately from Gershgorin's circle theorem \cite{GER, VAR}, that $\BFDelta$ has at least $q_\BFr-6q$ eigenvalues equal to zero (due to rows of zeros), at most $2q$ eigenvalues with norm at most 1 (due to the last $2q$ rows), and at most $4q$ eigenvalues with norm at most 2 (due to the remaining rows.)

We now bound the absolute value of the $(j,j')$ entry of the $(i,i')$ block of $\BFOmega^\top \,\BFDelta \,\BFOmega$, which is
\[
\frac{\BFomega_{i,j}^\top}{\|\BFomega_{i,j}\|} \BFDelta \frac{\BFomega_{i',j'}}{\|\BFomega_{i',j'}\|}.
\]
Using the symmetry of $\BFDelta$, there are three relevant cases to consider corresponding to:
a) $i,i'\in\{1,\dots,q\}$;
b) $i\in\{1,\dots,q\}$, $i'=q+1$; or
c) $i=i'=q+1$.

\paragraph{Case a):} If $i\neq i'$, then we immediately conclude that $\BFomega_{i,j}^\top\BFDelta\BFomega_{i',j'} = 0$; assume then that $i=i'$.
Using the fact that the matrices $\BFD_{i,i}$ are symmetric and idempotent in combination with the Cauchy-Schwarz inequality,
\[
\BFomega_{i,j}^\top\BFDelta\BFomega_{i,j'} =
\BFomega_{\zeta_{i,j}}^\top\BFD_{i,i}\BFomega_{\zeta_{i,j'}} =
\BFomega_{\zeta_{i,j}}^\top\BFD_{i,i}^T\BFD_{i,i}\BFomega_{\zeta_{i,j'}} \leqslant
\left(\BFomega_{\zeta_{i,j}}^\top\BFD_{i,i}\BFomega_{\zeta_{i,j}} \; \BFomega_{\zeta_{i,j'}}^\top\BFD_{i,i}\BFomega_{\zeta_{i,j'}} \right)^{1/2}.
\]
Furthermore, using the fact that $\zeta_{i,j}^{2r_i}=1$, and the definition of $\BFD_{i,i}$,
\[
\BFomega_{\zeta_{i,j}}^\top\BFD_{i,i}\BFomega_{\zeta_{i,j}} =
\big(1+\zeta_{i,j}^{2r_i-1}\big)^2 + \big(\zeta_{i,j}^{r_i-1} + \zeta_{i,j}^{r_i}\big)^2 =
\big(1+\zeta_{i,j}^{-1}\big)^2 + \big(\zeta_{i,j}^{-1} + 1\big)^2 =
2\big(1+\zeta_{i,j}^{-1}\big)^2.
\]
Hence, using the above and~\eqref{eq:norm}, we conclude that
\[
\frac{\BFomega_{\zeta_{i,j}}^\top\BFD_{i,i} \BFomega_{\zeta_{i,j}}}{\BFomega_{\zeta_{i,j}}^\top\BFomega_{\zeta_{i,j}}} =
\frac{2\big(1+\zeta_{i,j}^{-1}\big)^2}{2 r_i\, \zeta_{i,j}^{-1}} =
\frac{1+2\zeta_{i,j}^{-1}+\zeta_{i,j}^{-2}}{r_i\, \zeta_{i,j}^{-1}} =
\frac{2 + \zeta_{i,j}+\zeta_{i,j}^{-1}}{r_i} = 
\frac4{r_i},
\]
Putting everything together we find that, irrespectively of $j,j'$, if $i,i'\in\{1,\dots,q\}$, then
\[
\left|\frac{\BFomega_{i,j}^\top}{\|\BFomega_{i,j}\|} \BFDelta \frac{\BFomega_{i',j'}}{\|\BFomega_{i',j'}\|}\right| \le
\frac4{r_i}\, \delta_{i,i'}.
\]

\paragraph{Case b):} In this case, because of the block structure of $\BFDelta$ and the two vectors,
\begin{align*}
\BFomega_{i,j}^\top\BFDelta\BFomega_{q+1,j'} &=
\big(\BFzero, \cdots, \BFzero, \BFomega_{\zeta_{i,j}}^\top\BFD_{i,i}, \BFzero, \cdots, \BFzero, \BFomega_{\zeta_{i,j}}^\top \BFD_{i,q+1} \big)\BFomega_{q+1,j'} \\ &=
\BFomega_{\zeta_{i,j}}^\top \BFD_{i,q+1} (\dots, 0, 1, 0, \dots)^\top =
\BFomega_{\zeta_{i,j}}^\top \BFe_{i,j'},
\end{align*}
where 
the $1$ is on the $j'$-th position in a vector from $\{0,1\}^{2q}$, and
$\BFe_{i,j'} \in \{-1,0\}^{r_i-1}$ represents the $j'$ column of $\BFD_{i,q+1}$ which is a vector that has at most one $-1$;
if $\BFe_{i,j'}$ has a $-1$ entry, then the location of this entry depends on $i$ and on the ordering of the last $2q$ vertices in $\bar G_\BFr$.
From this we see that $\BFomega_{i,j}^\top\BFDelta\BFomega_{q+1,j'}$ is either $0$, or it is an entry of the vector $\BFomega_{\zeta_{i,j}}$.
Combining the above with~\eqref{eq:norm}, we conclude that irrespectively of $j,j'$, if $i\in\{1,\dots,q\}$ and $i'=q+1$, then
\[
\left|\frac{\BFomega_{i,j}^\top}{\|\BFomega_{i,j}\|} \BFDelta \frac{\BFomega_{i',j'}}{\|\BFomega_{i',j'}\|}\right| \leqslant
\max_{i=1,\dots,q}\max_{j=1,\dots,r_i-1}\max_{k=1,\dots,2q}
\frac{\big|\zeta_{i,j}^{k-1/2} + \zeta_{i,j}^{-(k-1/2)}\big|}{\sqrt{2 r_i}} \leqslant \max_{i=1,\dots,q}\sqrt{\frac2{r_i}},
\]
where we use the fact that the modulus of $\zeta$ is at most $1$.

\paragraph{Case c):} By construction and by the definition of $\BFD_\star$,
\[
\BFomega_{q+1,j}^\top\BFDelta\BFomega_{q+1,j'} =
(\BFD_{\star})_{j,j'} \in\{-1,0\}.
\]
Using the fact that $\BFomega_{q+1,j}^\top\BFomega_{q+1,j}=1$ we conclude that if $i=i'=q+1$, then
\[
\left|\frac{\BFomega_{i,j}^\top}{\|\BFomega_{i,j}\|} \BFDelta \frac{\BFomega_{i',j'}}{\|\BFomega_{i',j'}\|}\right| \leqslant
\delta_{j,j'}, 
\]
since the diagonal elements of $\BFD_{\star}$ are all $0$.
This concludes the proof.
\endproof

\begin{lemma}
\label{lem:series_sine}
For any $\lambda>0$, define
\[
h_j := 
h_j(\lambda) = 
\frac1{1+\lambda \sin\Big(\frac\pi2\frac{j-1}r\Big)^2}, 
\quad j=1,\dots,r-1.
\]
Assume that $\lambda=o(r^2)$ as $r\to\infty$. Then, as $\lambda\to\infty$, for any $s\in\mathbb{N}$, $t\in\mathbb{N}_0$, 
we have
\[
\sum_{j=1}^{r-1} h_j^s(1-h_j)^t = r \lambda^{-1/2} \kappa_{s,t}\{1+ o(1)\}, 
\qquad\text{where}\quad 
\kappa_{s,t} := \frac{\Gamma(s-\frac{1}{2})\Gamma(t+\frac{1}{2})}{\pi \, \Gamma(t+s)}.
\]
\end{lemma}
\proof{}
For $1\leqslant u \leqslant r+1$, denote $g(u):=h_u^s(1-h_u)^t$. For each $s$ and $t$ there exists $u_{\max}$ such that $g(u)$ is increasing on $u\in[1,u_{\max}]$ and decreasing on $u\in[u_{\max},r+1]$. Observe that
\[
\int_1^{r}g(u)\,{\rm d}u-g(u_{\max}) \leqslant \sum_{j=1}^{r-1} h_j^s(1-h_j)^t \leqslant \int_1^{r}g(u)\,{\rm d}u+g(u_{\max}),
\]
where clearly $\vert g(u_{\max}) \vert\leqslant 1$. Since $\int_r^{r+1}g(u)du\leqslant 1$, we conclude
\[
\sum_{j=1}^{r-1} h_j^s(1-h_j)^t=\int_1^{r+1}g(u)\,{\rm d}u+O(1),\quad \text{as } \lambda\to\infty.
\]
Moreover, we have
\[
\int_1^{r+1}g(u)\,{\rm d}u=\frac{r\lambda^t\Gamma(t+\frac{1}{2})}{\sqrt{\pi}\Gamma(t+1)}{}_2F_1\left(t+\frac{1}{2},s+t;t+1;-\lambda\right),
\]
where ${}_2F_1$ denotes the hypergeometric function. Using that \cite[Eqn.\ 15.3.7]{AS}
\[
{}_2F_1\left(t+\frac{1}{2},s+t;t+1;-\lambda\right)=\lambda^{-t-\frac{1}{2}}\frac{\Gamma(s-\frac{1}{2})\Gamma(t+1)}{\sqrt{\pi}\Gamma(t+s)}(1+o(1)),\quad \text{as } \lambda\to\infty,
\]
we obtain
\[
\int_1^{r+1}g(u)\,{\rm d}u= r \lambda^{-1/2} \kappa_{s,t}\{1+ o(1)\}.
\]
Since $\lambda=o(r^2)$, the conclusion follows.
\endproof


\section{Proof of main result}
\label{app:proofs} 
 
 In the proof of Theorem 
\ref{theo:moments}, the following lemma is used.

\begin{lemma}
\label{lem:traces}
Suppose that~\eqref{eq:our_covariance} holds.
Then, for any $s,t\in\mathbb{N}_0$, and $u,v\in\{0,1\}$ such that $u+2s\geqslant1$,
\begin{align*}
\Big| \tr\Big\{ 
\BFH^u
\big(\BFH^\top\BFH\big)^{s}&
\big(\BFI_{q_\BFr}-\BFH\big)^v
\big\{\big(\BFI_{q_\BFr}-\BFH\big)^\top\big(\BFI_{q_\BFr}-\BFH\big)\big\}^{t}
\Big\} \:-\\ 
&\frac{1}{2}\left(\frac\lambda n\right)^{-1/2} {\kappa_{u+2s,v+2t}} \sum_{i=1}^q \frac{r_i}{\sigma_i}
\{1+o(1)\}
\Big| \leqslant
O(q),
\end{align*}
where $\BFH$ abbreviates $\BFH(\lambda,\BFtheta)$, 
and $\kappa_{u+2s,v+2t}$ is defined in Lemma~$\ref{lem:series_sine}$.
\end{lemma}

\proof{}
Define $\BFh := (h_{1,1}, \dots, h_{1,r_1+1}, h_{2,1}, \dots ,h_{2,r_2+1},\dots,)^\top\in\mathbb{R}^{q_\BFr}$, where
\[
h_{i,j} := 
h_{i,j}(\lambda,n) = 
\frac1{\displaystyle 1+4\,\frac\lambda n \sigma_i^2 \sin\Big(\frac{\pi(j-1)}{2\,r_i}\Big)^2}, 
\qquad j=1,\dots,r_i-1, \quad i=1,\dots, q+1.
\]
We then have that for the matrices $\BFDelta$ and $\BFOmega$ from Proposition~\ref{prop:eigendecomposition}, and $\BFSigma_\BFr^{(\BFn)}(\BFtheta)$ as defined in~\eqref{eq:our_covariance},
\[
\BFH(\lambda,\BFtheta) = 
\Big(\BFI_{q_\BFr} + \lambda \BFSigma_\BFr^{(\BFn)}(\BFtheta)\bar{\BFL}_\BFr \Big)^{-1} =
\Big(\BFOmega\,\diag\{\BFh\}^{-1}\BFOmega^\top + \lambda \BFSigma_\BFr^{(\BFn)}(\BFtheta)\BFDelta \Big)^{-1}.
\]
Using the fact that, with $\BFA$, $\BFB$ denoting two square matrices so that $\BFA+\BFB$ and $\BFA$ are invertible, $(\BFA+\BFB)^{-1} = \BFA^{-1} - (\BFA+\BFB)^{-1}\BFB\BFA^{-1}$, we conclude that 
\[
\BFH(\lambda,\BFtheta) = 
\BFOmega\,\diag\{\BFh\}\BFOmega^\top - \BFdelta,\]
where
\[
\BFdelta := \lambda
\BFH(\lambda,\BFtheta)\,
\BFSigma_\BFr^{(\BFn)}(\BFtheta)\,
\BFDelta\,
\BFOmega\,\diag\{\BFh\}\BFOmega^\top.
\]
Note that $0\leqslant h_{i,j}\leqslant1$, $j=1,\dots,r_i-1$, $i=1,\dots, q+1$.
Also, by definition, the singular values of $\BFH$ are between $0$ and $1$.
By Weyl's inequalities \cite{WEYL} we then conclude that for any $\lambda>0$, the singular values of $\BFdelta$ must be between $-1$ and $1$.
Furthermore, since by Proposition~\ref{prop:eigendecomposition} we know that $\BFDelta$ (and consequently $\BFdelta$) has at least $q_\BFr-6q$ rows and columns of zeroes, and therefore at least that many singular values equal to zero.
This means that
\[
\Big|\tr\big\{\BFdelta(\BFdelta^T\BFdelta)^p\big\}\Big| =
\Big|\tr\big\{\BFdelta^T(\BFdelta^T\BFdelta)^p\big\}\Big| \le
\Big|\tr\big\{(\BFdelta^T\BFdelta)^p\big\}\Big| \le
\Big|\tr\big\{\BFdelta\big\}\Big| \leqslant 6q,
\qquad p\in\mathbb{N},
\]
by Von Neumann's trace inequality \cite{MIR}.
From this we conclude that for any $s,t\in\mathbb{N}_0$, $u,v\in\{0,1\}$, such that $s+t+u+v>0$,
\begin{align*}
\Big| \tr\Big\{ 
\BFH^u
\big(\BFH^T\BFH\big)^{s}
\big(\BFI_{q_\BFr}-\BFH\big)^v&
\big\{\big(\BFI_{q_\BFr}-\BFH\big)^T\big(\BFI_{q_\BFr}-\BFH\big)\big\}^{t}
\Big\} \:-\\ 
&\sum_{i=1}^{q+1}\sum_{j=1}^{r_i-1}h_{i,j}^{u+2s}(1-h_{i,j})^{v+2t} \Big| \leqslant 
O(q),
\end{align*}
where $\BFH$ abbreviates $\BFH(\lambda,\BFtheta)$.
The statement now follows from the result above together with Lemma~\ref{lem:series_sine}, in combination with the triangle inequality.
\endproof

\setcounter{theorem}{0}
\begin{theorem}
Suppose that 
\begin{equation}
\BFX^{(\BFn)} \sim 
{\mathscr N}\Big(\BFmu_{0,\BFr},\; n^{-1}\diag\big(\sigma_{0,1}^2\BFone_{r_1}, \dots, \sigma_{0,q}^2\BFone_{r_q}\big)\Big) \tag{\ref{eq:Xn}}
\end{equation}
for some $\BFmu_{0,\BFr}\in\mathscr{M}_{\BFr}(C)$.
Consider then $n\in\mathbb{N}$, and $\BFr\in\mathbb{N}_0^q$ such that $n = o\big(\min_{i=1,\dots,q}r_i\big)^2$ and define the collection
\[
\Lambda_{n,\BFr}: = 
\Big\{\lambda>0: n = o(\lambda), \lambda = o\big(\min_{i=1,\dots,q}r_i\big)^2\Big\}.
\]
Consider also
\[
\hat\lambda := 
\arg\min_{\lambda\in\Lambda_{n,\BFr}} {\rm GCV}(\lambda),
\]
for ${\rm GCV}(\lambda)$ as defined in~\eqref{eq:GCV}, as well as $\hat\BFtheta = (\hat\sigma_1^2, \dots, \hat\sigma_q^2)$ with each $\hat\sigma_i^2 = \hat\sigma_i^2(\hat\lambda)$ defined as in~\eqref{eq:variances_estimators}, and finally, define \[\hat\BFmu_\BFr := \hat\BFmu_\BFr\big(\hat\lambda, \BFSigma_\BFr^{(\BFn)}(\hat\BFtheta)\big).\]
Then, as long as either $n\to\infty$ or $\min_{i=1,\dots,q}r_i\to\infty$,
$\hat\BFmu_\BFr$ is consistent in probability for $\BFmu_{0,\BFr}$, and each
$\hat\sigma_i^2$ is consistent in probability for $\sigma_{0,i}^2$.
\end{theorem}
\proof{}
In what follows, we make use of the following well known result:
If $\BFY\sim {\mathscr N}\big(\BFy,\, \BFS\big)$, and $\BFM$ is square matrix of appropriate dimension, then
\begin{align*}
\mathbb{E}\BFY^\top\BFM\BFY& =
\BFy^\top\BFM\BFy + \tr\big(\BFM\BFS\big),\\\mathbb{V}\BFY^T\BFM\BFY &=
4\BFy^\top\BFM\BFS\BFM\BFy + 2\tr\big(\BFM\BFS\BFM\BFS\big).
\end{align*}
Suppose that $\BFX^{(\BFn)}$ is distributed according to \eqref{eq:Xn}, where we use the subscript `$0$' to distinguish between the true underlying parameters of the distribution, and arbitrary elements of the underlying parameter sets.

We first derive a risk bound for the estimator of $\BFmu_{0,\BFr}$.
Using $\BFH$ to abbreviate $\BFH(\lambda,\BFtheta)$,
\begin{align*}
\mathbb{E}\big\|\hat{\BFmu}_{\BFr} - \BFmu_{0,\BFr}\big\|^2 &=
\mathbb{E}\big\|\BFH(\BFX^{(\BFn)}-\BFmu_{0,\BFr}) - (\BFI_{q_\BFr}-\BFH)\BFmu_{0,\BFr}\big\|^2\\ & =
\BFmu_{0,\BFr}^\top(\BFI_{q_\BFr}-\BFH)^\top(\BFI_{q_\BFr}-\BFH)\BFmu_{0,\BFr} +
\frac1n\, \tr\big(\BFH\diag\big(\sigma_{0,1}^2\BFone_{r_1}, \dots, \sigma_{0,q}^2\BFone_{r_q}\big)\BFH^\top\big)\\ & \leqslant 
\sigma_{\max}^2\, \frac\lambda n\, P(\BFmu_{0,\BFr}) +
\frac{\sigma_{0,\max}^2}n\tr\big(\BFH^\top\BFH\big)
\end{align*}
since $\BFzero \leqslant \big(\BFI_{q_\BFr}-\BFH\big)^\top\big(\BFI_{q_\BFr}-\BFH\big) \leqslant \lambda n^{-1} \diag\big(\sigma_{0,1}^2\BFone_{r_1}, \dots, \sigma_{0,q}^2\BFone_{r_q}\big) \bar{\BFL}_{\BFr} \leqslant \sigma_{0,\max}^2 \lambda n^{-1} \bar{\BFL}_{\BFr}$, where $\BFA\leqslant\BFB$ means that $\BFB-\BFA$ is positive semi-definite, and where $P(\BFmu_{0,\BFr}) = \BFmu_{0,\BFr}^\top\bar{\BFL}_{\BFr}\BFmu_{0,\BFr}$.

Using Lemma~\ref{lem:traces} and the assumption that $\BFmu_{0,\BFr}\in\mathscr{M}_{\BFr}(C)$ (implying that that~\eqref{eq:bound_P} holds), we conclude that the previous upper bound is majorized by
\[
2\sigma_{\max}^2 \frac{q\,\lambda}{n\min_{i=1,\dots,q} r_i}C^2 + \frac{\sigma_{0,\max}^2}{\sigma_{\min}}\frac1n \left(\frac\lambda n\right)^{-1/2}\kappa_{2,0} \sum_{i=1}^q r_i,
\]
for all appropriately large $\lambda$.
Equating the derivative to $0$, it directly follows that this upper bound is minimized for
\begin{equation}
\label{eq:oracle_order}
\lambda = 
O\left\{\left(\frac{\sigma_{0,\max}^2}{\sigma_{\min}\sigma_{\max}^2\, q\, C^2}\right)^2\, n \Big(\min_{i=1,\dots,q}r_i\, \sum_{i=1}^qr_i\Big)^2\right\}^{1/3},
\end{equation}
which, since $\sum_{i=1}^qr_i\geqslant q \min_{i=1,\dots,q}r_i$, leads to the upper bound
\begin{equation}
\label{eq:risk_order}
\mathbb{E}\big\|\hat{\BFmu}_{\BFr} - \BFmu_{0,\BFr}\big\|^2 \leqslant 
O\left\{ \left(\frac{\sigma_{\max}^2\, \sigma_{0,\max}^4\, C^2}{\sigma_{\min}^2}\right)^{1/3}\sum_{i=1}^qr_i \Big(n \min_{i=1,\dots,q}r_i\Big)^{-2/3}
\right\}.
\end{equation}
We conclude that if the working variances $\sigma_i^2$ are bounded and $\lambda$ is picked as in~\eqref{eq:oracle_order}, then
$(\sum_{i=1}^qr_i)^{-1}\mathbb{E}\big\|\hat{\BFmu}_{\BFr} - \BFmu_{0,\BFr}\big\|^2$
converges to zero as long as either $n$ or $\min_{i=1,\dots,q}r_i$ converges to infinity.

\medskip
The next step is to establish a risk bound for the estimators of $\sigma_{0,i}^2$.
Using the notation that $\BFH_i$ is the principal sub-matrix of $\BFH(\lambda,\BFone)$ corresponding to the sub-edges of the $i$-th edge in $E$, and that, likewise, $\BFX_i^{(\BFn)} \sim {\mathscr N}(\BFmu_{0,i},\; n^{-1}\sigma_{0,i}^2\BFI_{r_i+1})$ is the data collected at sub-edges of the $i$-th edge in $E$, we have
\[
\hat\sigma_i^2 =
\frac{\big\{\BFX_i^{(\BFn)}\big\}^\top \big(\BFI_{r_i+1}-\BFH_i\big)^\top \big(\BFI_{r_i+1}-\BFH_i\big) \BFX_i^{(\BFn)}}{\tr\big(\BFI_{r_i+1}-\BFH_i\big)/n}, 
\qquad i=1,\dots,q.
\]
We directly have that for each $i=1,\dots,q$,
\[
\mathbb{E}\hat\sigma_i^2 = 
\frac{\BFmu_{0,i}^\top \big(\BFI_{r_i+1}-\BFH_i\big)^\top \big(\BFI_{r_i+1}-\BFH_i\big) \BFmu_{0,i}}{\tr\big(\BFI_{r_i+1}-\BFH_i\big)/n} + \sigma_{0,i}^2\,
\frac{\tr\big(\big(\BFI_{r_i+1}-\BFH_i\big)^\top \big(\BFI_{r_i+1}-\BFH_i\big)\big)}{\tr\big(\BFI_{r_i+1}-\BFH_i\big)},
\]
so that the absolute value of the associated bias is at most
\[
\big|{\rm Bias}(\hat\sigma_i^2)\big| =
\big|\mathbb{E}\hat\sigma_i^2 - \sigma_{0,i}^2\big| \leqslant
\frac{\lambda\, P(\BFmu_{0,i})}{\tr\big(\BFI_{r_i+1}-\BFH_i\big)} +
\sigma_{0,i}^2\,
\frac{\tr\big(\BFH_i^\top \big(\BFI_{r_i+1}-\BFH_i\big)\big)}{\tr\big(\BFI_{r_i+1}-\BFH_i\big)},
\qquad i=1,\dots,q.
\]
As for the variance of the estimator, after straightforward simplifications it equals, for any $i=1,\dots,q$,
\[
\mathbb{V}\hat\sigma_i^2 = 
4 \sigma_{0,i}^2 \frac{\lambda\, P(\BFmu_{0,i})}{\tr\big(\BFI_{r_i+1}-\BFH_i\big)^2} + 2 \sigma_{0,i}^2 \frac{\tr\big(\{\big(\BFI_{r_i+1}-\BFH_i\big)^\top \big(\BFI_{r_i+1}-\BFH_i\big)\}^2\big)}{\tr\big(\BFI_{r_i+1}-\BFH_i\big)^2}.
\]
Using Lemma~\ref{lem:traces} and the assumption on $\BFmu$ from Section~\ref{sec:parameter_models:mean}, we conclude that for $n,\BFr$, and $\lambda$ such that $n = o(\lambda)$ and $\lambda = o\big(\min_{i=1,\dots,q}r_i\big)^2$, since $\sum_{i=1}^qr_i\geqslant q \min_{i=1,\dots,q}r_i$,
\begin{align*}
\big|{\rm Bias}(\hat\sigma_i^2)\big| &\leqslant
O\left\{\lambda\Big(\min_{i=1,\dots,q}r_i\Big)^{-2}\right\} +
O\left\{\frac{(\lambda/n)^{-1/2}}{\min_{i=1,\dots,q}r_i}\right\} = o(1)
\quad\text{and}\quad\\
\mathbb{V}\hat\sigma_i^2 &\leqslant
O\left\{\lambda\Big(\min_{i=1,\dots,q}r_i\Big)^{-3}\right\} +
O\left\{\Big(\min_{i=1,\dots,q}r_i\Big)^{-1}\right\} = o(1),
\end{align*}
so that the risk of each of the estimators of $\sigma_{0,i}^2$ converges to zero.

In particular, we see that if $n=o\big(\min_{i=1,\dots,q}r_i\big)^2$, then choosing $\lambda$ as in~\eqref{eq:oracle_order} leads to the $\hat\sigma_i^2$ all being consistent.

\medskip

Noting that generalized cross validation provides us with a data-driven choice of $\lambda$ that is consistent for the minimizer of the risk, the statement of the theorem follows by an application of the law of total probability.

We will make this argument explicit for $\hat{\BFmu}_{\BFr}$. Let $\lambda_0$ denote the minimizer of the risk. Then, for arbitrary $\epsilon,\eta>0$ and $\lambda^*\in[\lambda_0(1-\eta), \lambda_0(1+\eta)]$,
\begin{align*}
\mathbb{P}\Big(\big(\sum_{i=1}^qr_i\big)^{-1}\big\|\hat{\BFmu}_{\BFr} - \BFmu_{0,\BFr}\big\|^2 > \epsilon \Big)
&\leqslant 2 \mathbb{P}\Big(\big(\sum_{i=1}^qr_i\big)^{-1}\left(\big\|\hat{\BFmu}_{\BFr}(\lambda^*) - \hat{\BFmu}_{\BFr}(\lambda_0)\big\|^2\right) > \frac{\epsilon}{2} \Big)\\
&\phantom{\leqslant}+2 \mathbb{P}\Big(\big(\sum_{i=1}^qr_i\big)^{-1}\left(\big\|\hat{\BFmu}_{\BFr}(\lambda_0) - \BFmu_{0,\BFr}\big\|^2 \right) > \frac{\epsilon}{2} \Big) + o(1),
\end{align*}
where, besides the law of total probability, we used the consistency of $\hat{\lambda}$ for $\lambda_0$ so that \[\mathbb{P}\Big(\frac{\hat{\lambda}}{\lambda_0}\notin [1-\eta, 1+\eta] \Big) = o(1),\:\:\:\mbox{and}\:\:\: \mathbb{P}\Big(\frac{\hat{\lambda}}{\lambda_0}\notin [1-\eta, 1+\eta] \Big)>\frac{1}{2}\] for large enough $n$, and using $\mathbb{P}(X+Y>\epsilon) \leqslant \mathbb{P}(X>\frac{\epsilon}{2}) + \mathbb{P}(Y>\frac{\epsilon}{2})$ for $X,Y\geqslant0$. Also note that $\hat\BFmu_\BFr(\cdot)$ abbreviates $\hat\BFmu_\BFr\big(\cdot, \BFSigma_\BFr^{(\BFn)}(\hat\BFtheta)\big)$.

We have already seen that $\big(\sum_{i=1}^qr_i\big)^{-1}\big\|\hat{\BFmu}_{\BFr}(\lambda_0) - \BFmu_{0,\BFr}\big\|^2=o_P(1)$. The conclusion follows by noting that the remaining term also goes to zero because the eigenvalues of the two matrices involved are close to one another for small $\eta$.
\endproof

\begin{corollary}
\label{cor:consistency_of_mu}
Assume the setting of Theorem \ref{theo:moments} and define $\hat{\BFmu} := \BFS_\BFr\hat{\BFmu}_{\BFr}$ for $\BFS_\BFr$ as defined in \eqref{eq:projection}. Then, as long as either $n\to\infty$ or $\min_{i=1,\dots,q}r_i\to\infty$,
$\hat\BFmu$ is consistent in probability for $\BFmu_{0} := \BFS_\BFr\BFmu_{0,\BFr}$.
\end{corollary}
\proof{}
From Theorem \ref{theo:moments}, it suffices to show that $\big\|\hat{\BFmu} - \BFmu_{0}\big\|^2$ can be upper bounded by an appropriate multiple of $\big\|\hat{\BFmu}_{\BFr} - \BFmu_{0,\BFr}\big\|^2$. Specifically, observe that
\[
\big\|\hat{\BFmu} - \BFmu_{0}\big\|^2 = \big\|\BFS_\BFr\hat{\BFmu}_{\BFr} - \BFS_\BFr\BFmu_{0,\BFr}\big\|^2 \leqslant \vertiii{\BFS_\BFr}^2 \big\|\hat{\BFmu}_{\BFr} - \BFmu_{0,\BFr}\big\|^2,
\]
where $\vertiii{\cdot}$ denotes the operator norm. Since we may write
\[
\BFS_\BFr^\top\BFS_\BFr = \diag\big\{\BFJ_{r_1}, \dots, \BFJ_{r_q}\big\},
\]
where $\BFJ_{i} \in\{1\}^{i\times i}$ denotes the $i \times i$ matrix of ones, we have that
\[
\vertiii{\BFS_\BFr}^2 = \max_{i=1,\dots,q}r_i.
\]
Combining this with \eqref{eq:risk_order}, we conclude that if the working variances $\sigma_i^2$ are bounded and $\lambda$ is picked as in~\eqref{eq:oracle_order}, then
$\big(\max_{i=1,\dots,q}r_i\sum_{i=1}^qr_i\big)^{-1}\mathbb{E}\big\|\hat{\BFmu} - \BFmu_{0}\big\|^2$
converges to zero as long as either $n$ or $\min_{i=1,\dots,q}r_i$ goes to infinity. The conclusion follows by an application of the law of total probability as can be seen in the proof of Theorem \ref{theo:moments}.
\endproof

\clearpage
\bibliographystyle{plain}

\end{document}